\documentclass[a4paper,reqno]{amsart}
\usepackage{amsmath, amsfonts, amsthm, amssymb}
\usepackage[maxbibnames=9]{biblatex}
\usepackage{mathrsfs}
\usepackage{enumerate}
\usepackage{xcolor}
\usepackage{booktabs}
\usepackage{placeins}
\usepackage{hyperref}
\usepackage{mathtools}
\usepackage[expansion=false]{microtype}

\DeclareFieldFormat{eprint:arxiv}{arXiv\addcolon\space\mbox{\ifhyperref{\href{https://arxiv.org/abs/#1}{\nolinkurl{#1}}}{\nolinkurl{#1}}}}
\numberwithin{equation}{section}

\allowdisplaybreaks
\hypersetup{
    bookmarks=true,         
    unicode=false,          
    pdftoolbar=true,        
    pdfmenubar=true,        
    pdffitwindow=false,     
    pdfstartview={FitH},    
    pdftitle={A general framework for universal bounds on the analytic functional calculus},
    pdfauthor={Mihai Putinar, Georgios Tsikalas},
    pdfkeywords={Crouzeix's conjecture, complete spectral set, numerical range, holomorphic functional calculus, H-infinity calculus, reproducing kernel, quantum annulus, Faber polynomial},
    pdfnewwindow=true,      
    colorlinks=false,       
    linkcolor=green,        
    citecolor=green,        
    filecolor=green,        
    urlcolor=green,         
    urlbordercolor={1 1 1}  
}

\newtheorem{theorem}{Theorem}[section]
\newtheorem{lemma}[theorem]{Lemma}
\newtheorem{proposition}[theorem]{Proposition}
\newtheorem{corollary}[theorem]{Corollary}
\newtheorem*{theoremA}{Theorem A}
\newtheorem*{theoremB}{Theorem B}
\newtheorem*{theoremC}{Theorem C}
\newtheorem*{theoremD}{Theorem D}
\theoremstyle{definition}
\newtheorem{definition}[theorem]{Definition}

\theoremstyle{remark}
\newtheorem{remark}[theorem]{Remark}

\newcommand{\C}{\mathbb{C}}
\newcommand{\Z}{\mathbb{Z}}
\newcommand{\RR}{\mathbb{R}}
\newcommand{\D}{\mathbb{D}}
\newcommand{\T}{\mathbb{T}}
\newcommand{\Ann}{\mathbb{A}}
\newcommand{\ip}[2]{\langle #1,\, #2\rangle}
\newcommand{\norm}[1]{\lVert #1\rVert}
\newcommand{\hs}[1]{\lVert #1\rVert_{\mathrm{HS}}}
\newcommand{\tr}{\operatorname{tr}}
\newcommand{\Rp}{\operatorname{Re}}
\newcommand{\Ip}{\operatorname{Im}}
\newcommand{\B}{\mathcal{B}}
\newcommand{\Sc}{\mathcal{S}_2}
\newcommand{\E}{\mathcal{H}}
\newcommand{\Kc}{\mathcal{K}}
\newcommand{\Bf}{\mathfrak{B}}
\newcommand{\Co}[1]{\mathcal{C}_{#1}}

\newcommand{\Om}{\Omega}
\newcommand{\Sbar}[1]{\overline{S_{#1}}}

\begin{document}

\title[\small Universal bounds on the analytic functional calculus]{\large A general framework for universal bounds on the analytic functional calculus}

\author[\small PUTINAR]{MIHAI PUTINAR}
\address{DEPARTMENT OF MATHEMATICS, UNIVERSITY OF CALIFORNIA,\newline SANTA BARBARA, CA}
\email{putinar@math.ucsb.edu}
\author[\small TSIKALAS]{GEORGIOS TSIKALAS}
\address{DEPARTMENT OF MATHEMATICS, UNIVERSITY OF CALIFORNIA,\newline SANTA BARBARA, CA}
\email{gtsikalas@ucsb.edu}

\subjclass[2020]{Primary 47A25; Secondary 47A60, 47A12, 46E22}
\keywords{Crouzeix's conjecture, complete spectral set, numerical range, holomorphic functional calculus, $H^\infty$-calculus, reproducing kernel, quantum annulus, Faber polynomial}
\small
\begin{abstract}
\small
We establish a general principle for bounding the matrix-valued analytic functional calculus $F\mapsto F(T)$ of a Hilbert space operator $T$, independently of the matrix size. It rests on two properties of an auxiliary Hilbert space of holomorphic functions: the accretivity of its moment map at $T$ and a uniform bound for its multipliers. The mechanism of OpenAI's direct proof of the complete Crouzeix conjecture is the instance in which these properties come from the positivity of the double-layer potential and from norm estimates for Faber series of a convex domain. The main advantage of our method is that the auxiliary space can be chosen freely, to fit the class of operators, and a well-chosen space can give sharp constants, also on domains that are not convex, simply connected or bounded. On every intersection of two disks of the Riemann sphere, for operators with both closed disks as spectral sets, we obtain the complete bound $2$, which is sharp for the quantum annulus and the strip. For operators, bounded or m-sectorial, with numerical range in a sector, our bounds improve on the complete Crouzeix inequality for every opening; they are optimal for the right angle, answering a question of Crouzeix, and nearly optimal for acute angles. On Jordan domains with regular real-analytic boundary, a space built from the exterior conformal map and the Faber polynomials of the domain gives complete bounds for operators whose double-layer potential is bounded below. They are optimal for the disk, recover the complete Crouzeix inequality on convex domains, and improve on earlier bounds, including bounds of Crouzeix and Greenbaum, also on nonconvex domains.
\end{abstract}
\maketitle

\par \normalsize

\section{Introduction}\label{sec:intro}

\subsection*{Spectral sets} Let $\Om\subseteq\C$ be a domain and let $\Kc$ be a class of bounded Hilbert space operators $T$ with $\sigma(T)\subset\Om$. The \emph{complete spectral constant}, or \emph{complete constant}, $K_{\mathrm{cb}}(\Kc)$ is the least $K\in[1,\infty]$ such that
\[
\norm{F(T)}\le K\sup_{z\in\Om}\norm{F(z)}
\]
for every $T\in\Kc$, every $d\ge1$ and every bounded holomorphic $F:\Om\to M_d$, where $M_d$ is the space of complex $d\times d$ matrices and $\norm{F(z)}$ is the operator norm. The \emph{scalar constant} $K(\Kc)$ is defined in the same way with $d=1$. Similarly, a compact subset $S\supseteq\sigma(T)$ of the Riemann sphere is a \emph{complete $K$-spectral set} for $T$ if $\norm{F(T)}\le K\max_S\norm F$ for every $d$ and every $M_d$-valued rational function $F$ with poles off $S$, and a \emph{$K$-spectral set} if this holds for $d=1$. For $K=1$ we speak of \emph{complete spectral sets} and \emph{spectral sets}. Throughout, the adjective \emph{complete} indicates that matrix-valued functions of every size are allowed, and not only scalar ones. Spectral sets were introduced by von Neumann, who proved that the closed unit disk is a spectral set for every contraction (von Neumann's inequality) \cite{vN}; this inequality is the starting point of the theory. By Sz.-Nagy's dilation theorem \cite{SN53}, the closed unit disk is even a complete spectral set for every contraction. By Paulsen's theorem \cite[p.~226]{Pa84b}, $S$ is a complete $K$-spectral set for $T$ if and only if $T$ is similar, by an invertible $G$ with $\norm G\norm{G^{-1}}\le K$, to an operator for which $S$ is a complete spectral set; for $S=\overline\D$, to a contraction. Passing from scalar bounds to complete ones is not automatic, and where it is possible it has usually required new arguments. The closed unit disk can be a $K$-spectral set for an operator for which it is not a complete $K'$-spectral set for any $K'$, by Pisier's counterexample to Halmos's similarity problem \cite{Pi97} and Paulsen's theorem. A spectral set that is a closed annulus is a complete spectral set, by a theorem of Agler \cite{Ag85}, but there are triply connected domains whose closures are spectral sets without being complete spectral sets \cite{DM05,AHR08}.

Bounds of the form $\norm{f(A)}\le K\sup_\Om|f|$, where $\overline\Om$ contains the numerical range of a matrix or operator $A$, are used to estimate the convergence rate of GMRES and of rational Krylov subspace methods for approximating $f(A)b$ \cite{CG19}. For sectors, they bound the powers $r(A)^n$ of rational functions $r$ with $|r|\le1$ on the sector uniformly in $n$, which yields stability criteria for time-discretization schemes \cite{Cr08}. They also bound the $H^\infty$-calculus of sectorial operators \cite{Ha06}, which applies to differential operators with complex coefficients (Section~\ref{sec:examples}).

A simple inspection of nilpotent matrices dispels any hope of bounding the analytic functional calculus in terms of the spectrum alone. A surprising benefit of recent advances on the optimal constant for the numerical range of an operator is that they adapt to a much wider context, which includes notable classes of unbounded operators. It is the aim of the present article to single out a unifying framework for bounding the complete constants $K_{\mathrm{cb}}(\Kc)$ directly, with no separate passage from scalar to matrix-valued functions. The classes we consider are defined by operator inequalities $\Rp r(T)\ge0$ for holomorphic functions $r$ on $\Om$, and they cover the classical hypotheses of the subject. For convex $\Om$, the supporting half-planes give the numerical-range hypothesis $W(T)\subseteq\overline\Om$ of Crouzeix's conjecture, and Theorem~D treats weaker hypotheses, lower bounds for the double-layer potential, also for nonconvex $\Om$. For intersections of disks $D_i$ of the Riemann sphere, among them annuli, strips and sectors, we use the conditions $\norm{b_i(T)}\le1$, with $b_i$ a conformal map of $D_i$ onto $\D$, which we call the \emph{disk norm conditions}; by von Neumann's inequality, they say that the closed disks $\overline{D_i}$ are spectral sets for $T$, and through the Cayley transform they take the form $\Rp r(T)\ge0$ (Section~\ref{ssec:test}). For the annulus $\{R^{-1}<|z|<R\}$, they define the quantum annulus $\{\norm T\le R,\ \norm{T^{-1}}\le R\}$, and for a sector they say that the numerical range lies in its closure, a hypothesis that extends to unbounded m-sectorial operators. These classes contain, for instance, Schr\"odinger operators with bounded complex potentials, after a shift, divergence form operators with complex coefficients and hyponormal operators (Section~\ref{sec:examples}).

\subsection*{The moment map theorem}
Crouzeix conjectured that the numerical range is a $2$-spectral set, and proved this for $2\times2$ matrices \cite{Cr04}; the complete bound $2$ was then proved for $2\times2$ matrices \cite{BCD06} and for weighted shift matrices \cite{Ch13} (see also \cite{CG25}). The scalar conjecture was proved by Jin \cite{Jin} and by Lorist and Schwenninger \cite{LS} (see also \cite{Luo}), and the complete one for matrices of order at most $3$ by Åhag, Czyż and Virtanen \cite[Theorem~1.1]{ACV}. OpenAI proved the complete Crouzeix conjecture in two manuscripts with separate proofs \cite{OAI,OAI2}. The direct proof \cite[Theorem~1.1]{OAI} gives the complete bound $2$ for matrices $A$ and matrix polynomials $F$. The similarity proof \cite[Theorem~1.2]{OAI2} constructs, for a matrix $A$, an invertible $S$ with $\norm S\norm{S^{-1}}\le2$ such that $Sf(A)S^{-1}$ is a contraction, where $f$ maps a bounded convex domain with regular real-analytic boundary containing $W(A)$ conformally onto $\D$. The complete bound $2$ extends to bounded operators by \cite[Corollary~6.3]{OAI2}. Both proofs combine the Faber polynomials of the domain with the positivity of the double-layer potential, which already underlies the earlier bounds of \cite{DD} and the bound $1+\sqrt2$ of \cite{CP17,RS18} (see also \cite{Cr07}) and the refinement of the latter for convex domains in \cite[Theorem~3]{MMOR25}. The complete bound $1+\sqrt2$ was stated in \cite{CP17} and proved for matrices in \cite[Section~2]{CK25}. The arguments of \cite{CP17} and \cite{LS} are extended to other sets in \cite{CG19} and \cite{CG26}. An abstract form of the argument of \cite{CP17}, for homomorphisms of uniform algebras, is developed in \cite{COR23,SdV24}, and the abstract version of Crouzeix's conjecture formulated in \cite{COR23} (see also \cite[Question~4.1]{RS18}), which implies the scalar conjecture, is proved in \cite{BOV26}.

The direct proof \cite[Section~3]{OAI} works on a bounded convex domain with regular real-analytic boundary containing $W(A)$, with unit vectors $x,y$ such that $F(A)x=\norm{F(A)}\,y$ and $F(A)^*y=\norm{F(A)}\,x$. It compares, in two steps, the norms of the linear functionals $G\mapsto\ip{x}{G(A)x}$, $G\mapsto\ip{y}{G(A)y}$ and $G\mapsto\ip{y}{G(A)x}$ on matrix polynomials $G=\sum_n\Phi_nC_n$ $(C_n\in M_d)$, with $\Phi_n$ the Faber polynomials of the domain (recalled before Theorem~D), normed by weighted $\ell^2$ norms of the coefficients $C_n$ (Remark~\ref{rem:engine}). In the first step, the norm estimates for Faber series of \cite[Proposition~2.2]{OAI}, applied to $FG$ and $GF$, bound $\norm{F(A)}$ times the norms of the first two functionals by the norm of the third when $\norm F\le1$ on the boundary. In the second, the positivity of the double-layer potential bounds the norm of the third functional by those of the first two. Together, the two steps give $\norm{F(A)}\le2$.

We raise the direct proof of \cite{OAI} to a general, unifying principle, as follows. Consider a Hilbert space $\E$ of holomorphic functions on $\Om$ with continuous point evaluations, orthonormal basis $(e_j)$ and reproducing kernel $k(z,w)=\sum_je_j(z)\overline{e_j(w)}$, in place of the Faber series. Let $\Sc(H)$ be the class of Hilbert--Schmidt operators on a Hilbert space $H$, with $\ip XY_{\mathrm{HS}}=\tr(X^*Y)$. For $T\in\B(H)$ with $\sigma(T)\subset\Om$ and $u,v\in H$, the functional $Q\mapsto\ip{u}{Q(T)v}$ on $\E$ has squared norm $\sum_j|\ip{u}{e_j(T)v}|^2=\ip{uu^*}{\Psi_T(vv^*)}_{\mathrm{HS}}$, where
\[
\Psi_T(X)=\sum_je_j(T)\,X\,e_j(T)^*\qquad(X\in\Sc(H)).
\]
We call $\Psi_T$ the \emph{moment map} of $\E$ at $T$; it depends only on $k$ (Lemma~\ref{lem:cauchy}), and it is \emph{accretive} if
\[
\Rp\ip X{\Psi_T(X)}_{\mathrm{HS}}\ge0\qquad(X\in\Sc(H)).
\]

The proof of Theorem~A below follows the two-step scheme of the direct proof of \cite{OAI}; we describe it here for holomorphic $F:\Om\to\C$ with $\sup_\Om|F|\le1$ for which $\norm{F(T)}$ is attained (the proof itself does not need this), with unit vectors $x,y\in H$ such that $F(T)x=\norm{F(T)}\,y$ and $F(T)^*y=\norm{F(T)}\,x$, and with the functionals $Q\mapsto\ip{x}{Q(T)x}$, $Q\mapsto\ip{y}{Q(T)y}$ and $Q\mapsto\ip{y}{Q(T)x}$ on $\E$. In the first step, the norm estimates for Faber series are replaced by the multiplier constant of $\E$,
\begin{multline*}
M(\E)=\sup\big\{\norm{Ff}:\ d\ge1,\ f\in\E\otimes\C^d,\ \norm f\le1,\\* F:\Om\to M_d\text{ holomorphic},\ \sup_\Om\norm F\le1\big\},
\end{multline*}
where $\E\otimes\C^d$ is the space of $\C^d$-valued functions with entries in $\E$ (Definition~\ref{def:mult}, Remark~\ref{rem:usualmult}). In the second step, the bound of the norm of the third functional by those of the first two follows, in terms of $\Psi_T$, from the Cauchy--Schwarz inequality for the Hermitian form $\ip{X}{(\Rp\Psi_T)(Y)}_{\mathrm{HS}}$, which is positive semidefinite exactly when $\Psi_T$ is accretive (Lemma~\ref{lem:momentineq}). In \cite{OAI}, this bound comes from a trace inequality for the double-layer potential, integrated over the boundary; for the kernel $2+\sum_{n\ge1}\Phi_n(z)\overline{\Phi_n(w)}$, that integrated inequality is a special case of the accretivity of $\Psi_T$ (Remark~\ref{rem:convex}). The other differences between the two proofs are in Remark~\ref{rem:engine}. The argument also works with two spaces, one for the third functional and one for the first two, as with the two norms of \cite{OAI}; Remark~\ref{rem:twospaces} explains why we use one. For $M_d$-valued $F$, one takes $x,y\in\C^d\otimes H$ and the functionals on $\E\otimes M_d$, with the Hilbert--Schmidt norm on $M_d$; their squared norms are still given by $\Psi_T$, with $uu^*$ replaced by $\sum_ru_ru_r^*$ for $u=(u_1,\dots,u_d)$ (Lemmas~\ref{lem:coef} and~\ref{lem:momentineq}). So the second step needs only the accretivity of $\Psi_T$ on $\Sc(H)$, for every $d$, and Theorem~A gives complete bounds directly (Remark~\ref{rem:complete}).

Call $\E$ \emph{nondegenerate} if it contains functions $q_n\to1$ locally uniformly on $\Om$.

\begin{theoremA}[Theorem~\ref{thm:engine}, Corollary~\ref{cor:main}]
Let $\E$ be nondegenerate and let $\Kc$ be a class of operators with spectrum in $\Om$ such that $\Psi_T$ is accretive for every $T\in\Kc$. Then
\[
\norm{F(T)}\le\sqrt2\,M(\E)\sup_\Om\norm F
\]
for every $T\in\Kc$, every $d\ge1$ and every bounded holomorphic $F:\Om\to M_d$.
\end{theoremA}

In fact, we obtain something stronger. Nondegeneracy is used only to provide, for each $T$ with $\sigma(T)\subset\Om$, some $q\in\E$ with $\Rp q(T)\ge\delta I$ for some $\delta>0$ (see after Definition~\ref{def:nondeg}). Theorem~\ref{thm:engine} assumes only this, together with the accretivity of $\Psi_T$, and gives $\norm{F(T)}^2\le2\,m_\ell(F)\,m_r(F)$, where $m_\ell(F)$ and $m_r(F)$ are the norms of $Q\mapsto FQ$ and $Q\mapsto QF$ on $\E\otimes M_d$, both at most $M(\E)\sup_\Om\norm F$ (Section~\ref{ssec:kernels}, Remark~\ref{rem:usualmult}).

To apply Theorem~A to a class $\Kc$, one needs a nondegenerate space $\E$ whose moment maps are accretive at every $T\in\Kc$, and a bound for $M(\E)$. For the complete Crouzeix inequality on a convex domain $\Om$ with regular real-analytic boundary, the function theory of $\Om$ enters only through the kernel $2+\sum_{n\ge1}\Phi_n(z)\overline{\Phi_n(w)}$ and two facts about it: its moment maps are accretive when $W(T)\subseteq\overline\Om$, by the positivity of the double-layer potential, and the multiplier constant of its space is at most $\sqrt2$, by the norm estimates for Faber series of \cite[Proposition~2.2]{OAI}; Theorem~A then gives the constant $2$ (Theorem~D, Remark~\ref{rem:convex}).

\subsection*{Intersections of two disks}
For holomorphic $r:\Om\to\{\Rp z>0\}$ with $\Rp r(T)\ge0$, the moment map of the kernel $k_r=1/(r(z)+\overline{r(w)})$ is the inverse of the accretive map $X\mapsto r(T)X+Xr(T)^*$, hence accretive, and so is that of $k_{1/r}$, since $\Rp\big(r(T)^{-1}\big)\ge0$ (Section~\ref{sec:positivity}). As the moment map depends linearly on the kernel (Lemma~\ref{lem:cauchy}), sums of such kernels have accretive moment maps throughout a class defined by several such inequalities, for instance by the disk norm conditions $\norm{b_i(T)}\le1$, with $r_i=(1+b_i)/(1-b_i)$ (Corollary~\ref{cor:sums}).

\begin{theoremB}[Lemma~\ref{lem:mobius}, Propositions~\ref{prop:annulus}, \ref{prop:strip} and~\ref{prop:sectorbound}]
Let $D_1,D_2$ be open disks of the Riemann sphere such that $\Om:=D_1\cap D_2$ is a nonempty subset of $\C$ different from $D_1$ and $D_2$, and let $b_i:D_i\to\D$ be conformal. For every $T\in\B(H)$ with $\sigma(T)\subset\Om$ and $\norm{b_i(T)}\le1$ ($i=1,2$), every $d\ge1$ and every bounded holomorphic $F:\Om\to M_d$,
\[
\norm{F(T)}\le2\sup_\Om\norm F .
\]
\end{theoremB}

A M\"obius map reduces $\Om$ to an annulus $\Ann_R=\{R^{-1}<|z|<R\}$, a strip $\Pi_a=\{|\Ip z|<a\}$ or a sector $S_\alpha=\{z\ne0:|\arg z|<\alpha\}$ of opening $2\alpha$, with $0<\alpha<\pi/2$ (disjoint, tangent or crossing circles), and the disk norm conditions to $\norm T,\norm{T^{-1}}\le R$, to $\norm{\Ip T}\le a$, and to $W(T)\subseteq\Sbar\alpha$, the closed sector. For the annulus and the strip, $\E$ is the reproducing kernel Hilbert space of a sum of two kernels, one for each disk, chosen to match the symmetry of the domain: $\frac12(k_r+k_{1/r})$ for each of the two disks of the annulus, and $k_r$ for each of the two half-planes $\{\Rp r>0\}$, with $r$ affine, of the strip. Explicitly, the two kernels are
\begin{equation}\label{eq:twokernels}
\sum_{n\in\Z}R^{-2|n|}(z\bar w)^n\qquad\text{and}\qquad\frac{4a}{4a^2+(z-\bar w)^2}.
\end{equation}
For the sector, Theorem~C below uses instead, for openings $2\alpha\le\pi/2$, the space of a kernel built from the product of the kernels $k_r$ of the two half-planes, and, for $2\alpha\ge\pi/2$, a Hardy space of the sector; these spaces give smaller constants, so Theorem~B follows from Theorem~C.

The choice of space is essential. A Hardy space of $\Om$, with the $L^2$ norm on the two boundary curves normalized as in Sections~\ref{ssec:annulus}--\ref{ssec:sector}, has multiplier constant $1$; but for the quantum annulus, the strip and sectors of opening below $\pi/2$, no nondegenerate space with $M(\E)=1$ has accretive moment maps throughout the class, since the scalar constants exceed $\sqrt2$ (for sectors, they are at least $\pi\sin\alpha/(2\alpha)$ \cite{CD03}, \cite[Theorem~1]{BC06}, \cite[Lemma~4.1]{Cr16}). The space $\E$ we use instead has the same functions as this Hardy space but a different norm (Sections~\ref{ssec:annulus} and~\ref{ssec:strip}). For the annulus, the squared norms of $f(z)=\sum_nc_nz^n$ in $\E$ and in the Hardy space are
\[
\sum_nR^{2|n|}|c_n|^2\qquad\text{and}\qquad\sum_n(R^{2n}+R^{-2n})|c_n|^2,
\]
and for the strip, with $f(z)=\int\hat f(t)e^{itz}\,dt$, they are
\[
\int|\hat f(t)|^2e^{2a|t|}\,dt\qquad\text{and}\qquad\int|\hat f(t)|^2(e^{2at}+e^{-2at})\,dt.
\]
In both cases, the weight of $\E$ is between $\frac12$ and $1$ times that of the Hardy space, so $M(\E)\le\sqrt2$ (Lemma~\ref{lem:compare}), which gives the constant $2$. For the annulus, $\E$ appears, up to normalization, in \cite{ARS08}, and, up to the rescaling $z\mapsto z/R$, in \cite[Proposition~3.1]{Ts22}, where this bound is proved for scalar multipliers and shown to be sharp.

By Theorem~B and Corollary~\ref{cor:qa}, for every $R>1$, the closed annulus $\overline{\Ann_R}$ is a complete $2$-spectral set for every $T$ with $\norm T\le R$ and $\norm{T^{-1}}\le R$; the constant $2$ is optimal, since the scalar constant already equals $2$, by the lower bound of \cite{Ts}. So the complete constant of the quantum annulus is $2$ for every $R>1$; the earlier complete bounds, recalled below, all exceed $2$. For the strip, the scalar and the complete constant both equal $2$ (Proposition~\ref{prop:stripopt}); the lower bound comes from a continuous analogue of the weighted shifts used for the annulus in \cite{Ts}.

The question of the scalar constant of the quantum annulus goes back to Shields \cite{Sh}. Badea, Beckermann and Crouzeix proved that every intersection of two disks is a complete $(2+2/\sqrt3)$-spectral set under the disk norm conditions \cite{BBC07,BBC}. For the quantum annulus, further scalar bounds are in \cite{CG19,Pa25} (with two further proofs of the latter in \cite{PPT26}) and, for matrices, in \cite{CrQA}, and complete bounds larger than $2$ are in \cite{DP86,Cr12,JT24,To26}. The scalar bounds of \cite{CrQA,Pa25} and the complete bound of \cite{To26} tend to $2$ as $R\to\infty$. The scalar bound $2$ was proved by Jin \cite[Section~7.3]{Jin} and by Lorist and Schwenninger \cite[Section~3(v)]{LS}, and Jin \cite[Section~7.4]{Jin} obtained the scalar bound $2$ for every intersection of two disks (see also \cite[Section~4.1 and Theorem~3]{CG26}). See \cite{MP23,PT25} for dilations and models and \cite{BY23,Ts22,JT24} for related classes.

For the strip, where $\norm{\Ip T}\le a$ means that $W(T)$ lies in the closed strip, Crouzeix and Delyon proved the bound $2+2/\sqrt3$ \cite{CD03} (see also \cite[Theorem~7.1.16]{Ha06} and \cite[Section~3]{Cr16}), and its complete form follows from \cite[Section~3, (3.8)--(3.10)]{LeM03}. Crouzeix regarded the strip as a natural candidate for the largest constant in his conjecture and asked whether its complete constant is at most $2$ \cite[Sections~7 and~11]{Cr16}. It is, by \cite{OAI,OAI2}, since the strip is convex, and Theorem~B gives a second proof; this case of Theorem~B also follows from the annulus case with $R=e^{a/n}$, applied to the operators $e^{iT/n}$ and to $2\pi n$-periodic approximations of $F$, as $n\to\infty$ (Remark~\ref{rem:stripannulus}).

\subsection*{Sectors, lenses and crescents}
For sectors, the complete bound $2$ also follows from \cite{OAI,OAI2}, since the closed sector is convex. The following theorem lowers this bound for every opening.

\begin{theoremC}[Propositions~\ref{prop:sectorbound} and~\ref{prop:sectorlower}, Corollary~\ref{cor:sectorclosed}]
Let $0<\alpha<\pi/2$ and $c_\alpha=\max\{\sqrt2,\sqrt{4-8\alpha/\pi}\}$. For every bounded operator $A$ with $W(A)\subseteq\Sbar\alpha=\{z:|\arg z|\le\alpha\}\cup\{0\}$, every $d\ge1$ and every $M_d$-valued rational function $F$ with poles off $\Sbar\alpha$,
\[
\norm{F(A)}\le c_\alpha\sup_{\Sbar\alpha}\norm F .
\]
For $\alpha=\pi/4$ the constant $\sqrt2$ is optimal. For $\alpha<\pi/4$ the optimal constant is at least $2^{1-2\alpha/\pi}$.
\end{theoremC}

For $\alpha\le\pi/4$, $c_\alpha=\sqrt{4-8\alpha/\pi}$ is less than $1.031\cdot2^{1-2\alpha/\pi}$. The bound and the lower bounds extend to m-sectorial $A$, that is, closed, possibly unbounded, $A$ with $W(A)\subseteq\Sbar\alpha$ whose resolvent set contains $\C\setminus\Sbar\alpha$, for every $M_d$-valued rational $F$ bounded on $\Sbar\alpha$ (Corollary~\ref{cor:unbounded}; the lower bounds already hold for matrices), and, by a M\"obius transformation, to lenses and crescents $\Om$ of vertex angle $2\alpha$ bounded by two circles: $\norm{F(T)}\le c_\alpha\sup_\Om\norm F$ for $T\in\B(H)$ satisfying the disk norm conditions with $\sigma(T)\subset\Om$ and bounded holomorphic $F:\Om\to M_d$ (Corollary~\ref{cor:lens}).

For $\alpha\le\pi/4$, we use, instead of the sum of the kernels $k_r$ of the two half-planes, their product times $2\sin(2\alpha)z\bar w$, the kernel
\[
\frac{2\sin(2\alpha)\,z\bar w}{z^2-2\cos(2\alpha)\,z\bar w+\bar w^2}.
\]
With $f(z)=\int\hat f(\xi)z^{i\xi}\,d\xi$, the squared norms of $f$ in the space $\E$ of this kernel and in the Hardy space $H^2(S_\alpha,dr/r)$, defined with the measure $dr/r$ on the two rays, are
\[
\int|\hat f(\xi)|^2\,\frac{\sinh(\pi\xi)}{\sinh((\pi-2\alpha)\xi)}\,d\xi\qquad\text{and}\qquad\int|\hat f(\xi)|^2\,2\cosh(2\alpha\xi)\,d\xi;
\]
the ratio of the second weight to the first lies between $1$ and $2-4\alpha/\pi$, so $M(\E)\le\sqrt{2-4\alpha/\pi}$. For bounded $A$ with $\sigma(A)\subset S_\alpha$ and $W(A)\subseteq\Sbar\alpha$, the moment map of this space at $A$ has an explicit inverse, built from $X\mapsto AX(A^{-1})^*$, and its accretivity follows from a sum-of-squares identity and an inequality for operators with numerical range in a sector (Lemma~\ref{lem:inverse}, Proposition~\ref{prop:rightangle}). For $\pi/4\le\alpha<\pi/2$, we use $H^2(S_\alpha,dr/r)$ itself, with $M(\E)=1$; the change of variable $z\mapsto z^{\pi/(4\alpha)}$ and Kato's theorem on fractional powers \cite{Ka} transfer accretivity from the right-angle sector (Lemma~\ref{lem:fracpow}, Proposition~\ref{prop:obtuse}).

By Lemma~\ref{lem:mobius} (for convex lenses, see \cite[Theorem~1]{BC06}), the best constant for lenses and crescents of vertex angle $2\alpha$ under the disk norm conditions equals that for $S_\alpha$ under $W(A)\subseteq\Sbar\alpha$. Earlier upper bounds for this constant are in \cite{CD03,BCD06}, \cite[Theorem~5 and Corollary~6]{BC07} and \cite[Section~3]{Cr16}, and Crouzeix \cite{Cr25} improved them to
\begin{equation}\label{eq:crouzeixbound}
1-\frac\alpha\pi+\sqrt{2-\frac{4\alpha}\pi+\frac{\alpha^2}{\pi^2}}.
\end{equation}
The bound \eqref{eq:crouzeixbound} decreases in $\alpha$, equals $\sqrt2$ at $\alpha=\alpha_*=\pi(5-3\sqrt2)/2$, and, as he noted, also holds for the complete constant.

For the right angle, the constant for $2\times2$ matrices $A$ with $W(A)\subseteq\Sbar{\pi/4}$ and rational $f$ was known to be $\sqrt2$ \cite[Theorem~4.2]{Cr04}, \cite[Sections~3 and~4]{Cr16}. Crouzeix reported numerical evidence for $4\times4$ matrices and asked whether $\sqrt2$ holds for matrices of every size \cite[Sections~8 and~11]{Cr16}, and Zeng proved it for the operators $\lambda I+N$ with $N^2=0$ \cite[Theorem~1.1]{Ze26} and for special families of rational functions, and stated that the general case remains open. Jin \cite[Section~7.6]{Jin} lists the constant of a fixed lens among open questions. Theorem~C and its extension to lenses and crescents (Corollary~\ref{cor:lens}) answer Crouzeix's question, and Jin's in the right-angle case: for the right-angle sector, and for lenses and crescents bounded by orthogonal circles, the scalar and the complete constant are both $\sqrt{2}$.

For $0<\alpha<\pi/4$, Theorem~C gives the complete bound $\sqrt{4-8\alpha/\pi}$, which is below $2$ and below \eqref{eq:crouzeixbound}, and no other choice of the space in Theorem~A gives less (Proposition~\ref{prop:spaceopt}). Theorem~C also gives the lower bound $2^{1-2\alpha/\pi}$, which exceeds the lower bound $\pi\sin\alpha/(2\alpha)$ of \cite{CD03}, \cite[Theorem~1]{BC06} and \cite[Lemma~4.1]{Cr16}, and, numerically, the constant for $2\times2$ matrices determined in \cite[Theorem~4.2]{Cr04} (Remark~\ref{rem:sectorcompare}). Thus the scalar and the complete constant both lie between $2^{1-2\alpha/\pi}$ and $\sqrt{4-8\alpha/\pi}$, and both tend to $2$, the constant of the strip, as $\alpha\to0$. For $\pi/4<\alpha<\pi/2$, Theorem~C gives the complete bound $\sqrt2$, and our method cannot give less, since Theorem~A never yields a constant below $\sqrt2$ (Remark~\ref{rem:floor}). For these $\alpha$, the scalar and the complete constant are bounded below by the constant for $2\times2$ matrices of \cite[Theorem~4.2]{Cr04}, which is at least $\pi\sin\alpha/(2\alpha)$ (Remark~\ref{rem:sectorcompare}). The bound $\sqrt2$ improves on \eqref{eq:crouzeixbound} exactly when $\alpha<\alpha_*$.

\subsection*{Domains with real-analytic boundary}
Theorem~D needs more of the classical function theory of the domain than Theorems~B and~C: its exterior conformal map, its Faber polynomials and the double-layer potential of its boundary. Let $\Om$ be a bounded Jordan domain with regular real-analytic boundary (Section~\ref{sec:convex}), and let $\psi$ be the exterior conformal map of $\Om$, from $\{|t|>1\}$ onto $\C\setminus\overline\Om$, with $\psi(t)/t\to c>0$ as $t\to\infty$. For $t\in\T$ and $z\in\Om$ let
\[
r_t(z)=\frac{t\psi'(t)}{\psi(t)-z}=1+\sum_{n\ge1}\Phi_n(z)\,t^{-n};
\]
the coefficients $\Phi_n$ are the Faber polynomials of $\Om$, and for the disk, $\psi(t)=t$ and $\Phi_n(z)=z^n$. If $\Om$ is convex, $r_t$ has positive real part on $\Om$, and for $\sigma(T)\subset\Om$, $\Rp r_t(T)\ge0$ for all $t$ exactly when $W(T)\subseteq\overline\Om$, since $1/r_t$ is affine and $\{\Rp(1/r_t)>0\}$ is the open half-plane bounded by the tangent line at $\psi(t)$ that contains $\Om$. Instead of a finite sum of kernels $k_r$, we use for this family the kernel
\[
1+\int_\T r_t(z)\overline{r_t(w)}\,d\omega(t)=2+\sum_{n\ge1}\Phi_n(z)\overline{\Phi_n(w)},
\]
with $\omega$ the normalized arclength measure: since $r_t$ is holomorphic in $t$ near $\{|t|\ge1\}$ with $r_t\to1$ as $t\to\infty$ (Lemma~\ref{lem:convexclass}), its moment map $\Psi_T$ satisfies (Proposition~\ref{prop:pencil})
\[
\Rp\ip{X}{\Psi_T(X)}_{\mathrm{HS}}=\frac12\int_\T\ip{X}{N_T(t)\,X\,N_T(t)}_{\mathrm{HS}}\,d\omega(t),\qquad N_T(t)=r_t(T)+r_t(T)^*,
\]
which is nonnegative when every $N_T(t)$ is positive; $N_T(t)$ is a positive multiple of the operator-valued double-layer potential of \cite{DD,CP17} (see \eqref{eq:dlp}). The multiplier constant is bounded by the norm estimates for Faber series of \cite{OAI}, which extend, with the same proof, to nonconvex $\Om$ (Lemma~\ref{lem:faber}).

For convex $\Om$, both steps of the direct proof of \cite{OAI} rest on the positivity of the double-layer potential, in two forms: the first, through the norm estimates for Faber series, on the nonnegativity of its kernel $k_\partial(t,u)=2\Rp r_t(\psi(u))$ $(t,u\in\T,\ t\ne u)$ on $\partial\Om$, in the parametrization $\zeta=\psi(t)$, and the second on the positivity of the operator-valued potential $N_T(t)=2\Rp r_t(T)$, which follows from $W(T)\subseteq\overline\Om$. For nonconvex $\Om$ and some $t\in\T$, $\Rp r_t$ takes negative values on $\Om$ and on $\partial\Om$ (Lemma~\ref{lem:rhoOm}(a)), so neither form of positivity is available. Theorem~D replaces the second by the lower bound $\Rp r_t(T)\ge(1-\frac\rho2)I$ for some $\rho\ge\sqrt2$, and the first by the norm $V(\Om)$ of the integral operator with kernel $k_\partial$ on the closed span of the functions $t^n$, $n\ge1$, in $L^2(\T)$. In terms of the Grunsky coefficients $b_{nj}$ of $\psi$, $V(\Om)$ is the norm on $\ell^2$ of the matrix $(nb_{nj})_{n,j\ge1}$, and it is at most $1$ for convex $\Om$ (Lemma~\ref{lem:rhoOm}(c)); but convexity is not necessary: there is a sequence of nonconvex domains, all contained in a fixed disk, with $V(\Om)\le1$ for each of them while the minima of their kernels $k_\partial$ tend to $-\infty$ (Lemma~\ref{lem:rhoOm}(d) and the example after it).

\begin{theoremD}[Theorem~\ref{thm:convex}]
Let $\Om$ be a bounded Jordan domain with regular real-analytic boundary, and let $\rho\ge\sqrt2$. For every $T\in\B(H)$ with $\sigma(T)\subset\Om$ and
\[
\Rp r_t(T)\ge\Big(1-\frac\rho2\Big)I\qquad(t\in\T),
\]
every $d\ge1$ and every bounded holomorphic $F:\Om\to M_d$,
\[
\norm{F(T)}\le\max\big\{\rho,\sqrt{2+2V(\Om)^2}\big\}\sup_\Om\norm F .
\]
For $\rho\ge2$, the constant is $\rho$ if $V(\Om)\le1$, in particular if $\Om$ is convex; for $\Om=\D$, it is $\rho$ for every $\rho\ge\sqrt2$.
\end{theoremD}

For $\rho=2$, the hypothesis means that $W(T)$ lies in the intersection of the closed half-planes $\{\Rp(1/r_t)\ge0\}$, $t\in\T$, a convex subset of $\overline\Om$ that equals $\overline\Om$ for convex $\Om$. For convex $\Om$, Theorem~D with $\rho=2$ and Corollary~\ref{cor:convexclosed}(b) therefore give the complete Crouzeix inequality of \cite{OAI,OAI2} for these domains: $\overline\Om$ is a complete $2$-spectral set for every $T$ with $W(T)\subseteq\overline\Om$. For every $\Om$ and every $\rho\in[\sqrt2,2]$, the complete Crouzeix inequality, applied to $\overline{W(T)}\subseteq\overline\Om$, gives the constant $2$ under the hypothesis of Theorem~D (Remark~\ref{rem:rho2}), and Theorem~D improves on this constant exactly when $\rho<2$ and $V(\Om)<1$.

For the disk, where $r_t(T)=(I-\bar tT)^{-1}$, the hypothesis characterizes, for every $\rho>0$ and among the operators with $\sigma(T)\subset\D$, the \emph{$\rho$-contractions}, the operators with $T^n=\rho P_HU^n|_H$ $(n\ge1)$ for a unitary $U$ on a Hilbert space containing $H$, where $P_H$ is the orthogonal projection onto $H$ \cite{SzNF}. By analogy, we say that $T$ has a \emph{Faber $\rho$-dilation} if
\[
\Phi_n(T)=\rho P_HU^n|_H\qquad(n\ge1)
\]
for such a $U$; by Lemma~\ref{lem:rhoconvex}, an operator $T$ with $\sigma(T)\subset\Om$ has one exactly when it satisfies the hypothesis of Theorem~D. For $\rho>2$, if $\Om$ is convex or $\rho$ is large enough, the hypothesis means $\norm{(\lambda-T)^{-1}}\le\mathrm{dist}(\lambda,\overline\Om)^{-1}$ for every $\lambda$ on the curve $\{\psi(t)+t\psi'(t)/(\rho-2):t\in\T\}$ (Lemma~\ref{lem:rhoconvex}(b)). For example, hyponormal operators with spectrum in $\Om$ have Faber $\rho$-dilations for every $\rho\ge\rho_\Om$, a constant determined by the double-layer potential of $\partial\Om$, with $\rho_\Om=2$ exactly for convex $\Om$ (Lemma~\ref{lem:rhoOm}(a), Proposition~\ref{prop:hyponormal}). Special cases of these characterizations are known (Remark~\ref{rem:special}). Without assuming $\sigma(T)\subset\Om$, a Faber $\rho$-dilation of $T$ forces $\sigma(T)\subseteq\overline\Om$ and makes $\overline\Om$ a complete $K$-spectral set for $T$, with the constant $K$ of Theorem~D (Corollary~\ref{cor:convexclosed}(a)); for the disk, this is the theorem of Okubo and Ando \cite{OA} for $\rho\ge\sqrt2$ (for other proofs and refinements, see \cite{Dr08,SdV,HdV,HM} and Remark~\ref{rem:rho}).

\subsection*{What is new}
In the direct proof of \cite{OAI}, both ingredients of Theorem~A are dictated by the convex domain: the space is that of its Faber series, and the accretivity of the moment map comes from the positivity of its double-layer potential (Remarks~\ref{rem:engine} and~\ref{rem:convex}). Theorem~A, in contrast, applies to every nondegenerate space $\E$: the space can be chosen to fit the class, which enters only through the accretivity of the moment map, a single operator inequality. Bounding a complete constant thus becomes the construction of a nondegenerate space with accretive moment maps and a small multiplier constant, when one exists. The spaces we use are the following, with the results that give their accretivity; in each row, the bound is $\sqrt2$ times the bound for $M(\E)$, by Theorem~A.
\begin{center}
\small
\setlength{\tabcolsep}{2.75pt}
\begin{tabular}{@{}lllll@{}}
\toprule
Class & Kernel & Accretivity & $M(\E)$ & Bound\\
\midrule\addlinespace[3pt]
$\norm T\le R$, $\norm{T^{-1}}\le R$ & $\sum_{n\in\Z}R^{-2|n|}(z\bar w)^n$ & Cor.~\ref{cor:sums} & $\sqrt2$ & $2$ (sharp)\\\addlinespace[8pt]
$\norm{\Ip T}\le a$ & $\dfrac{4a}{4a^2+(z-\bar w)^2}$ & Cor.~\ref{cor:sums} & $\sqrt2$ & $2$ (sharp)\\\addlinespace[8pt]
$W(T)\subseteq\Sbar\alpha$, $0<\alpha<\frac\pi4$ & $\dfrac{2\sin(2\alpha)\,z\bar w}{z^2-2\cos(2\alpha)\,z\bar w+\bar w^2}$ & Prop.~\ref{prop:rightangle} & $\sqrt{2-\frac{4\alpha}\pi}$ & $\sqrt{4-\frac{8\alpha}\pi}$\\\addlinespace[8pt]
$W(T)\subseteq\Sbar{\pi/4}$ & $\dfrac{2z\bar w}{z^2+\bar w^2}$ & Prop.~\ref{prop:rightangle} & $1$ & $\sqrt2$ (sharp)\\\addlinespace[8pt]
$W(T)\subseteq\Sbar\alpha$, $\frac\pi4<\alpha<\frac\pi2$ & $\dfrac{2s\,z^s\bar w^s}{z^{2s}+\bar w^{2s}}$, $s=\pi/(4\alpha)$ & Prop.~\ref{prop:obtuse} & $1$ & $\sqrt2$\\\addlinespace[8pt]
Faber $\rho$-dilation, $\rho\ge\sqrt2$ & $2+\frac4{\rho^2}\sum_{n\ge1}\Phi_n(z)\overline{\Phi_n(w)}$ & Prop.~\ref{prop:pencil} & $\le C_{\Om,\rho}/\sqrt2$ & $C_{\Om,\rho}$\\
\bottomrule
\end{tabular}
\end{center}
\medskip

Here $\sigma(T)\subset\Om$, where $\Om$ is $\Ann_R$, $\Pi_a$ or $S_\alpha$ in the first five rows and, in the last, a bounded Jordan domain with regular real-analytic boundary, with Faber polynomials $\Phi_n$ and $C_{\Om,\rho}=\max\{\rho,\sqrt{2+2V(\Om)^2}\}$ the constant of Theorem~D. For spectra in $\overline\Om$, see Corollaries~\ref{cor:qa}, \ref{cor:sectorclosed} and~\ref{cor:convexclosed}.

To our knowledge, the following are new:
\begin{enumerate}[(i)]
\item Theorem~A for every nondegenerate reproducing kernel Hilbert space of holomorphic functions;
\item for two disks with disjoint boundary circles, the complete bound $2$, which is optimal for the quantum annulus, for every $R>1$, by the lower bound of \cite{Ts}, and whose scalar version is due to Jin \cite[Sections~7.3 and~7.4]{Jin} and, for the quantum annulus, also to Lorist and Schwenninger \cite{LS} (Theorem~B);
\item for the strip, the lower bound $2$ for the scalar constant, so that the scalar and the complete constant both equal $2$ (Proposition~\ref{prop:stripopt});
\item for sectors of opening $2\alpha\in[\pi/2,2\alpha_*)$ and lenses and crescents of vertex angle $2\alpha$, the complete bound $\sqrt2$, which is optimal at the right angle, by the lower bound of \cite[Lemma~4.1]{Cr16}, and improves on the bound \eqref{eq:crouzeixbound} of \cite{Cr25} (Theorem~C and Corollary~\ref{cor:lens});
\item for sectors of opening $2\alpha<\pi/2$ and lenses and crescents of vertex angle $2\alpha$, the complete bound $\sqrt{4-8\alpha/\pi}$, which improves on the bound $2$ of \cite{OAI,OAI2} and on \eqref{eq:crouzeixbound}, and the lower bound $2^{1-2\alpha/\pi}$, which improves on the lower bound $\pi\sin\alpha/(2\alpha)$ of \cite{CD03}, \cite[Theorem~1]{BC06} and \cite[Lemma~4.1]{Cr16}, and, numerically, on the constant for $2\times2$ matrices of \cite[Theorem~4.2]{Cr04}; together they determine the scalar and the complete constant up to a factor $1.031$ (Theorem~C and Corollary~\ref{cor:lens});
\item for $2\alpha\in(0,2\alpha_*)$, the same bounds for m-sectorial operators $A$ with $W(A)\subseteq\Sbar\alpha$ and, if they are injective, for their $H^\infty(S_\theta)$-calculus, $\theta>\alpha$, by adapting \cite[Section~3]{Cr25} and \cite[Theorem~3.13]{MR25} (Corollary~\ref{cor:unbounded} and Remark~\ref{rem:hinfty}). Examples are Schr\"odinger operators with bounded complex potentials, after a shift, and divergence form operators with complex coefficients (Section~\ref{sec:examples});
\item for Jordan domains $\Om$ with regular real-analytic boundary and operators with a Faber $\rho$-dilation, $\rho>2$, the complete bound $\max\{\rho,\sqrt{2+2V(\Om)^2}\}$, which is $\rho$ when $V(\Om)\le1$, as for every convex $\Om$ and for some nonconvex $\Om$ (Theorem~D, Lemma~\ref{lem:rhoOm}(c), (d) and the example after it). The complete bound $\rho$ was known for the disk \cite{OA}. For scalar functions, the bound $\rho$ holds for every such $\Om$, with no condition on $V(\Om)$ \cite[Theorem~1]{CG26}. Earlier complete bounds follow from \cite[Theorem~2]{DEY17}, with a constant that is not explicit (Remark~\ref{rem:special}), and from \cite[Theorem~2]{CG19}, extended as in \cite[Section~2]{CK25}, with an explicit constant larger than $\rho$, which the constant of Theorem~D improves on when $V(\Om)\le1$ or $\rho\ge\sqrt{2+2V(\Om)^2}$; for the other cases, see Remark~\ref{rem:CG};
\item for the same domains, operators with a Faber $\rho$-dilation and $\sqrt2\le\rho<2$, the same complete bound (Theorem~D). For $\rho\le2$, \cite[Corollary~6.3]{OAI2} gives the complete bound $2$ for every such $\Om$ (Remark~\ref{rem:rho2}); for $\rho<2$, Theorem~D improves on it when $V(\Om)<1$, and, already for scalar functions, on the bound $1+\sqrt{\rho/2}$ of \cite[Theorem~1]{CG26} when $2+2V(\Om)^2<(1+\sqrt{\rho/2})^2$, for instance for ellipses whose axes have ratio at most $3+2\sqrt2$ (Remark~\ref{rem:CG}). For the disk, it is the bound $\rho$ of \cite{OA}.
\end{enumerate}

\subsection*{Open questions}
For sectors of opening $2\alpha<\pi/2$, the bounds $2^{1-2\alpha/\pi}$ and $\sqrt{4-8\alpha/\pi}$ agree at both ends of this range, and we do not know the optimal constant; to go below $\sqrt{4-8\alpha/\pi}$ with our method, one would need spaces that depend on the operator (Proposition~\ref{prop:spaceopt}), two spaces (Remark~\ref{rem:twospaces}), or a sharper form of the second step of the proof of Theorem~A, where the factor $\sqrt2$ enters; by Remark~\ref{rem:floor}, such a step must use more than the accretivity of the moment map. For openings $2\alpha>\pi/2$, the optimal constant is also unknown; it lies between the constant for $2\times2$ matrices of \cite[Theorem~4.2]{Cr04} and the minimum of $\sqrt2$ and \eqref{eq:crouzeixbound}. For Faber $\rho$-dilations with $\rho>2$, we have not been able to obtain the complete bound $\rho$ with our method without the condition $V(\Om)\le1$ or $\rho\ge\sqrt{2+2V(\Om)^2}$. The constant of Theorem~D is optimal for the disk (Theorem~\ref{thm:convex}) and, when $\rho=2$, for ellipses, by \cite[Theorem~6]{Jin}; for $\rho\le2$ and $V(\Om)>1$ it is not (Remark~\ref{rem:rho2}). We have not attempted to determine whether it is optimal in the other cases.

\subsection*{Organization}
Section~\ref{sec:setting} sets up the framework and Section~\ref{sec:engine} proves Theorem~A. Section~\ref{sec:positivity} proves that the moment map of the kernel $k_r=1/(r(z)+\overline{r(w)})$, for holomorphic $r:\Om\to\{\Rp z>0\}$, is accretive when $\Rp r(T)\ge0$, and so are those of sums of such kernels and of the kernels $1+\int_\T r_t(z)\overline{r_t(w)}\,d\omega(t)$ when $\Rp r_t(T)\ge0$ for all $t\in\T$; it also treats changes of variable. The accretivity for sectors of opening at most $\pi/2$ needs a different argument, given in Section~\ref{ssec:sector} (Proposition~\ref{prop:rightangle}). Section~\ref{sec:twodisks} proves Theorems~B and~C, with separate subsections for the annulus, the strip and sectors, and Section~\ref{sec:convex} proves Theorem~D. Sections~\ref{sec:twodisks} and~\ref{sec:convex} apply the scheme of Sections~\ref{sec:engine} and~\ref{sec:positivity} in each configuration: a kernel whose moment maps are accretive throughout the class, a bound for the multiplier constant of its space, and Theorem~A. A reader interested in one configuration can skip the parts of Sections~\ref{sec:twodisks} and~\ref{sec:convex} that concern the others; each of Sections~\ref{ssec:annulus}--\ref{ssec:sector} needs only Section~\ref{ssec:normal}, and Section~\ref{ssec:sector} also uses the Hardy space of the strip from Section~\ref{ssec:strip}. Section~\ref{sec:convex} depends only on Sections~\ref{sec:setting}--\ref{sec:positivity}, but it uses more machinery than the rest of the paper: Faber polynomials, Grunsky coefficients and the double-layer potential of $\partial\Om$. Section~\ref{sec:examples} treats Schr\"odinger operators, divergence form operators and hyponormal operators. Appendix~\ref{app:spaceopt} proves Proposition~\ref{prop:spaceopt}.
\section{Setting}\label{sec:setting}

\subsection{Conventions}\label{ssec:conv}
Throughout, $\Om\subseteq\C$ is a domain. Hilbert spaces are complex and inner products are linear in the second variable. For a Hilbert space $H$, $\B(H)$ denotes the bounded operators and $\Sc(H)$ the Hilbert--Schmidt operators, with $\ip{X}{Y}_{\mathrm{HS}}=\tr(X^*Y)$. For $x,y\in H$, $xy^*$ is the rank-one operator $v\mapsto\ip{y}{v}x$. We write $\Rp Y=\frac12(Y+Y^*)$ and $\Ip Y=\frac1{2i}(Y-Y^*)$, and $Y\ge0$ means that $Y$ is self-adjoint and positive semidefinite. The numerical range of $T\in\B(H)$ is $W(T)=\{\ip{x}{Tx}:\norm x=1\}$. For bounded operators on or between Hilbert--Schmidt classes, in particular on $\Sc(H)$, we write the adjoint with a dagger, so that $\Rp\mathcal A=\frac12(\mathcal A+\mathcal A^\dagger)$; operators on $H$, or on a larger Hilbert space, keep the star. A bounded operator $\mathcal A$ is \emph{accretive} if $\Rp\mathcal A\ge0$. If $\mathcal A$ is invertible, then $\Rp(\mathcal A^{-1})=(\mathcal A^{-1})^\dagger(\Rp\mathcal A)\mathcal A^{-1}$, so inverses of invertible accretive operators are accretive.

For $B\in\B(H)$ let $\Co B(X)=BXB^*$ on $\Sc(H)$. Then $\norm{\Co B}\le\norm B^2$, $\Co B^\dagger=\Co{B^*}$ and $\Co B\Co{B'}=\Co{BB'}$. Each $\Co B$ maps positive operators to positive operators. If $P\ge0$, then $\Co P\ge0$ as an operator on $\Sc(H)$, since $\ip{X}{\Co P(X)}_{\mathrm{HS}}=\tr(X^*PXP)=\norm{P^{1/2}XP^{1/2}}_{\mathrm{HS}}^2$.

The open unit disk is $\D$, and $\T=\partial\D$ carries the normalized arclength measure $\omega$. For $0<\alpha\le\pi$, $S_\alpha=\{z\ne0:|\arg z|<\alpha\}$ is the open sector of half-aperture $\alpha$ (opening $2\alpha$), and $\Sbar\alpha=\{z:|\arg z|\le\alpha\}\cup\{0\}$ is its closure; $\arg$, $\log$ and powers $z^s$ are principal branches.

\subsection{Holomorphic calculus}\label{ssec:calc}
For $T\in\B(H)$ with $\sigma(T)\subset\Om$ and $f$ holomorphic on $\Om$, $f(T)$ is defined by the Riesz--Dunford integral
\[
f(T)=\frac1{2\pi i}\int_\Gamma f(z)R_T(z)\,dz,\qquad R_T(z)=(z-T)^{-1},
\]
where $\Gamma$ is a finite union of piecewise smooth closed curves in $\Om\setminus\sigma(T)$ with winding number $1$ about $\sigma(T)$ and $0$ about $\C\setminus\Om$. The calculus is linear, multiplicative and independent of $\Gamma$, satisfies the spectral mapping theorem and the composition rule $(g\circ\varphi)(T)=g(\varphi(T))$ for holomorphic $\varphi$ on $\Om$ and $g$ on an open set containing $\varphi(\Om)$, and $f_n(T)\to f(T)$ in norm if $f_n\to f$ locally uniformly \cite[Section~VII.4]{Con}. For matrix-valued $F$, $F(T)$ acts on $\C^d\otimes H$, with the matrix factor first: $F(T)=\sum_{r,s}E_{rs}\otimes F_{rs}(T)$, with $E_{rs}$ the matrix units. By multiplicativity, $(FG)(T)=F(T)G(T)$.

\subsection{Herglotz functions and operator classes}\label{ssec:test}

\begin{definition}\label{def:test}
A \emph{Herglotz function} on $\Om$ is a holomorphic $r:\Om\to\C$ with $\Rp r>0$, and a \emph{Schur function} is a holomorphic $b:\Om\to\D$. For a family $\Bf$ of Herglotz functions, the \emph{class} of $\Bf$ is
\[
\Kc(\Bf)=\{T\in\B(H):\ \sigma(T)\subset\Om,\ \Rp r(T)\ge0\text{ for every }r\in\Bf\}.
\]
The Hilbert space $H$ is arbitrary and is suppressed from the notation.
\end{definition}

Norm conditions take this form through the Cayley transform, by the following standard identity.

\begin{lemma}\label{lem:cayley}
Let $Y$ be a bounded operator on a Hilbert space with $I-Y$ invertible. Then
\[
\Rp\big[(I+Y)(I-Y)^{-1}\big]=\big((I-Y)^{-1}\big)^*\,(I-Y^*Y)\,(I-Y)^{-1},
\]
so $\Rp[(I+Y)(I-Y)^{-1}]\ge0$ if and only if $\norm Y\le1$.
\end{lemma}

A Schur function $b$ gives the Herglotz function $r=(1+b)/(1-b)$, and every Herglotz function $r$ arises in this way, from $b=(r-1)/(r+1)$. If $\sigma(T)\subset\Om$, then $\sigma(b(T))=b(\sigma(T))\subset\D$ by the spectral mapping theorem, so $I-b(T)$ is invertible and, by multiplicativity, $r(T)=(I+b(T))(I-b(T))^{-1}$; by Lemma~\ref{lem:cayley},
\begin{equation}\label{eq:cayleyfact}
\norm{b(T)}\le1\iff\Rp r(T)\ge0\iff\Rp\big(r(T)^{-1}\big)\ge0,
\end{equation}
where the second equivalence holds because inverses of invertible accretive operators are accretive; here $r(T)^{-1}=(1/r)(T)$, and $1/r$ is again a Herglotz function.

Every normal operator $T$ with $\sigma(T)\subset\Om$, on any Hilbert space, lies in every class $\Kc(\Bf)$. Indeed, for $r\in\Bf$ the operator $r(T)$ is normal, its spectrum $r(\sigma(T))$ lies in the open right half-plane, and the numerical range of a normal operator lies in the closed convex hull of its spectrum; hence $\Rp r(T)\ge0$.

\subsection{Kernels, moment maps and multipliers}\label{ssec:kernels}
Throughout, $\E$ is a Hilbert space of holomorphic functions on $\Om$ with continuous point evaluations and reproducing kernel $k$, so that $f(z)=\ip{k(\cdot,z)}{f}$ for $f\in\E$. We call $k$ a \emph{holomorphic kernel} and write $\E=\E_k$; its diagonal is locally bounded, by the uniform boundedness principle. Conversely, every positive kernel on $\Om$ that is holomorphic in $z$ and antiholomorphic in $w$ is a holomorphic kernel: by Hartogs' theorem \cite[Theorem~2.2.8]{Ho} its diagonal is locally bounded, so norm convergence in its reproducing kernel Hilbert space $\E_k$ \cite[Chapter~2]{PR} implies locally uniform convergence, and $\E_k$, the closed span of the holomorphic functions $k(\cdot,w)$, consists of holomorphic functions. $\E$ is separable (a function orthogonal to $k(\cdot,w)$ for $w$ in a countable dense set vanishes), and we fix an orthonormal basis $(e_j)_j$. Put $\E_d=\E\otimes M_d$ with the Hilbert--Schmidt norm on $M_d$: an element is $Q=\sum_jQ_je_j$ with $\norm{Q}_{\E_d}^2=\sum_j\hs{Q_j}^2$, and $Q(z)=\sum_jQ_je_j(z)$.

\begin{lemma}\label{lem:conv}
Let $\E$ and its orthonormal basis $(e_j)$ be as above, and let $T\in\B(H)$ with $\sigma(T)\subset\Om$.
\begin{enumerate}[(a)]
\item $k(z,w)=\sum_je_j(z)\overline{e_j(w)}$ for $z,w\in\Om$; in particular $\sum_j|e_j(z)|^2=k(z,z)$.
\item $C_T^2:=\sum_j\norm{e_j(T)}^2<\infty$.
\item The \emph{moment map} $\Psi_T(X)=\sum_je_j(T)Xe_j(T)^*$ converges absolutely in $\B(\Sc(H))$ and maps positive operators to positive operators.
\item For $Q=\sum_jQ_je_j\in\E_d$, the series $\sum_jQ_j\otimes e_j(T)$ converges absolutely in norm, its sum is $Q(T)$, and $\norm{Q(T)}\le C_T\norm{Q}_{\E_d}$.
\end{enumerate}
\end{lemma}

\begin{proof}
(a) is the expansion of a reproducing kernel in an orthonormal basis \cite[Chapter~2]{PR}.

(b) Let $\Gamma$ be as in Section~\ref{ssec:calc}, with length $L$, and let $c_\Gamma=\max_\Gamma\norm{R_T}$. For $f$ holomorphic on $\Om$, the Cauchy--Schwarz inequality gives $\norm{f(T)}^2\le(2\pi)^{-2}c_\Gamma^2L\int_\Gamma|f|^2\,|dz|$. Summing over $f=e_j$ and using (a),
\[
C_T^2\le(2\pi)^{-2}c_\Gamma^2L\int_\Gamma k(z,z)\,|dz|<\infty,
\]
as $k(z,z)$ is bounded on the compact set $\Gamma$.

(c) The map $\Co{e_j(T)}$ has norm at most $\norm{e_j(T)}^2$ on $\Sc(H)$ and preserves positivity (Section~\ref{ssec:conv}), so the series converges absolutely by (b).

(d) $\sum_j\norm{Q_j\otimes e_j(T)}\le\sum_j\hs{Q_j}\norm{e_j(T)}\le C_T\norm{Q}_{\E_d}$. The partial sums of $\sum_jQ_je_j$ converge to $Q$ in $\E_d$, hence uniformly on $\Gamma$, because $\hs{Q(z)}\le\norm{Q}_{\E_d}k(z,z)^{1/2}$ by (a) and $k(z,z)$ is bounded on $\Gamma$; so their Riesz--Dunford integrals converge to $Q(T)$.
\end{proof}

\begin{lemma}\label{lem:cauchy}
Let $T\in\B(H)$ with $\sigma(T)\subset\Om$, and let $\Gamma$ be as in Section~\ref{ssec:calc}. Then, for $X\in\Sc(H)$,
\begin{equation}\label{eq:cauchy}
\Psi_T(X)=\frac1{4\pi^2}\int_\Gamma\int_\Gamma k(z,w)\,R_T(z)\,X\,R_T(w)^*\,dz\,\overline{dw},
\end{equation}
where $\overline{dw}=\overline{\gamma'(s)}\,ds$ for a parametrization $w=\gamma(s)$. In particular $\Psi_T$ depends only on $k$.
\end{lemma}

\begin{proof}
Taking adjoints in the Riesz--Dunford integral gives, for $f$ holomorphic on $\Om$,
\begin{equation}\label{eq:fXf}
f(T)\,X\,f(T)^*=\frac1{4\pi^2}\int_\Gamma\int_\Gamma f(z)\overline{f(w)}\,R_T(z)XR_T(w)^*\,dz\,\overline{dw}.
\end{equation}
Apply this to $f=e_j$ and sum over $j$; by Lemma~\ref{lem:conv}(c) the left sides sum to $\Psi_T(X)$. By Lemma~\ref{lem:conv}(a) and the Cauchy--Schwarz inequality, the partial sums of $\sum_je_j(z)\overline{e_j(w)}$ converge to $k(z,w)$ and are bounded by $k(z,z)^{1/2}k(w,w)^{1/2}$, which is bounded on $\Gamma\times\Gamma$, so dominated convergence gives \eqref{eq:cauchy}.
\end{proof}

Lemma~\ref{lem:cauchy} shows that $\Psi_T$ is additive in the kernel: the moment map of $k_1+k_2$ is the sum of the moment maps of $k_1$ and $k_2$, and the constant kernel $c\ge0$ has moment map $X\mapsto cX$. Most kernels below are given as integrals of products, and the following lemma identifies both their spaces and their moment maps.

\begin{lemma}\label{lem:intkernel}
Let $(\Lambda,\nu)$ be a $\sigma$-finite measure space and let $\phi_\lambda$ $(\lambda\in\Lambda)$ be holomorphic functions on $\Om$, with $(\lambda,z)\mapsto\phi_\lambda(z)$ measurable and $\int_\Lambda|\phi_\lambda(z)|^2\,d\nu(\lambda)$ locally bounded in $z$. Put
\[
k(z,w)=\int_\Lambda\phi_\lambda(z)\overline{\phi_\lambda(w)}\,d\nu(\lambda).
\]
Then $k$ is a holomorphic kernel, and for $T\in\B(H)$ with $\sigma(T)\subset\Om$ the integral of the maps $X\mapsto\phi_\lambda(T)\,X\,\phi_\lambda(T)^*$ with respect to $\nu$ converges absolutely in the norm of $\B(\Sc(H))$ and equals $\Psi_T$.
\end{lemma}

\begin{proof}
The kernel is positive, holomorphic in $z$ by the Cauchy integral formula and Fubini's theorem, and antiholomorphic in $w$ by symmetry, hence a holomorphic kernel (Section~\ref{ssec:kernels}). As in the proof of Lemma~\ref{lem:conv}(b), with the same $\Gamma$, $L$ and $c_\Gamma$, $\int_\Lambda\norm{\phi_\lambda(T)}^2\,d\nu\le(2\pi)^{-2}c_\Gamma^2L\int_\Gamma k(z,z)\,|dz|<\infty$, so the integral converges absolutely. Integrating \eqref{eq:fXf} for $f=\phi_\lambda$ over $\Lambda$, Fubini's theorem, which applies since $\int_\Lambda|\phi_\lambda(z)\phi_\lambda(w)|\,d\nu\le k(z,z)^{1/2}k(w,w)^{1/2}$, and \eqref{eq:cauchy} show that the integral equals $\Psi_T$.
\end{proof}

The space of such a kernel has the standard description by a Kolmogorov decomposition \cite[Chapter~2]{PR}, here $z\mapsto\overline{\phi_\cdot(z)}\in L^2(\nu)$: if $g\mapsto\int_\Lambda g\,\phi_\lambda\,d\nu$ is injective on $L^2(\nu)$, then $\E_k$ consists of the functions $\int_\Lambda g\,\phi_\lambda\,d\nu$, $g\in L^2(\nu)$, with norm $\norm g_{L^2(\nu)}$. Below, $\Lambda$ is $\RR$ or a subset of $\Z$, with Lebesgue, respectively counting, measure $d\lambda$, and $d\nu=\beta^{-2}d\lambda$ for a measurable function $\beta>0$ on $\Lambda$; writing $\hat f=\beta^{-2}g$, this says that, if $\hat f\mapsto\int\hat f\phi_\lambda\,d\lambda$ is injective on $\{\hat f:\int|\hat f|^2\beta^2\,d\lambda<\infty\}$, then
\begin{equation}\label{eq:intspace}
\E_k=\Big\{\int\hat f\phi_\lambda\,d\lambda:\ \int|\hat f|^2\beta^2\,d\lambda<\infty\Big\},\qquad\Big\|\int\hat f\phi_\lambda\,d\lambda\Big\|_{\E_k}^2=\int|\hat f|^2\beta^2\,d\lambda .
\end{equation}
For the Laurent, Fourier and Mellin expansions of Section~\ref{sec:twodisks}, injectivity holds by the uniqueness of Laurent coefficients and of the Fourier transform.

\begin{definition}\label{def:mult}
For a bounded holomorphic $F:\Om\to M_d$, let $m_\ell(F)$ be the norm of $Q\mapsto FQ$ on $\E_d$, and $m_r(F)$ that of $Q\mapsto QF$, with the value $\infty$ if the product does not always lie in $\E_d$. The \emph{multiplier constant} of $\E\ne\{0\}$ is
\begin{multline*}
M(\E)=\sup\big\{\max(m_\ell(F),m_r(F)):\ d\ge1,\\* F:\Om\to M_d\text{ holomorphic},\ \sup_\Om\norm F\le1\big\}\in[1,\infty].
\end{multline*}
\end{definition}

\begin{remark}\label{rem:usualmult}
Writing $Q\in\E_d$ by its columns, $m_\ell(F)$ is the norm of multiplication by $F$ on the space $\E\otimes\C^d$ of $\C^d$-valued functions with entries in $\E$, and, by transposition, $m_r(F)$ is that of $F^{\mathsf T}$, and $\sup_\Om\norm{F^{\mathsf T}}=\sup_\Om\norm F$.
\end{remark}

In all uses of the following lemma, the norm $\norm\cdot_{\mathcal N}$ is a Hardy norm.

\begin{lemma}\label{lem:compare}
Let $(\mu_\iota)$ be positive measures on $\Om$, indexed by a directed set, and put $\norm f_{\mathcal N}^2=\sup_\iota\int_\Om|f|^2\,d\mu_\iota$. Suppose that $\iota\mapsto\int_\Om|f|^2\,d\mu_\iota$ is nondecreasing for every $f\in\E$, that $\E$ consists of the holomorphic functions $f$ on $\Om$ with $\norm f_{\mathcal N}<\infty$, and that, for some $0<c_-\le c_+<\infty$,
\[
c_-\norm f_\E^2\le\norm f_{\mathcal N}^2\le c_+\norm f_\E^2\qquad(f\in\E).
\]
Then $M(\E)\le\sqrt{c_+/c_-}$.
\end{lemma}

\begin{proof}
Extend $\norm\cdot_{\mathcal N}$ to $M_d$-valued functions with $\hs\cdot$ in place of $|\cdot|$. By monotonicity, $\norm Q_{\mathcal N}^2=\sum_{r,s}\norm{Q_{rs}}_{\mathcal N}^2$ for the entries $Q_{rs}\in\E$ of $Q$; since $\norm Q_{\E_d}^2=\sum_{r,s}\norm{Q_{rs}}_\E^2$, the two inequalities hold on $\E_d$. If $\sup\norm F\le1$ and $Q\in\E_d$, then $\hs{F(z)Q(z)}\le\hs{Q(z)}$ pointwise, so $\norm{FQ}_{\mathcal N}\le\norm{Q}_{\mathcal N}<\infty$; hence $FQ\in\E_d$ and $\norm{FQ}_{\E_d}^2\le c_-^{-1}\norm{Q}_{\mathcal N}^2\le(c_+/c_-)\norm{Q}_{\E_d}^2$. Likewise for $QF$.
\end{proof}

\begin{definition}\label{def:nondeg}
$\E$ is \emph{nondegenerate} if there are $q_n\in\E$ with $q_n\to1$ uniformly on compact subsets of $\Om$.
\end{definition}

If $\E$ is nondegenerate and $\sigma(T)\subset\Om$, then $\norm{q_n(T)-I}\to0$ by the Riesz--Dunford integral, so $\Rp q(T)\ge\frac12I$ for some $q\in\E$.

\section{The moment map theorem}\label{sec:engine}

Throughout this section, $\E$, $k$, $(e_j)$ and $\E_d$ are as in Section~\ref{ssec:kernels}. For $T\in\B(H)$ with $\sigma(T)\subset\Om$, recall that $C_T=\big(\sum_j\norm{e_j(T)}^2\big)^{1/2}<\infty$ (Lemma~\ref{lem:conv}(b)). For $d\ge1$ and $x=(x_1,\dots,x_d)$, $y=(y_1,\dots,y_d)$ in $\C^d\otimes H$, put
\[
\mu_j(x,y)=\big[\ip{x_r}{e_j(T)y_s}\big]_{r,s=1}^d\in M_d,\qquad J(x,y)=\Big(\sum_j\hs{\mu_j(x,y)}^2\Big)^{1/2},
\]
and $P_x=\sum_rx_rx_r^*$, a positive finite-rank operator on $H$.

\begin{lemma}\label{lem:coef}
For $x,y\in\C^d\otimes H$,
\[
J(x,y)=\sup\big\{|\ip{x}{Q(T)y}|:\ Q\in\E_d,\ \norm{Q}_{\E_d}\le1\big\}\le C_T\norm{x}\norm{y}.
\]
\end{lemma}

\begin{proof}
Let $Q=\sum_jQ_je_j\in\E_d$. Then $(Q(T)y)_r=\sum_j\sum_s(Q_j)_{rs}e_j(T)y_s$ by Lemma~\ref{lem:conv}(d), so $\ip{x}{Q(T)y}=\sum_j\sum_{r,s}(Q_j)_{rs}\,\mu_j(x,y)_{rs}$ and $|\ip{x}{Q(T)y}|\le J(x,y)\norm Q_{\E_d}$ by the Cauchy--Schwarz inequality; equality is approached by $Q_j=\overline{\mu_j(x,y)}$ (entrywise conjugate) for $j$ in a finite set and $Q_j=0$ otherwise, normalized. Since $|\mu_j(x,y)_{rs}|\le\norm{e_j(T)}\norm{x_r}\norm{y_s}$, summing the squares over $r$, $s$ and $j$ gives $J(x,y)\le C_T\norm x\norm y$.
\end{proof}

\begin{lemma}\label{lem:momentineq}
$J(x,y)^2=\ip{P_x}{\Psi_T(P_y)}_{\mathrm{HS}}$. If $\Psi_T$ is accretive, then
\[
J(x,y)^2+J(y,x)^2\le2\,J(x,x)\,J(y,y)\qquad(x,y\in \C^d\otimes H).
\]
\end{lemma}

\begin{proof}
$\tr\big(x_rx_r^*\,e_j(T)y_sy_s^*e_j(T)^*\big)=|\ip{x_r}{e_j(T)y_s}|^2$, and $X\mapsto\tr(P_xX)=\ip{P_x}{X}_{\mathrm{HS}}$ is bounded on $\Sc(H)$, so the sum over $j,r,s$ gives the identity, by Lemma~\ref{lem:conv}(c). In particular $\ip{P_x}{\Psi_T(P_y)}_{\mathrm{HS}}\ge0$ for all $x,y$.

Let $\mathcal R=\Rp\Psi_T=\frac12(\Psi_T+\Psi_T^\dagger)$, which is $\ge0$ by hypothesis. Since $\ip{P_y}{\Psi_T(P_x)}_{\mathrm{HS}}=J(y,x)^2$ is real,
\[
\ip{P_x}{\mathcal R(P_y)}_{\mathrm{HS}}=\tfrac12\big(J(x,y)^2+J(y,x)^2\big),\qquad \ip{P_x}{\mathcal R(P_x)}_{\mathrm{HS}}=J(x,x)^2 .
\]
The Cauchy--Schwarz inequality for the positive form $\ip{\cdot}{\mathcal R(\cdot)}_{\mathrm{HS}}$ finishes the proof.
\end{proof}

Recall from Definition~\ref{def:mult} that $m_\ell(F)$ and $m_r(F)$ are the norms of left and right multiplication by $F$ on $\E_d$.

\begin{theorem}\label{thm:engine}
Let $T\in\B(H)$ with $\sigma(T)\subset\Om$. Assume that $\Psi_T$ is accretive on $\Sc(H)$ and that $\Rp q(T)\ge\delta I$ for some $q\in\E$ and $\delta>0$. Then for every $d\ge1$ and every bounded holomorphic $F:\Om\to M_d$ with $m_\ell(F),m_r(F)<\infty$,
\[
\norm{F(T)}^2\le2\,m_\ell(F)\,m_r(F).
\]
\end{theorem}

\begin{proof}
Let $B=F(T)$ and $\gamma=\norm B$; we may assume $\gamma>0$. Choose unit vectors $\xi_n\in \C^d\otimes H$ with $\norm{B\xi_n}\to\gamma$. Put $\eta_n=\gamma^{-1}B\xi_n$ and $\Delta_n=B^*\eta_n-\gamma\xi_n$. Since $\ip{B^*\eta_n}{\xi_n}=\gamma\norm{\eta_n}^2$,
\[
\norm{\Delta_n}^2=\norm{B^*\eta_n}^2-2\gamma^2\norm{\eta_n}^2+\gamma^2\le\gamma^2(1-\norm{\eta_n}^2)\to0,\qquad \norm{\eta_n}\to1 .
\]
For $Q\in\E_d$, $(QF)(T)=Q(T)B$ and $(FQ)(T)=BQ(T)$, so
\begin{align*}
\gamma\ip{\eta_n}{Q(T)\eta_n}&=\ip{\eta_n}{(QF)(T)\xi_n},\\
\gamma\ip{\xi_n}{Q(T)\xi_n}&=\ip{\eta_n}{(FQ)(T)\xi_n}-\ip{\Delta_n}{Q(T)\xi_n}.
\end{align*}
Take the supremum over $\norm{Q}_{\E_d}\le1$. By Lemma~\ref{lem:coef}, the left sides become $\gamma J(\eta_n,\eta_n)$ and $\gamma J(\xi_n,\xi_n)$; the first terms on the right are at most $J(\eta_n,\xi_n)$ times $\norm{QF}_{\E_d}\le m_r(F)$, respectively $\norm{FQ}_{\E_d}\le m_\ell(F)$; and $|\ip{\Delta_n}{Q(T)\xi_n}|\le C_T\norm{\Delta_n}$. With $a_n=J(\eta_n,\xi_n)$ and $\epsilon_n=C_T\norm{\Delta_n}\to0$,
\begin{equation}\label{eq:twotests}
\gamma\, J(\eta_n,\eta_n)\le m_r(F)\,a_n,\qquad \gamma\, J(\xi_n,\xi_n)\le m_\ell(F)\,a_n+\epsilon_n .
\end{equation}
Applying Lemma~\ref{lem:coef} to $Q=qI_d/(\sqrt d\norm q_\E)$, which has norm $1$ in $\E_d$, gives, for every $\xi\in \C^d\otimes H$,
\[
J(\xi,\xi)\ge\frac{\Rp\ip{\xi}{(I_d\otimes q(T))\xi}}{\sqrt d\,\norm q_\E}\ge\kappa\norm\xi^2,\qquad\kappa=\frac{\delta}{\sqrt d\,\norm q_\E}.
\]
Since $\norm{\xi_n}=1$, the second inequality in \eqref{eq:twotests} gives $m_\ell(F)a_n\ge\gamma\kappa-\epsilon_n\ge\gamma\kappa/2$ for all large $n$. Hence $m_\ell(F)>0$ and $a_n\ge c_0:=\gamma\kappa/(2m_\ell(F))>0$ for all large $n$. Lemma~\ref{lem:momentineq} with $x=\xi_n$ and $y=\eta_n$, with its nonnegative term $J(\xi_n,\eta_n)^2$ dropped, and \eqref{eq:twotests} give
\[
\gamma^2a_n^2\le2\,\big(\gamma J(\xi_n,\xi_n)\big)\big(\gamma J(\eta_n,\eta_n)\big)\le2\big(m_\ell(F)a_n+\epsilon_n\big)m_r(F)a_n .
\]
Divide by $a_n^2$, use $a_n\ge c_0$, so that $\epsilon_n/a_n\to0$, and let $n\to\infty$.
\end{proof}

\begin{corollary}\label{cor:main}
Let $\Kc$ be a class of operators with spectrum in $\Om$, and let $\E$ be nondegenerate with $\Psi_T$ accretive for every $T\in\Kc$. Then for every $T\in\Kc$, every $d\ge1$ and every bounded holomorphic $F:\Om\to M_d$,
\[
\norm{F(T)}^2\le2\,m_\ell(F)\,m_r(F)\le2\,M(\E)^2\sup_\Om\norm F^2 .
\]
Consequently
\[
K_{\mathrm{cb}}(\Kc)\le\sqrt2\,\inf\big\{M(\E):\ \E\text{ nondegenerate, }\Psi_T\text{ accretive for every }T\in\Kc\big\}.
\]
\end{corollary}

\begin{proof}
By Definition~\ref{def:nondeg} and the paragraph following it, there is $q\in\E$ with $\Rp q(T)\ge\frac12I$. We may assume $F\ne0$; then $Fq\ne0$, so taking $Q=qI_d$ gives $m_\ell(F),m_r(F)>0$. If both are finite, Theorem~\ref{thm:engine} gives the first inequality; otherwise its right side is infinite. The second inequality is Definition~\ref{def:mult}, and the last assertion follows from the definition of $K_{\mathrm{cb}}(\Kc)$ (Section~\ref{sec:intro}).
\end{proof}

\begin{remark}\label{rem:complete}
Theorem~\ref{thm:engine} bounds $M_d$-valued $F$ for every $d$ under hypotheses on $T$ that do not involve $d$. In its proof, Lemma~\ref{lem:momentineq} is applied to vectors $x,y\in\C^d\otimes H$. Since $J(x,y)^2=\ip{P_x}{\Psi_T(P_y)}_{\mathrm{HS}}$ with $P_x,P_y\in\Sc(H)$, the lemma follows from the Cauchy--Schwarz inequality for the Hermitian form $\ip{X}{(\Rp\Psi_T)(Y)}_{\mathrm{HS}}$ on $\Sc(H)$, positive semidefinite since $\Psi_T$ is accretive, at $X=P_x$ and $Y=P_y$; this inequality holds for all $X,Y\in\Sc(H)$ with constant $1$, so the accretivity of $\Psi_T$ on $\Sc(H)$ suffices for every $d$. In the conclusion $\norm{F(T)}^2\le2\,m_\ell(F)\,m_r(F)$, the multiplier norms are at most $M(\E)\sup_\Om\norm F$ for every $d$, and the constant $\kappa$ of the proof, which depends on $d$, does not appear. So Corollary~\ref{cor:main} gives complete bounds whenever $M(\E)<\infty$, with no separate passage from scalar to matrix-valued functions. In Sections~\ref{sec:twodisks} and~\ref{sec:convex}, $M(\E)$ is bounded by comparing the norm of $\E$ with a boundary $L^2$ norm, which handles all $d$ at once (Lemma~\ref{lem:compare} and the proof of Theorem~\ref{thm:convex}).
\end{remark}

\begin{remark}\label{rem:engine}
The proof of Theorem~\ref{thm:engine} follows the two-step scheme of \cite[Section~3, proof of Theorem~3.2]{OAI} described in Section~\ref{sec:intro}, with the vectors $\xi_n$ and $\eta_n$ in place of $x$ and $y$. The essential difference is in the second step. Since $J(x,y)^2=\ip{P_x}{\Psi_T(P_y)}_{\mathrm{HS}}$, the comparison of $J(y,x)$ with $J(x,x)$ and $J(y,y)$ is an inequality for the moment map, which Lemma~\ref{lem:momentineq} proves from the accretivity of $\Psi_T$ alone; so the argument is the same for every kernel, and the inequalities that define a class are used only to prove the accretivity of $\Psi_T$ (Sections~\ref{sec:positivity}--\ref{sec:convex}). In \cite{OAI}, for a matrix $A$ with $W(A)$ in a bounded convex domain $\Om$ with regular real-analytic boundary, the comparison is obtained from $\Rp r_t(A)\ge0$ $(t\in\T)$ by a trace inequality integrated over $\T$; the integrated inequality is $\Rp\ip{X}{\Psi_A(X)}_{\mathrm{HS}}\ge0$ at $X=P_x-P_y$, for the moment map $\Psi_A$ of the kernel $2+\sum_{n\ge1}\Phi_n(z)\overline{\Phi_n(w)}$ (Remark~\ref{rem:convex}). In the first step, \cite{OAI} uses two equivalent norms on the matrix polynomials $G=\sum_{n\ge0}\Phi_nC_n$ $(C_n\in M_d,\ \Phi_0=1)$, namely $\big(\sum_nw_n\hs{C_n}^2\big)^{1/2}$ with $w_n=1$ for all $n$, and with $w_0=1$ and $w_n=2$ for $n\ge1$. The norm of $G\mapsto\ip{y}{G(A)x}$ is taken with respect to the first, those of $G\mapsto\ip{x}{G(A)x}$ and $G\mapsto\ip{y}{G(A)y}$ with respect to the second, and the norm estimates for Faber series \cite[Proposition~2.2]{OAI} show that multiplication by an $F$ with $\sup_\Om\norm F\le1$ is a contraction from the second norm to the first; Theorem~\ref{thm:engine} uses a single space $\E$, with the norms $m_\ell(F)$ and $m_r(F)$ of multiplication by $F$ on $\E_d$ (see Remark~\ref{rem:twospaces}). Finally, \cite{OAI} keeps $J(y,x)$ away from $0$ by evaluating $G\mapsto\ip{y}{G(A)x}$ at $G=F$, which requires $F$ to lie in the space; this is not available for the spaces of the strip and the sectors in Section~\ref{sec:twodisks}, which do not contain the constants. The proof of Theorem~\ref{thm:engine} uses instead a function $q\in\E$ with $\Rp q(T)\ge\delta I$, which nondegeneracy provides.
\end{remark}

\begin{remark}\label{rem:floor}
Since $M(\E)\ge1$, Corollary~\ref{cor:main} never gives a constant below $\sqrt2$; for $\rho$-contractions with $1\le\rho<\sqrt2$, contractions included, this is weaker than the Okubo--Ando constant $\rho$ (Remark~\ref{rem:rho}). For $\Sbar{\pi/4}$ the constant $\sqrt2$ is optimal (Corollary~\ref{cor:sectorclosed}). The factor $\sqrt2$ enters through Lemma~\ref{lem:momentineq}, which bounds the sum $J(x,y)^2+J(y,x)^2$, while the proof of Theorem~\ref{thm:engine} needs a bound for $J(\eta_n,\xi_n)^2$ only. Under the hypotheses of Theorem~\ref{thm:engine}, this cannot be improved, already when $M(\E)=1$: for $\E=H^2(\D)$, the moment map $X\mapsto X+2NXN^*$ at $T=\sqrt2N$, with $N=\big(\begin{smallmatrix}0&1\\0&0\end{smallmatrix}\big)$, is accretive, since $\Rp\ip X{\Psi_T(X)}_{\mathrm{HS}}=\hs X^2+2\Rp\big(\overline{X_{11}}X_{22}\big)$, and $F(z)=z$ gives $\norm{F(T)}=\sqrt2=\sqrt2\,M(\E)\sup_\D|F|$.
\end{remark}

\begin{remark}\label{rem:twospaces}
The proof of Theorem~\ref{thm:engine} applies with two spaces: $\E_1$, with kernel $k_1$, for the functional $Q\mapsto\ip{\eta_n}{Q(T)\xi_n}$, and $\E_2$, with kernel $k_2$, for the other two. Let $J_1,J_2$ and $\Psi^{(1)}_T,\Psi^{(2)}_T$ be defined as $J$ and $\Psi_T$ for $\E_1$ and $\E_2$, and let $m_\ell(F)$ and $m_r(F)$ now be the norms of $Q\mapsto FQ$ and $Q\mapsto QF$ from $\E_2\otimes M_d$ to $\E_1\otimes M_d$. If $\Rp q(T)\ge\delta I$ for some $q\in\E_2$ and $\delta>0$, and, for some $c>0$,
\[
J_1(y,x)^2\le c\,J_2(x,x)\,J_2(y,y)\qquad(d\ge1,\ x,y\in\C^d\otimes H),
\]
then the proof gives $\norm{F(T)}^2\le c\,m_\ell(F)\,m_r(F)$. For one space, Lemma~\ref{lem:momentineq} gives this inequality with $c=2$ when $\Psi_T$ is accretive. For two spaces and $a>0$, the identity
\begin{multline*}
a\ip X{\Psi^{(2)}_TX}+a\ip Y{\Psi^{(2)}_TY}-\ip X{\Psi^{(1)}_TY}-\ip Y{\Psi^{(1)}_TX}\\
=\tfrac12\ip{X+Y}{\big(a\Psi^{(2)}_T-\Psi^{(1)}_T\big)(X+Y)}+\tfrac12\ip{X-Y}{\big(a\Psi^{(2)}_T+\Psi^{(1)}_T\big)(X-Y)}
\end{multline*}
in $\Sc(H)$, at $X=sP_x$ and $Y=s^{-1}P_y$ $(s>0)$, gives it with $c=2a$ when $a\Psi^{(2)}_T+\Psi^{(1)}_T$ and $a\Psi^{(2)}_T-\Psi^{(1)}_T$ are accretive. If instead the moment map of a kernel $k$ is accretive at $T$, and, for some $b_1,b_2>0$, $b_1k-k_1$ and $b_2k_2-k$ are positive kernels, then the inequality holds with $c=2b_1b_2$, and the norms of $Q\mapsto FQ$ and $Q\mapsto QF$ on $\E_k\otimes M_d$ are at most $\sqrt{b_1b_2}\,m_\ell(F)$ and $\sqrt{b_1b_2}\,m_r(F)$; so, if $\E_k$ is nondegenerate, Corollary~\ref{cor:main} for $\E_k$ already gives $\norm{F(T)}^2\le2b_1b_2\,m_\ell(F)\,m_r(F)$. This covers the two norms of \cite{OAI} (Remark~\ref{rem:engine}), for which the inequality comes from the accretivity of the moment map of the kernel $k=2+\sum_{n\ge1}\Phi_n(z)\overline{\Phi_n(w)}$ of Remark~\ref{rem:convex}: they are those of $\E_1=\E_{k-1}$ and $\E_2=\E_{k/2}$, so $b_1=1$, $b_2=2$ and $c=4$. It also covers the accretivity of $a\Psi^{(2)}_T\pm\Psi^{(1)}_T$ when $ak_2-k_1$ is a positive kernel, with $k=k_2$, $b_1=a$ and $b_2=1$, since $\Psi^{(2)}_T$ is then accretive. Sections~\ref{sec:positivity}--\ref{sec:convex} verify the accretivity of single moment maps, so we use one space. If the inequality is verified in another way, for instance through the accretivity of $a\Psi^{(2)}_T\pm\Psi^{(1)}_T$ when $ak_2-k_1$ is not a positive kernel, two spaces might give smaller constants; we have not found a case where they do.
\end{remark}

\section{Accretive moment maps}\label{sec:positivity}

Throughout this section, $T\in\B(H)$ and $\sigma(T)\subset\Om$.

\subsection{Herglotz functions}\label{ssec:cayley}

\begin{proposition}\label{prop:herglotz}
Let $r$ be a Herglotz function on $\Om$. Then
\[
k_r(z,w)=\frac1{r(z)+\overline{r(w)}}=\int_0^\infty e^{-t\,r(z)}\,\overline{e^{-t\,r(w)}}\,dt
\]
is a holomorphic kernel. The map $\mathcal L(X)=r(T)X+Xr(T)^*$ is invertible on $\Sc(H)$, and the moment map of $k_r$ is $\Psi_T=\mathcal L^{-1}$. If $\Rp r(T)\ge0$, then $\Psi_T$ is accretive.
\end{proposition}

\begin{proof}
Lemma~\ref{lem:intkernel} with the functions $e^{-t\,r}$ $(t>0)$, for which $\int_0^\infty|e^{-t\,r(z)}|^2dt=(2\Rp r(z))^{-1}$, shows that $k_r$ is a holomorphic kernel with
\[
\Psi_T(X)=\int_0^\infty e^{-t\,r(T)}Xe^{-t\,r(T)^*}dt;
\]
moreover, $\norm{e^{-t\,r(T)}}\le ce^{-\delta t}$ for some $c,\delta>0$, since $\sigma(r(T))$ is a compact subset of the open right half-plane. Integrating
\[
\frac{d}{dt}\big(e^{-t\,r(T)}Xe^{-t\,r(T)^*}\big)=-\mathcal L\big(e^{-t\,r(T)}Xe^{-t\,r(T)^*}\big)=-e^{-t\,r(T)}\mathcal L(X)e^{-t\,r(T)^*}
\]
over $(0,\infty)$ gives $\mathcal L\Psi_T=\Psi_T\mathcal L=I$. If $\Rp r(T)\ge0$, then $\Rp\ip{X}{\mathcal L(X)}_{\mathrm{HS}}=\norm{(\Rp r(T))^{1/2}X}_{\mathrm{HS}}^2+\norm{X(\Rp r(T))^{1/2}}_{\mathrm{HS}}^2\ge0$, so $\mathcal L$ is accretive, and so is its inverse $\Psi_T$.
\end{proof}

Since $\Rp r(T)\ge0$ implies $\Rp\big[(1/r)(T)\big]\ge0$ by \eqref{eq:cayleyfact}, Proposition~\ref{prop:herglotz} applies to $r$ and to $1/r$; with the additivity noted after Lemma~\ref{lem:cauchy}, we obtain the following.

\begin{corollary}\label{cor:sums}
Let $\Bf$ be a finite family of Herglotz functions on $\Om$, let $c_0,c_r,c_r'\ge0$ $(r\in\Bf)$, and let
\[
k=c_0+\sum_{r\in\Bf}\big(c_rk_r+c_r'k_{1/r}\big).
\]
Then $\Psi_T$ is accretive for every $T\in\Kc(\Bf)$.
\end{corollary}

For a Schur function $b$ and $r=(1+b)/(1-b)$, a direct computation gives
\[
k_r(z,w)=\frac{(1-b(z))(1-\overline{b(w)})}{2(1-b(z)\overline{b(w)})},\qquad k_{1/r}(z,w)=\frac{(1+b(z))(1+\overline{b(w)})}{2(1-b(z)\overline{b(w)})},
\]
the second because $1/r=(1-b)/(1+b)$. Their average is
\begin{equation}\label{eq:schurkernel}
\tfrac12(k_r+k_{1/r})(z,w)=\frac12\,\frac{1+b(z)\overline{b(w)}}{1-b(z)\overline{b(w)}}=\frac12+\sum_{n\ge1}b(z)^n\,\overline{b(w)^n},
\end{equation}
so, by Corollary~\ref{cor:sums} and \eqref{eq:cayleyfact}, its moment map is accretive at every $T$ with $\sigma(T)\subset\Om$ and $\norm{b(T)}\le1$.

\subsection{Families of holomorphic functions}\label{ssec:pencil}

\begin{proposition}\label{prop:pencil}
Let $(r_t)_{t\in\T}$ be holomorphic functions on $\Om$. Suppose that for every compact $C\subset\Om$ there is $\tau_C<1$ such that $(t,z)\mapsto r_t(z)$ extends to $\{|t|>\tau_C\}\times C$ as a bounded continuous function, holomorphic in $t$, that tends to $1$ as $t\to\infty$, uniformly in $z\in C$. Then
\[
k(z,w)=1+\int_\T r_t(z)\overline{r_t(w)}\,d\omega(t)
\]
is a holomorphic kernel on $\Om$, and for every $T$ with $\sigma(T)\subset\Om$,
\begin{equation}\label{eq:pencil}
\Psi_T+\Psi_T^\dagger=\int_\T\Co{N_T(t)}\,d\omega(t),\qquad N_T(t)=r_t(T)+r_t(T)^*,
\end{equation}
where $\Co{N_T(t)}(X)=N_T(t)XN_T(t)$. In particular $\Psi_T$ is accretive whenever $\Rp r_t(T)\ge0$ for all $t\in\T$.
\end{proposition}

\begin{proof}
By Lemma~\ref{lem:intkernel}, applied to $\omega$ and the functions $r_t$, which are bounded on compact sets uniformly in $t$, and by the additivity of the moment map in the kernel (Section~\ref{ssec:kernels}), $k$ is a holomorphic kernel and, writing $r_t$ for $r_t(T)$,
\[
\Psi_T(X)=X+\int_\T r_tXr_t^*\,d\omega(t),\qquad \Psi_T^\dagger(X)=X+\int_\T r_t^*Xr_t\,d\omega(t),
\]
the second since $\Co B^\dagger=\Co{B^*}$.
Next, $\int_\T r_tXr_t\,d\omega=X$. First, for $u,v\in\Om$ the function $g(t)=r_t(u)r_t(v)$ is holomorphic on $\{|t|>\tau\}$ for some $\tau<1$, and bounded with $g(\infty)=1$; its Laurent series contains only nonpositive powers of $t$ and has constant term $1$, so $\int_\T g\,d\omega=1$. Now, with $\Gamma$ as in Section~\ref{ssec:calc},
\[
r_tXr_t=\frac1{(2\pi i)^2}\int_\Gamma\int_\Gamma r_t(z)r_t(w)\,R_T(z)XR_T(w)\,dz\,dw ;
\]
by Fubini's theorem, integrating over $t\in\T$ first replaces $r_t(z)r_t(w)$ by $1$, and then $(2\pi i)^{-1}\int_\Gamma R_T(z)\,dz=I$ gives $X$. Taking adjoints, $\int_\T r_t^*Xr_t^*\,d\omega=X$. Expanding $N_T(t)XN_T(t)$ into its four terms and integrating gives $(\Psi_T(X)-X)+(\Psi_T^\dagger(X)-X)+X+X$, which is \eqref{eq:pencil}. If $N_T(t)\ge0$, then $\Co{N_T(t)}\ge0$ on $\Sc(H)$ (Section~\ref{ssec:conv}), and so is its average.
\end{proof}

\begin{remark}\label{rem:pencilconst}
The constant $1$ in $k$ cannot be lowered: for $\Om=\D$ and $r_t(z)=(1-\bar tz)^{-1}$, the condition $\Rp r_t(T)\ge0$ for all $t\in\T$ means $\Rp(I-\bar tT)\ge0$ for all $t$ (Section~\ref{ssec:conv}), that is, $W(T)\subseteq\overline\D$, and $k=1+(1-z\bar w)^{-1}$. At $T=2N$, where $N=\big(\begin{smallmatrix}0&1\\0&0\end{smallmatrix}\big)$, we have $W(T)=\overline\D$, and for $X=\mathrm{diag}(1,-1)$ the moment map of the kernel $\kappa+(1-z\bar w)^{-1}$ at $T$ satisfies $\Rp\ip{X}{\Psi_T(X)}_{\mathrm{HS}}=2\kappa-2$, since $f(T)=f(0)I+2f'(0)N$, so it is not accretive when $0\le\kappa<1$.
\end{remark}

Unlike the kernel $k_r$ of Proposition~\ref{prop:herglotz}, which is a holomorphic kernel only when $\Rp r>0$ on $\Om$, the kernel $k$ of Proposition~\ref{prop:pencil} is a holomorphic kernel for every family as in that proposition: the $r_t$ need not be Herglotz functions, and the sign of $\Rp r_t$ enters only through the hypothesis $\Rp r_t(T)\ge0$. If the $r_t$ are Herglotz functions, every normal operator with spectrum in $\Om$ satisfies this hypothesis (Section~\ref{ssec:test}). In general, the hypothesis also confines $\sigma(T)$ to the set where $\Rp r_t\ge0$ for all $t\in\T$, since $\sigma(r_t(T))=r_t(\sigma(T))$ and the spectrum of an operator with nonnegative real part lies in the closed right half-plane. Section~\ref{sec:convex} uses families of both kinds.

\subsection{Changes of variable}\label{ssec:covariance}

\begin{lemma}\label{lem:covariance}
Let $\varphi:\Om\to\Om'$ be holomorphic, let $k'$ be a holomorphic kernel on $\Om'$, let $c>0$, and put $k(z,w)=c\,k'(\varphi(z),\varphi(w))$. Then $k$ is a holomorphic kernel on $\Om$, and for $T\in\B(H)$ with $\sigma(T)\subset\Om$,
\[
\Psi^{k}_T=c\,\Psi^{k'}_{\varphi(T)},
\]
where the superscript indicates the kernel. In particular, $\Psi^k_T$ is accretive if and only if $\Psi^{k'}_{\varphi(T)}$ is.
\end{lemma}

\begin{proof}
For an orthonormal basis $(e'_j)$ of $\E_{k'}$, Lemma~\ref{lem:conv}(a) gives $k(z,w)=c\sum_je'_j(\varphi(z))\overline{e'_j(\varphi(w))}$, and Lemma~\ref{lem:intkernel}, applied to the functions $e'_j\circ\varphi$, gives $\Psi^k_T=c\sum_j\Co{e'_j(\varphi(T))}=c\,\Psi^{k'}_{\varphi(T)}$.
\end{proof}

\section{Choosing the space for two disks}\label{sec:twodisks}

By Theorem~A, the complete constant of a class $\Kc$ is at most $\sqrt2\,M(\E)$ for every nondegenerate reproducing kernel Hilbert space $\E$ of holomorphic functions on $\Om$ whose moment maps are accretive throughout $\Kc$. This section proves Theorems~B and~C, choosing such spaces for the upper bounds. A M\"obius transformation reduces an intersection of two open disks of the Riemann sphere to an annulus, a strip or a sector (Section~\ref{ssec:normal}), which Sections~\ref{ssec:annulus}, \ref{ssec:strip} and~\ref{ssec:sector} treat separately. For the annulus and the strip, $\E$ is the space of a sum of two kernels, one for each disk, each of the form $k_r(z,w)=1/(r(z)+\overline{r(w)})$ or $\frac12(k_r+k_{1/r})$ for a Herglotz function $r$ (Section~\ref{ssec:cayley}), chosen to match the symmetry of the domain, and $M(\E)\le\sqrt2$. For the sector, $\E$ is the space with kernel $2\sin(2\alpha)\,z\bar w/(z^2-2\cos(2\alpha)\,z\bar w+\bar w^2)$, with $M(\E)\le\sqrt{2-4\alpha/\pi}$, when $0<\alpha\le\pi/4$, and the Hardy space $H^2(S_\alpha,dr/r)$, with $M(\E)=1$, when $\pi/4\le\alpha<\pi/2$; this gives the upper bound of Theorem~C, and Theorem~B for crossing circles.

\subsection{Reduction to three cases}\label{ssec:normal}
Let $D_1,D_2$ be open disks of the Riemann sphere such that $\Om=D_1\cap D_2$ is a nonempty subset of $\C$ different from $D_1$ and $D_2$. Let $b_i:D_i\to\D$ be conformal.

\begin{lemma}\label{lem:mobius}
There is a M\"obius transformation $m$ that maps $\Om$ onto $\Ann_R=\{R^{-1}<|z|<R\}$ for some $R>1$, onto $\Pi_a=\{|\Ip z|<a\}$ for some $a>0$, or onto $S_\alpha$ for some $0<\alpha<\pi/2$, according as the boundary circles of $D_1$ and $D_2$ are disjoint, tangent or crossing; in the last case $2\alpha$ is the angle of $\Om$ at the crossing points. Let $T\in\B(H)$ with $\sigma(T)\subset\Om$, and let $A=m(T)$. Then $F(T)=(F\circ m^{-1})(A)$ for every bounded holomorphic $F$ on $\Om$, and the conditions $\norm{b_1(T)}\le1$, $\norm{b_2(T)}\le1$ are equivalent, in the three cases, to $\norm A\le R$ and $\norm{A^{-1}}\le R$, to $\norm{\Ip A}\le a$, and to $W(A)\subseteq\Sbar\alpha$.
\end{lemma}

\begin{proof}
This is standard M\"obius geometry, see \cite[Section~4]{BBC}: for disjoint circles, see also \cite[proof of Theorem~10]{CG19}; otherwise, send a common point of the two circles to $\infty$, and translate, rotate and dilate. If $\norm{b_i(T)}\le1$ for one conformal map $b_i$ of $D_i$ onto $\D$, then this holds for every such map, since two such maps differ by an automorphism $\varphi$ of $\D$, and $\norm{\varphi(B)}\le1$ for every contraction $B$, by von Neumann's inequality \cite{vN}. As $b_i(T)=(b_i\circ m^{-1})(A)$, we may therefore take for $b_i\circ m^{-1}$ any conformal map of $m(D_i)$ onto $\D$: $b=z/R$ and $b=1/(Rz)$ for the two disks of the annulus, and $b=(r-1)/(r+1)$ for a half-plane $\{\Rp r>0\}$ with $r$ affine; then use \eqref{eq:cayleyfact}.
\end{proof}

For $\Ann_R$, $\Pi_a$ and $S_\alpha$ we use the \emph{defining families}
\begin{equation}\label{eq:families}
\Bf_{\Ann}=\Big\{\frac{R+z}{R-z},\ \frac{Rz+1}{Rz-1}\Big\},\qquad \Bf_\Pi=\{a-iz,\ a+iz\},\qquad \Bf_S=\{e^{-i\gamma}z,\ e^{i\gamma}z\},
\end{equation}
the first consisting of the Cayley transforms $(1+b)/(1-b)$ of $b=z/R$ and $b=1/(Rz)$, and the last with $\gamma=\pi/2-\alpha$. Their classes are
\begin{align*}
\Kc(\Bf_{\Ann})&=\{T:\ \sigma(T)\subset\Ann_R,\ \norm T\le R,\ \norm{T^{-1}}\le R\},\\
\Kc(\Bf_\Pi)&=\{T:\ \sigma(T)\subset\Pi_a,\ \norm{\Ip T}\le a\},\\
\Kc(\Bf_S)&=\{T:\ \sigma(T)\subset S_\alpha,\ W(T)\subseteq\Sbar\alpha\}.
\end{align*}
For $\Ann_R$ and $\Pi_a$, with $\Bf$ the defining family of the domain, $\E=\E_\Bf$ is the space of the kernel
\begin{equation}\label{eq:kB}
k_\Bf=\sum_{r\in\Bf}\tfrac12(k_r+k_{1/r})\quad\text{for }\Ann_R,\qquad k_\Bf=\sum_{r\in\Bf}k_r\quad\text{for }\Pi_a,
\end{equation}
where $k_r(z,w)=1/(r(z)+\overline{r(w)})$, so that $\Psi_T$ is accretive for every $T\in\Kc(\Bf)$ by Corollary~\ref{cor:sums}. The spaces for $S_\alpha$ are given in Section~\ref{ssec:sector}. In all three cases, we compare the norm of $\E$ with a Hardy norm $\norm\cdot_{\mathcal N}$ by Lemma~\ref{lem:compare}. Both norms are weighted $L^2$ norms of the coefficients in the expansions
\[
f(z)=\sum_{n\in\Z}c_nz^n,\qquad f(z)=\int_\RR\hat f(t)\,e^{itz}\,dt,\qquad f(z)=\int_\RR\hat f(\xi)\,z^{i\xi}\,d\xi,
\]
by \eqref{eq:intspace} for $\E$, and by the theorems of Parseval and Plancherel for $\norm\cdot_{\mathcal N}$. The comparison of the two norms thus becomes a comparison of weights.

\subsection{The quantum annulus}\label{ssec:annulus}
Let $R>1$. For $r=(1+b)/(1-b)$ we have $\frac12(k_r+k_{1/r})(z,w)=\frac12+\sum_{n\ge1}b(z)^n\overline{b(w)}^n$, by \eqref{eq:schurkernel}; the choices $b=z/R$ and $b=1/(Rz)$ give
\[
k_{\Bf_{\Ann}}(z,w)=\Big(\frac12+\sum_{n\ge1}\frac{(z\bar w)^n}{R^{2n}}\Big)+\Big(\frac12+\sum_{n\ge1}\frac{(z\bar w)^{-n}}{R^{2n}}\Big)=\sum_{n\in\Z}R^{-2|n|}(z\bar w)^n .
\]
By \eqref{eq:intspace}, with $\Lambda=\Z$, $\phi_n=z^n$ and $\beta(n)^2=R^{2|n|}$, $\E_{\Bf_{\Ann}}$ consists of the functions $f=\sum_{n\in\Z}c_nz^n$ with $\norm f_\E^2=\sum_nR^{2|n|}|c_n|^2<\infty$. This space appears, up to normalization, in \cite{ARS08}, and, up to the rescaling $z\mapsto z/R$, in \cite[Proposition~3.1]{Ts22}, where the bound $\sqrt2$ is proved for scalar multipliers and shown to be sharp. For the comparison we use the following Hardy norm, with equal weights on the two circles: $\norm f_{\mathcal N}^2$ is the supremum over $1<\tau<R$ of $\int_{|z|=\tau}|f|^2+\int_{|z|=1/\tau}|f|^2$, with normalized arclength. For $f=\sum_nc_nz^n$ holomorphic on $\Ann_R$, this sum of integrals is $\sum_n|c_n|^2(\tau^{2n}+\tau^{-2n})$, which increases, as $\tau\uparrow R$, to $\sum_n(R^{2n}+R^{-2n})|c_n|^2$.

\begin{proposition}\label{prop:annulus}
For every $T\in\Kc(\Bf_{\Ann})$, every $d\ge1$ and every bounded holomorphic $F:\Ann_R\to M_d$, $\norm{F(T)}\le2\sup_{\Ann_R}\norm F$.
\end{proposition}

\begin{proof}
By Corollary~\ref{cor:sums}, the moment map of $k_{\Bf_{\Ann}}$ at $T$ is accretive. Since $(R^{2n}+R^{-2n})/R^{2|n|}=1+R^{-4|n|}\in(1,2]$, Lemma~\ref{lem:compare} applies with $\norm\cdot_{\mathcal N}$, $c_-=1$ and $c_+=2$, so $M(\E_{\Bf_{\Ann}})\le\sqrt2$. The space contains $1$, so it is nondegenerate, and Corollary~\ref{cor:main} gives the bound $\sqrt2\cdot\sqrt2$.
\end{proof}

\begin{corollary}\label{cor:qa}
Let $R>1$ and let $T\in\B(H)$ satisfy $\norm T\le R$ and $\norm{T^{-1}}\le R$. For every $d\ge1$ and every $M_d$-valued rational function $F$ with poles off $\overline{\Ann_R}=\{R^{-1}\le|z|\le R\}$,
\[
\norm{F(T)}\le2\max_{\overline{\Ann_R}}\norm F .
\]
The constant $2$ is optimal for every $R>1$, even for scalar $F$, here and in Proposition~\ref{prop:annulus}: $K(\Kc(\Bf_{\Ann}))=K_{\mathrm{cb}}(\Kc(\Bf_{\Ann}))=2$.
\end{corollary}

\begin{proof}
Let $R'>R$ be so close to $R$ that $F$ has no poles in $\overline{\Ann_{R'}}$. Then $\sigma(T)\subseteq\overline{\Ann_R}\subset\Ann_{R'}$ and $\norm T,\norm{T^{-1}}\le R'$, so $T\in\Kc(\Bf_{\Ann})$ for $R'$ in place of $R$, and Proposition~\ref{prop:annulus} gives $\norm{F(T)}\le2\sup_{\Ann_{R'}}\norm F$. Let $R'\downarrow R$ and use the continuity of $F$ near $\overline{\Ann_R}$. By Proposition~\ref{prop:annulus}, $K(\Kc(\Bf_{\Ann}))\le K_{\mathrm{cb}}(\Kc(\Bf_{\Ann}))\le2$, and the lower bounds follow from \cite{Ts}: the weighted shifts used there are similar, by a diagonal operator, to the bilateral shift on $\ell^2(\Z)$, so their spectrum is $\T\subset\Ann_R$.
\end{proof}

\subsection{The strip}\label{ssec:strip}
Let $a>0$. Since $\big(2a\mp i(z-\bar w)\big)^{-1}=\int_0^\infty e^{-2at}e^{\pm itz}\,\overline{e^{\pm itw}}\,dt$,
\[
k_{\Bf_\Pi}(z,w)=\frac1{2a-i(z-\bar w)}+\frac1{2a+i(z-\bar w)}=\int_\RR e^{itz}\,\overline{e^{itw}}\,e^{-2a|t|}\,dt .
\]
By \eqref{eq:intspace}, with $\Lambda=\RR$, $\phi_t=e^{itz}$ and $\beta(t)^2=e^{2a|t|}$, $\E_{\Bf_\Pi}$ consists of the functions $f=\int_\RR\hat f(t)e^{itz}\,dt$ with $\norm f_\E^2=\int|\hat f(t)|^2e^{2a|t|}\,dt<\infty$. For the comparison we use the Hardy space $H^2(\Pi_a)$ of holomorphic $f$ on $\Pi_a$ with
\[
\norm f^2=\sup_{0<\vartheta<a}\int_\RR\big(|f(x+i\vartheta)|^2+|f(x-i\vartheta)|^2\big)\frac{dx}{2\pi}<\infty .
\]
By the Paley--Wiener theorem for the strip \cite[Theorem~IV]{PW}, it consists of the functions $f(z)=\int_\RR\hat f(t)e^{itz}\,dt$, with $\hat f(t)=\frac1{2\pi}\int_\RR f(x)e^{-itx}\,dx$ the Fourier transform of $f|_\RR$ in the $L^2$ sense, and the Fourier transform of $f(\cdot+iy)$ is $e^{-yt}\hat f(t)$ for $|y|<a$; so, by Plancherel's theorem on the lines $\Ip z=\pm\vartheta$, the integral in the definition of the norm is $\int|\hat f(t)|^2(e^{2\vartheta t}+e^{-2\vartheta t})\,dt$, which increases, as $\vartheta\uparrow a$, to $\norm f^2=\int|\hat f(t)|^2(e^{2at}+e^{-2at})\,dt$.

\begin{proposition}\label{prop:strip}
For every $T\in\Kc(\Bf_\Pi)$, every $d\ge1$ and every bounded holomorphic $F:\Pi_a\to M_d$, $\norm{F(T)}\le2\sup_{\Pi_a}\norm F$.
\end{proposition}

\begin{proof}
By Corollary~\ref{cor:sums}, the moment map of $k_{\Bf_\Pi}$ at $T$ is accretive. Since $(e^{2at}+e^{-2at})/e^{2a|t|}=1+e^{-4a|t|}\in(1,2]$, Lemma~\ref{lem:compare} applies with the norm of $H^2(\Pi_a)$, $c_-=1$ and $c_+=2$, so $M(\E_{\Bf_\Pi})\le\sqrt2$. The space is nondegenerate: $e^{-\epsilon z^2}\in\E_{\Bf_\Pi}$ for $\epsilon>0$, and it tends to $1$ locally uniformly as $\epsilon\downarrow0$. Corollary~\ref{cor:main} gives the bound $\sqrt2\cdot\sqrt2$.
\end{proof}

\begin{remark}\label{rem:stripannulus}
Proposition~\ref{prop:strip} also follows from Proposition~\ref{prop:annulus}. Let $T\in\Kc(\Bf_\Pi)$ and $n\ge1$. Since $\Rp(\pm iT)=\mp\Ip T\le aI$, the operator $U_n=e^{iT/n}$ satisfies $\norm{U_n},\norm{U_n^{-1}}\le e^{a/n}$, and $\sigma(U_n)=e^{i\sigma(T)/n}\subset\Ann_{e^{a/n}}$; so $U_n\in\Kc(\Bf_{\Ann})$ for $R=e^{a/n}$. Under $\zeta=e^{iz/n}$, the bounded holomorphic functions on $\Ann_{e^{a/n}}$ are the $2\pi n$-periodic ones on $\Pi_a$. For $\epsilon>0$, $F_n(z)=\sum_{k\in\Z}F(z+2\pi nk)\,e^{-\epsilon(z+2\pi nk)^2}$ is $2\pi n$-periodic, $\limsup_n\sup_{\Pi_a}\norm{F_n}\le e^{\epsilon a^2}\sup_{\Pi_a}\norm F$, and $F_n\to Fe^{-\epsilon z^2}$ uniformly near $\sigma(T)$. Proposition~\ref{prop:annulus} at $U_n$ gives $\norm{F_n(T)}\le2\sup_{\Pi_a}\norm{F_n}$; let $n\to\infty$, then $\epsilon\to0$.
\end{remark}

The lower bound for the strip is modelled on that of \cite{Ts} for the annulus. Under $\zeta=e^{iz}$, which maps $\Pi_a$ onto $\Ann_R$ with $R=e^a$, the functions $\cos(nz)/\cosh(na)$ become $(\zeta^n+\zeta^{-n})/(R^n+R^{-n})$, up to constant factors the functions of \cite{Ts}, where they are evaluated at bilateral weighted shifts whose weights equal $R$ at each of $n$ consecutive steps and $R^{-1}$ at each of the next $n$, periodically. The proof of Proposition~\ref{prop:stripopt} uses a continuous analogue of these shifts: the group of weighted translations generated by $-d/d\theta+V$ on $L^2(\RR/2n\Z)$, where $V$ equals $c$ and $-c$ on alternate half-periods.

\begin{proposition}\label{prop:stripopt}
For every $a>0$, the constant $2$ in Proposition~\ref{prop:strip} is optimal, even for scalar $F$: $K(\Kc(\Bf_\Pi))=K_{\mathrm{cb}}(\Kc(\Bf_\Pi))=2$.
\end{proposition}

\begin{proof}
By Proposition~\ref{prop:strip}, $K(\Kc(\Bf_\Pi))\le K_{\mathrm{cb}}(\Kc(\Bf_\Pi))\le2$. For the lower bound, fix $0<c<a$ and $n\ge1$. On $H=L^2(\RR/2n\Z)$ let $D=-i\,d/d\theta$, self-adjoint with domain the periodic Sobolev space $W^{1,2}(\RR/2n\Z)$, let $V$ be multiplication by the function equal to $c$ on $(0,n)$ and to $-c$ on $(n,2n)$, and put $L=-iD+V=-d/d\theta+V$. A direct computation shows that $L$ generates the group of weighted translations
\[
(e^{tL}g)(\theta)=\exp\Big(\int_{\theta-t}^{\theta}V\Big)\,g(\theta-t)\qquad(t\in\RR).
\]
Since $V$ has integral $0$ over a period, $e^{2nL}=I$, so $e^{-nL}=e^{nL}$; and $\norm{e^{nL}}=\max_\theta\exp\int_{\theta-n}^{\theta}V=e^{cn}$, attained at $\theta=n$.

For $m\ge1$ let $T_m=D\,\mathbf 1_{[-m,m]}(D)+iV\in\B(H)$. Then $\Ip T_m=V$, so $\norm{\Ip T_m}=c$ and $\sigma(T_m)\subseteq\overline{W(T_m)}\subseteq\{|\Ip z|\le c\}\subset\Pi_a$; thus $T_m\in\Kc(\Bf_\Pi)$. Put $L_m=-iT_m=-iD\,\mathbf 1_{[-m,m]}(D)+V$. Since $\frac{d}{dt}\norm{e^{tL_m}x}^2=2\ip{e^{tL_m}x}{Ve^{tL_m}x}$, we have $\norm{e^{tL_m}}\le e^{c|t|}$. For $x\in W^{1,2}(\RR/2n\Z)$ and $t>0$,
\[
e^{tL_m}x-e^{tL}x=\int_0^te^{(t-s)L_m}(L_m-L)e^{sL}x\,ds,\qquad (L_m-L)y=iD\,\mathbf 1_{\RR\setminus[-m,m]}(D)\,y,
\]
and $\norm{D\,\mathbf 1_{\RR\setminus[-m,m]}(D)\,y}\to0$ as $m\to\infty$, uniformly for $y$ in the set $\{e^{sL}x:0\le s\le t\}$, which is compact in $W^{1,2}(\RR/2n\Z)$, since, as $V$ is bounded, the graph norm of $L$ is equivalent to the norm of $W^{1,2}(\RR/2n\Z)$, and $s\mapsto e^{sL}x$ is continuous in the graph norm. Hence $e^{tL_m}x\to e^{tL}x$ for every $x\in H$, by density, since $\norm{e^{tL_m}}$ and $\norm{e^{tL}}$ are at most $e^{c|t|}$; likewise for $t<0$.

The function $F(z)=\cos(nz)/\cosh(na)$ is bounded and holomorphic on $\Pi_a$, with $|F|\le1$ there, since $|\cos(x+iy)|^2=\cos^2x+\sinh^2y\le\cosh^2y$. As $T_m=iL_m$,
\[
F(T_m)=\frac{e^{inT_m}+e^{-inT_m}}{2\cosh(na)}=\frac{e^{-nL_m}+e^{nL_m}}{2\cosh(na)}\longrightarrow\frac{e^{-nL}+e^{nL}}{2\cosh(na)}=\frac{e^{nL}}{\cosh(na)}
\]
applied to each $x\in H$, as $m\to\infty$. Since $\norm{F(T_m)}\le K(\Kc(\Bf_\Pi))$, this gives $\norm{e^{nL}x}\le K(\Kc(\Bf_\Pi))\cosh(na)\norm x$ for every $x$, so $K(\Kc(\Bf_\Pi))\ge\norm{e^{nL}}/\cosh(na)=e^{cn}/\cosh(na)$. Letting $c\uparrow a$ gives $K(\Kc(\Bf_\Pi))\ge2/(1+e^{-2an})$, and letting $n\to\infty$ gives $K(\Kc(\Bf_\Pi))\ge2$.
\end{proof}

\subsection{Sectors, lenses and crescents}\label{ssec:sector}
Let $0<\alpha<\pi/2$ and $\gamma=\pi/2-\alpha$. By Lemma~\ref{lem:mobius}, lenses and crescents of vertex angle $2\alpha$ reduce to the class $\Kc(\Bf_S)$. For $\alpha\le\pi/4$, we use a kernel built from the two functions $e^{\mp i\gamma}z$ of $\Bf_S$ in \eqref{eq:families}: its moment maps are accretive on $\Kc(\Bf_S)$ (Lemma~\ref{lem:inverse}, Proposition~\ref{prop:rightangle}), and a comparison of weights bounds its multiplier constant by $\sqrt{2-4\alpha/\pi}$ (Lemma~\ref{lem:prodweights}). For $\alpha\ge\pi/4$, we use the Hardy space $H^2(S_\alpha,dr/r)$, with multiplier constant $1$, whose moment maps are accretive on $\Kc(\Bf_S)$ by a change of variable from the right angle (Lemma~\ref{lem:fracpow}, Proposition~\ref{prop:obtuse}). Proposition~\ref{prop:sectorbound} combines the two cases into the complete bound $c_\alpha=\max\{\sqrt2,\sqrt{4-8\alpha/\pi}\}$ on $\Kc(\Bf_S)$; since $c_\alpha<2$, this gives the case of sectors in Theorem~B. Proposition~\ref{prop:spaceopt} shows that, for $\alpha<\pi/4$, no other space gives less in Theorem~A, and Proposition~\ref{prop:sectorlower} gives the lower bound $2^{1-2\alpha/\pi}$. Corollaries~\ref{cor:sectorclosed}--\ref{cor:unbounded} extend the bounds to the closed sector, to lenses and crescents, and to unbounded operators. We begin with $H^2(S_\alpha,dr/r)$, which is the space for $\alpha\ge\pi/4$ and provides the comparison norm $\norm\cdot_{\mathcal N}$ of Lemma~\ref{lem:compare} for $\alpha\le\pi/4$.

\begin{lemma}\label{lem:hardysector}
\begin{enumerate}[(i)]
\item Let $H^2(S_\alpha,dr/r)$ be the space of holomorphic $f$ on $S_\alpha$ with
\[
\norm f^2=\sup_{0<\vartheta<\alpha}\int_0^\infty\big(|f(re^{i\vartheta})|^2+|f(re^{-i\vartheta})|^2\big)\frac{dr}{2\pi r}<\infty.
\]
The map $f\mapsto f\circ\exp$ is unitary from $H^2(S_\alpha,dr/r)$ onto the space $H^2(\Pi_\alpha)$ of Section~\ref{ssec:strip}; more precisely, for each $\vartheta$ the integral over the rays $\arg z=\pm\vartheta$ equals the integral of the image over the lines $\Ip w=\pm\vartheta$.
\item Consequently $H^2(S_\alpha,dr/r)$ consists of the functions $f=\int\hat f(\xi)z^{i\xi}\,d\xi$ with $\norm f^2=\int|\hat f(\xi)|^2(e^{2\alpha\xi}+e^{-2\alpha\xi})\,d\xi$. The reproducing kernel of $H^2(S_\alpha,dr/r)$ is
\[
k_{S_\alpha}(z,w)=\int_\RR z^{i\xi}\,\overline{w^{i\xi}}\,\frac{d\xi}{2\cosh2\alpha\xi}=\frac{2s\,z^{s}\,\bar w^{s}}{z^{2s}+\bar w^{2s}},\qquad s=\frac\pi{4\alpha}.
\]
\end{enumerate}
\end{lemma}

\begin{proof}
(i) This is the change of variable $w=\log z$, as for Hardy spaces of general domains \cite[Chapter~10]{Du70}: it sends $S_\alpha$ onto $\Pi_\alpha$ and the ray $\arg z=\pm\vartheta$ onto the line $\Ip w=\pm\vartheta$, with $dr/r=dx$ for $w=x\pm i\vartheta$.

(ii) The description follows from (i) and the description of $H^2(\Pi_\alpha)$ in Section~\ref{ssec:strip}, since $e^{i\xi\log z}=z^{i\xi}$. So \eqref{eq:intspace}, with $\beta(\xi)^2=2\cosh2\alpha\xi$, gives the integral form of the kernel. For the closed form, $z^{i\xi}\overline{w^{i\xi}}=e^{i\xi u}$ with $u=\log z-\log\bar w$. Since $\int_\RR e^{-2\pi ix\xi}\operatorname{sech}(\pi\xi)\,d\xi=\operatorname{sech}(\pi x)$ for real $x$,
\begin{equation}\label{eq:sech}
\int_\RR e^{i\xi u}\operatorname{sech}(\pi\xi)\,d\xi=\operatorname{sech}(u/2)
\end{equation}
for real $u$, and hence for $|\Ip u|<\pi$, by analytic continuation. The substitution $\xi\mapsto2s\xi$ gives $\int e^{i\xi u}\,d\xi/(2\cosh2\alpha\xi)=s\operatorname{sech}(su)=2s/(e^{su}+e^{-su})$ for $|\Ip u|<2\alpha$, and here $\Ip u=\arg z+\arg w\in(-2\alpha,2\alpha)$. Finally $e^{su}=z^{s}\bar w^{-s}$, and $z^{2s}+\bar w^{2s}\ne0$ because both terms lie in the open right half-plane.
\end{proof}

For $0<\alpha\le\pi/4$ we use, instead of the sum $k_{r_+}+k_{r_-}$ of the kernels of the two functions $r_\pm(z)=e^{\mp i\gamma}z$ of $\Bf_S$, their product times $2\sin(2\alpha)z\bar w$. Since $\cos2\gamma=-\cos2\alpha$, we have $(e^{-i\gamma}z+e^{i\gamma}\bar w)(e^{i\gamma}z+e^{-i\gamma}\bar w)=z^2-2\cos(2\alpha)z\bar w+\bar w^2$, and we put
\[
\tilde k_{S_\alpha}(z,w)=2\sin(2\alpha)\,z\bar w\,k_{r_+}(z,w)\,k_{r_-}(z,w)=\frac{2\sin(2\alpha)\,z\bar w}{z^2-2\cos(2\alpha)\,z\bar w+\bar w^2}.
\]
For $\alpha=\pi/4$ this is the kernel $k_{S_{\pi/4}}$ of $H^2(S_{\pi/4},dr/r)$, by Lemma~\ref{lem:hardysector}(ii) with $s=1$. The moment map of $\tilde k_{S_\alpha}$ can be inverted explicitly, and its accretivity rests on the following inequality.

\begin{lemma}\label{lem:inverse}
Let $0<\theta\le\pi/2$, and let $\mathcal A$ be an invertible bounded operator on $\Sc(H)$ with $\ip X{\mathcal A(X)}_{\mathrm{HS}}\in\Sbar\theta$ for every $X\in\Sc(H)$. Then $\Rp(\mathcal A+\mathcal A^{-1})\ge2\cos\theta\,I$.
\end{lemma}

\begin{proof}
For $\epsilon>0$, $\mathcal A+\epsilon I$ satisfies the hypothesis, since $\Sbar\theta+[0,\infty)\subseteq\Sbar\theta$, and it is invertible, since $\Rp(\mathcal A+\epsilon I)\ge\epsilon I$; as $(\mathcal A+\epsilon I)^{-1}\to\mathcal A^{-1}$ when $\epsilon\to0$, it suffices to prove the claim when $G=\Rp\mathcal A\ge\epsilon I$ for some $\epsilon>0$. Put $K=(\mathcal A-\mathcal A^\dagger)/(2i)$ and $L=G^{-1/2}KG^{-1/2}$, so that $\mathcal A=G^{1/2}(I+iL)G^{1/2}$. The hypothesis says that $\cos\theta\,|\ip X{K(X)}_{\mathrm{HS}}|\le\sin\theta\,\ip X{G(X)}_{\mathrm{HS}}$, so $\cos\theta\,\norm L\le\sin\theta$, that is, $\cos^2\theta\,(I+L^2)\le I$, and
\[
\Rp(\mathcal A^{-1})=G^{-1/2}\,\Rp\big((I+iL)^{-1}\big)\,G^{-1/2}=G^{-1/2}(I+L^2)^{-1}G^{-1/2}\ge\cos^2\theta\,G^{-1}.
\]
Hence $\Rp(\mathcal A+\mathcal A^{-1})\ge G+\cos^2\theta\,G^{-1}\ge2\cos\theta\,I$, since $t+\cos^2\theta/t\ge2\cos\theta$ for $t>0$.
\end{proof}

For matrices and $\theta<\pi/2$, the inequality $\Rp(\mathcal A^{-1})\ge\cos^2\theta\,(\Rp\mathcal A)^{-1}$ in this proof is \cite[Lemma~3]{Li16}, where it is noted to be implicit in the proof of \cite[Theorem~3.1]{DL14}.

\begin{proposition}\label{prop:rightangle}
Let $0<\alpha\le\pi/4$, let $A\in\B(H)$ satisfy $\sigma(A)\subset S_\alpha$ and $W(A)\subseteq\Sbar\alpha$, and let $\Psi_A$ be the moment map of $\tilde k_{S_\alpha}$. Then $A$ and $\Psi_A$ are invertible,
\[
\Psi_A^{-1}=\frac{\mathcal M+\mathcal M^{-1}-2\cos(2\alpha)\,I}{2\sin2\alpha},\qquad \mathcal M(X)=AX(A^{-1})^*\quad(X\in\Sc(H)),
\]
and $\Psi_A$ is accretive.
\end{proposition}

\begin{proof}
Put $R_\vartheta=\Rp(e^{i\vartheta}A)$. Since $\alpha+\gamma=\pi/2$, $e^{i\vartheta}\Sbar\alpha$ lies in the closed right half-plane for $|\vartheta|\le\gamma$, so $W(A)\subseteq\Sbar\alpha$ gives $R_\vartheta\ge0$ for these $\vartheta$, in particular for $\vartheta=\pm\alpha$ and $\vartheta=\pm\gamma$, as $\alpha\le\gamma$. And $A$ is invertible, since $0\notin\sigma(A)$.

The kernel $\tilde k_{S_\alpha}$ is a holomorphic kernel, as a product of holomorphic kernels, and $(z^2-2\cos(2\alpha)z\bar w+\bar w^2)\,\tilde k_{S_\alpha}(z,w)=2\sin(2\alpha)z\bar w$. By \eqref{eq:fXf} and polarization, as in the proof of Lemma~\ref{lem:cauchy}, multiplying the integrand of \eqref{eq:cauchy}, for $T=A$, by $a(z)\overline{b(w)}$, for $a,b$ holomorphic on $S_\alpha$, replaces $\Psi_A(X)$ by $a(A)\Psi_A(X)b(A)^*$, and replacing $k$ by $a(z)\overline{b(w)}$ gives $a(A)Xb(A)^*$. Hence, for $X\in\Sc(H)$,
\[
A^2\Psi_A(X)-2\cos(2\alpha)A\Psi_A(X)A^*+\Psi_A(X)(A^*)^2=2\sin(2\alpha)AXA^*.
\]
Since $\cos2\gamma=-\cos2\alpha$, the left side is $\mathcal L_-\mathcal L_+(\Psi_A(X))$, where $\mathcal L_\pm(Y)=r_\pm(A)Y+Yr_\pm(A)^*$ is invertible by Proposition~\ref{prop:herglotz}, as $\Rp r_\pm>0$ on $S_\alpha$. As $A$ is invertible, so is $\Psi_A$, and applying $X\mapsto A^{-1}X(A^{-1})^*$ to both sides gives the formula for $\Psi_A^{-1}$, since $\mathcal M^{-1}(X)=A^{-1}XA^*$.

We show that $\ip X{\mathcal M(X)}_{\mathrm{HS}}\in\Sbar{2\alpha}$ for $X\in\Sc(H)$. Put $Y=X(A^*)^{-1}$; then $X=YA^*$, so $\ip{X}{\mathcal M(X)}_{\mathrm{HS}}=\tr(AY^*AY)$. For self-adjoint $P,Q$ the number $\tr(PY^*QY)$ is real, since its conjugate is $\tr(Y^*QYP)$. Multiplying out $\tr(e^{i\gamma}AY^*e^{-i\alpha}AY)$ with $e^{i\gamma}A=R_\gamma+iR_{-\alpha}$ and $e^{-i\alpha}A=R_{-\alpha}-iR_\gamma$, and $\tr(e^{i\alpha}AY^*e^{-i\gamma}AY)$ with $e^{i\alpha}A=R_\alpha+iR_{-\gamma}$ and $e^{-i\gamma}A=R_{-\gamma}-iR_\alpha$, and taking real parts, we obtain, as $\gamma-\alpha=\pi/2-2\alpha$,
\[
\Rp\big(e^{\pm i(\pi/2-2\alpha)}\ip{X}{\mathcal M(X)}_{\mathrm{HS}}\big)=\hs{R_{\mp\alpha}^{1/2}YR_{\pm\gamma}^{1/2}}^2+\hs{R_{\pm\gamma}^{1/2}YR_{\mp\alpha}^{1/2}}^2\ge0,
\]
that is, $\ip X{\mathcal M(X)}_{\mathrm{HS}}\in\Sbar{2\alpha}$. Lemma~\ref{lem:inverse}, with $\theta=2\alpha$, gives $\Rp(\mathcal M+\mathcal M^{-1})\ge2\cos(2\alpha)I$. So $\Psi_A^{-1}$ is accretive, and so is its inverse $\Psi_A$.
\end{proof}

The factor $2z\bar w$ is essential, already for $\alpha=\pi/4$, where $\tilde k_{S_{\pi/4}}=2z\bar w\,k_r$ with $r(z)=z^2$. By Proposition~\ref{prop:herglotz}, the moment map of $k_r$ at $A$ is the inverse of $\mathcal L(X)=A^2X+X(A^2)^*$, so it is accretive only if $\mathcal L$ is, which, tested on $X=xx^*$, requires $\Rp A^2\ge0$. This fails for $A_0=I+\sqrt2N$, with $N=\big(\begin{smallmatrix}0&1\\0&0\end{smallmatrix}\big)$: its numerical range is the disk $\{|z-1|\le\sqrt2/2\}\subseteq\Sbar{\pi/4}$, while $\Rp A_0^2=I+\sqrt2(N+N^*)$ has the eigenvalue $1-\sqrt2$.

\begin{lemma}\label{lem:prodweights}
Let $0<\alpha\le\pi/4$, and let $\nu_\alpha(\xi)=\sinh((\pi-2\alpha)\xi)/\sinh(\pi\xi)$ for $\xi\ne0$ and $\nu_\alpha(0)=1-2\alpha/\pi$. Then $\E_{\tilde k_{S_\alpha}}$ consists of the functions $f=\int_\RR\hat f(\xi)z^{i\xi}\,d\xi$ with $\norm f^2=\int|\hat f(\xi)|^2\nu_\alpha(\xi)^{-1}\,d\xi<\infty$, and, for $\xi\in\RR$,
\[
\nu_\alpha(\xi)^{-1}\le2\cosh(2\alpha\xi)\le\big(2-4\alpha/\pi\big)\,\nu_\alpha(\xi)^{-1}.
\]
\end{lemma}

\begin{proof}
For $z,w\in S_\alpha$ put $u=\log z-\log\bar w$; then $|\Ip u|<2\alpha$, $e^{i\xi u}=z^{i\xi}\,\overline{w^{i\xi}}$ and $\tilde k_{S_\alpha}(z,w)=\sin(2\alpha)/(\cosh u-\cos2\alpha)$. By \cite[formula~3.981(5)]{GR}, with $\pi-2\alpha$ and $\pi$ in place of $a$ and $b$, and since $\nu_\alpha$ is even,
\[
\int_\RR e^{i\xi u}\nu_\alpha(\xi)\,d\xi=2\int_0^\infty\nu_\alpha(\xi)\cos(u\xi)\,d\xi=\frac{\sin2\alpha}{\cosh u-\cos2\alpha}\qquad(u\in\RR).
\]
Both sides are holomorphic for $|\Ip u|<2\alpha$, since $\nu_\alpha(\xi)=O(e^{-2\alpha|\xi|})$ and $\cosh u\ne\cos2\alpha$ there, so the identity holds for these $u$. Hence $\tilde k_{S_\alpha}(z,w)=\int z^{i\xi}\,\overline{w^{i\xi}}\,\nu_\alpha(\xi)\,d\xi$, and \eqref{eq:intspace}, with $\beta^2=\nu_\alpha^{-1}$, gives the description of $\E_{\tilde k_{S_\alpha}}$. Finally, $2\cosh(2\alpha\xi)\sinh((\pi-2\alpha)\xi)=\sinh(\pi\xi)+\sinh((\pi-4\alpha)\xi)$, so
\[
2\cosh(2\alpha\xi)\,\nu_\alpha(\xi)=1+\frac{\sinh((\pi-4\alpha)\xi)}{\sinh(\pi\xi)}\in\Big[1,\,2-\frac{4\alpha}\pi\Big],
\]
because the last quotient is even and, with $c=1-4\alpha/\pi\in[0,1)$, $0\le\sinh(c\pi\xi)\le c\sinh(\pi\xi)$ for $\xi\ge0$, by the convexity of $\sinh$ on $[0,\infty)$.
\end{proof}

For $\pi/4\le\alpha<\pi/2$ we use the kernel $k_{S_\alpha}$ of $H^2(S_\alpha,dr/r)$ (Lemma~\ref{lem:hardysector}(ii)), and we reduce to $\alpha=\pi/4$ by the change of variable $z\mapsto z^s$, $s=\pi/(4\alpha)\in(\frac12,1]$, using Kato's theorem on fractional powers. Since $(z^s)^2=z^{2s}$ and $\overline{w^s}=\bar w^s$ on $S_\alpha$,
\begin{equation}\label{eq:kpower}
k_{S_\alpha}(z,w)=s\,k_{S_{\pi/4}}(z^s,w^s).
\end{equation}

\begin{lemma}\label{lem:fracpow}
Let $0<\alpha<\pi/2$ and $0<s\le1$. If $A\in\B(H)$ satisfies $\sigma(A)\subset S_\alpha$ and $W(A)\subseteq\Sbar\alpha$, then $W(A^s)\subseteq\Sbar{s\alpha}$.
\end{lemma}

\begin{proof}
The case $s=1$ is trivial; let $s<1$. Put $B_\pm=e^{\pm i\gamma}A$. Then $W(B_\pm)$ lies in the closed right half-plane, that is, $B_\pm$ is accretive, and $\sigma(B_\pm)$ lies in the open one, so $W(B_\pm^s)\subseteq\Sbar{s\pi/2}$ by Kato's theorem \cite[Theorem~2.4]{Ka}, see \cite[Corollary~7.1.13]{Ha06}. On $S_\alpha$ we have $|\arg z\pm\gamma|<\pi/2$, hence $(e^{\pm i\gamma}z)^s=e^{\pm is\gamma}z^s$ and $B_\pm^s=e^{\pm is\gamma}A^s$. Therefore
\[
W(A^s)\subseteq e^{-is\gamma}\,\Sbar{s\pi/2}\cap e^{is\gamma}\,\Sbar{s\pi/2}=\Sbar{s\pi/2-s\gamma}=\Sbar{s\alpha}.
\] This is the rotation argument of \cite[proof of Corollary~7.1.14]{Ha06}, where it is applied to the logarithm; for matrices with positive definite real part see also \cite[Lemma~2.2]{AS}, where that lemma is attributed to Kato and to \cite[Corollary~2.4]{Dr}.
\end{proof}

\begin{proposition}\label{prop:obtuse}
Let $\pi/4\le\alpha<\pi/2$ and $s=\pi/(4\alpha)$. For every $A\in\B(H)$ with $\sigma(A)\subset S_\alpha$ and $W(A)\subseteq\Sbar\alpha$, the moment map of $k_{S_\alpha}$ at $A$ is $s$ times that of $k_{S_{\pi/4}}$ at $A^s$, and it is accretive.
\end{proposition}

\begin{proof}
The map $z\mapsto z^s$ sends $S_\alpha$ onto $S_{\pi/4}$, so \eqref{eq:kpower} and Lemma~\ref{lem:covariance} give the first assertion, with $\sigma(A^s)\subset S_{\pi/4}$. By Lemma~\ref{lem:fracpow}, $W(A^s)\subseteq\Sbar{s\alpha}=\Sbar{\pi/4}$, and Proposition~\ref{prop:rightangle}, with $\alpha=\pi/4$, applies to $A^s$.
\end{proof}

For $0<\alpha<\pi/2$ let $c_\alpha$ be as in Theorem~C; thus $c_\alpha=\sqrt{4-8\alpha/\pi}$ if $\alpha\le\pi/4$, and $c_\alpha=\sqrt2$ if $\alpha\ge\pi/4$.

\begin{proposition}\label{prop:sectorbound}
Let $0<\alpha<\pi/2$. For every $A\in\B(H)$ with $\sigma(A)\subset S_\alpha$ and $W(A)\subseteq\Sbar\alpha$, every $d\ge1$ and every bounded holomorphic $F:S_\alpha\to M_d$,
\[
\norm{F(A)}\le c_\alpha\sup_{S_\alpha}\norm F .
\]
\end{proposition}

\begin{proof}
Let $\E$ be the space of $\tilde k_{S_\alpha}$ if $\alpha\le\pi/4$ and $H^2(S_\alpha,dr/r)$ if $\alpha\ge\pi/4$; the two coincide for $\alpha=\pi/4$. By Propositions~\ref{prop:rightangle} and~\ref{prop:obtuse}, $\Psi_A$ is accretive. By Lemmas~\ref{lem:hardysector} and~\ref{lem:prodweights}, $\E$ and $H^2(S_\alpha,dr/r)$ have the same elements, and the integrals over the rays $\arg z=\pm\vartheta$ in the norm of $H^2(S_\alpha,dr/r)$ increase with $\vartheta$, by Lemma~\ref{lem:hardysector}(i) and Section~\ref{ssec:strip}. So Lemma~\ref{lem:compare} applies to $\E$ with $\norm\cdot_{\mathcal N}$ the norm of $H^2(S_\alpha,dr/r)$, and with $c_-=1$ and $c_+=2-4\alpha/\pi$ if $\alpha\le\pi/4$; if $\alpha\ge\pi/4$, then $\E=H^2(S_\alpha,dr/r)$ and $c_-=c_+=1$. So $M(\E)\le c_\alpha/\sqrt2$. The space is nondegenerate: for $\epsilon>0$, the function $q_\epsilon(z)=z/((1+\epsilon z)(z+\epsilon))$ satisfies
\begin{equation}\label{eq:qeps}
|q_\epsilon(z)|^2\le\frac{|z|^2}{(1+\epsilon^2|z|^2)(|z|^2+\epsilon^2)}\qquad(\Rp z\ge0),
\end{equation}
since $|1+\epsilon z|^2\ge1+\epsilon^2|z|^2$ and $|z+\epsilon|^2\ge|z|^2+\epsilon^2$ for $\Rp z\ge0$; so $q_\epsilon$ lies in $H^2(S_\alpha,dr/r)$, and hence in $\E$, and $q_\epsilon\to1$ locally uniformly as $\epsilon\downarrow0$. Corollary~\ref{cor:main} gives the bound.
\end{proof}

\begin{proof}[Proof of Theorem~B]
By Lemma~\ref{lem:mobius}, we may assume that $\Om$ is $\Ann_R$, $\Pi_a$ or $S_\alpha$, and that $T\in\Kc(\Bf)$ for its defining family $\Bf$. These cases are Propositions~\ref{prop:annulus} and~\ref{prop:strip} and, since $c_\alpha<2$, Proposition~\ref{prop:sectorbound}.
\end{proof}

By the following proposition, for $0<\alpha<\pi/4$, no space gives less than $\sqrt{4-8\alpha/\pi}$ in Theorem~A, already when accretivity is required only at $2\times2$ matrices. Since the space of $\tilde k_{S_\alpha}$ gives this value (proof of Proposition~\ref{prop:sectorbound}), the bound $\sqrt2\inf M(\E)$ of Corollary~\ref{cor:main} for $\Kc(\Bf_S)$ equals $\sqrt{4-8\alpha/\pi}$.

\begin{proposition}\label{prop:spaceopt}
Let $0<\alpha<\pi/4$, and let $k\ne0$ be a holomorphic kernel on $S_\alpha$ whose moment map is accretive at every $A\in M_2$ with $\sigma(A)\subset S_\alpha$ and $W(A)\subseteq\Sbar\alpha$. Then $M(\E_k)\ge\sqrt{2-4\alpha/\pi}$.
\end{proposition}

We prove Proposition~\ref{prop:spaceopt} in Appendix~\ref{app:spaceopt}.

The proof of Proposition~\ref{prop:sectorlower} uses the $n\times n$ matrices similar, by the square root of a positive definite Toeplitz matrix, to $\mathrm{diag}(1,e^h,\dots,e^{(n-1)h})$, where the symbol of the Toeplitz matrix is chosen so that their numerical ranges lie in $\Sbar\alpha$, and the function $\cos(\pi\log z/h)/\cosh(\pi\alpha/h)$, whose values at the eigenvalues alternate in sign; letting $n\to\infty$ and then $h\to0$ gives the lower bound $2^{1-2\alpha/\pi}$ of Theorem~C.

\begin{proposition}\label{prop:sectorlower}
Let $0<\alpha<\pi/2$ and $\epsilon>0$. There are a matrix $A$ with $\sigma(A)\subset S_\alpha$ and $W(A)\subseteq\Sbar\alpha$ and a rational function $f$ with poles off $\Sbar\alpha$ such that $|f|\le1$ on $\Sbar\alpha$ and $\norm{f(A)}>2^{1-2\alpha/\pi}-\epsilon$. In particular, $K(\Kc(\Bf_S))\ge2^{1-2\alpha/\pi}$.
\end{proposition}

\begin{proof}
Fix $h>0$ and $a>h/2$, and let
\[
\varphi(z)=\frac{4\alpha}\pi\sum_{k\ \mathrm{odd}}\frac{\cos(kz)}{k\sinh(ka)},
\]
the sum over odd $k\ge1$, which converges locally uniformly on the strip $|\Ip z|<a$. Then $\varphi$ is holomorphic there, $2\pi$-periodic and real on $\RR$, $\varphi(\pi)=-\varphi(0)$, and $|\Ip\varphi|\le\alpha$: since $\frac4\pi\sum_{k\ \mathrm{odd}}\sin(kx)/k=\operatorname{sgn}(\sin x)$, the function $\Ip\varphi(x+iy)=-\frac{4\alpha}\pi\sum_{k\ \mathrm{odd}}\sin(kx)\sinh(ky)/(k\sinh(ka))$ is the Poisson integral on the strip of the function equal to $\mp\alpha\operatorname{sgn}(\sin x)$ on the lines $y=\pm a$.

Put $g=e^\varphi$, with Fourier coefficients $g_j=\frac1{2\pi}\int_0^{2\pi}g(x)e^{-ijx}\,dx$, so that $g(z)=\sum_jg_je^{ijz}$ for $|\Ip z|<a$. For $n\ge1$ let $P=(g_{j-l})_{j,l=0}^{n-1}$, $D=\mathrm{diag}(e^{jh})_{j=0}^{n-1}$ and $J=\mathrm{diag}((-1)^j)_{j=0}^{n-1}$, and for $u\in\C^n$ put $\hat u(x)=\sum_lu_le^{ilx}$. Since $g_je^{jh/2}$ are the Fourier coefficients of $g(x-ih/2)$,
\begin{gather*}
\ip u{Pu}=\frac1{2\pi}\int_0^{2\pi}g(x)|\hat u(x)|^2\,dx,\\
\ip u{D^{1/2}PD^{-1/2}u}=\frac1{2\pi}\int_0^{2\pi}g(x-ih/2)|\hat u(x)|^2\,dx.
\end{gather*}
As $g>0$ on $\RR$, $P$ is positive definite. Let $A=P^{-1/2}DP^{1/2}$. Then $\sigma(A)=\{e^{jh}:0\le j<n\}\subset S_\alpha$, and for $v\in\C^n$ and $u=D^{1/2}P^{-1/2}v$ we have $\ip v{Av}=\ip u{D^{1/2}PD^{-1/2}u}$, which lies in the closed convex cone $\Sbar\alpha$, since $|\arg g(x-ih/2)|=|\Ip\varphi(x-ih/2)|\le\alpha$. So $W(A)\subseteq\Sbar\alpha$.

The function
\[
F_h(z)=\frac{\cos(\pi\log z/h)}{\cosh(\pi\alpha/h)}
\]
is holomorphic on $S_\alpha$, with $|F_h|\le1$ there, since $|\cos(s+it)|\le\cosh t$. As $F_h(e^{jh})=(-1)^j/\cosh(\pi\alpha/h)$, we have $F_h(A)=P^{-1/2}JP^{1/2}/\cosh(\pi\alpha/h)$, whose adjoint maps $P^{1/2}w$ to $P^{1/2}Jw/\cosh(\pi\alpha/h)$ for $w\in\C^n$. For $w=((-1)^j)_{j=0}^{n-1}$ all entries of $Jw$ equal $1$, so
\[
\cosh^2(\pi\alpha/h)\,\norm{F_h(A)}^2\ge\frac{\ip{Jw}{PJw}}{\ip w{Pw}}=\frac{\int_0^{2\pi}g(x)\,\big|\sum_{j=0}^{n-1}e^{ijx}\big|^2\,dx}{\int_0^{2\pi}g(x)\,\big|\sum_{j=0}^{n-1}e^{ij(x+\pi)}\big|^2\,dx}.
\]
By Fej\'er's theorem, applied to $g$ at $0$ and at $\pi$, the right side tends to $g(0)/g(\pi)=e^{2\varphi(0)}$ as $n\to\infty$, and $\varphi(0)=\frac{4\alpha}\pi\sum_{k\ \mathrm{odd}}1/(k\sinh(ka))$. So, letting $n\to\infty$ and then $a\downarrow h/2$, the supremum $m_h$ of $\norm{F_h(A)}$ over $n\ge1$ and $a>h/2$ satisfies
\[
m_h\ge\frac{e^{(4\alpha/\pi)S(h)}}{\cosh(\pi\alpha/h)},\qquad S(h)=\sum_{k\ \mathrm{odd}}\frac1{k\sinh(kh/2)}.
\]

Let $Q(x)=1/(x\sinh(x/2))-2/x^2$, which extends continuously to $[0,\infty)$ and is $O(x^{-2})$ at infinity. Since $\sum_{k\ \mathrm{odd}}1/k^2=\pi^2/8$, Riemann sums give
\[
S(h)-\frac{\pi^2}{4h}=h\sum_{k\ \mathrm{odd}}Q(kh)\longrightarrow\frac12\int_0^\infty Q(x)\,dx\qquad(h\downarrow0).
\]
With $x=2y$, $\int_0^\infty Q(x)\,dx=\int_0^\infty\big(1/\sinh y-1/y\big)\,dy/y=-\log2$ \cite[formula~3.529(1)]{GR}. Hence, as $h\downarrow0$,
\[
\frac{4\alpha}\pi S(h)=\frac{\pi\alpha}h-\frac{2\alpha}\pi\log2+o(1),\qquad \frac{e^{(4\alpha/\pi)S(h)}}{\cosh(\pi\alpha/h)}\longrightarrow2\cdot2^{-2\alpha/\pi}=2^{1-2\alpha/\pi}.
\]
So there are $h>0$, $n\ge1$ and $a>h/2$ for which $\norm{F_h(A)}>2^{1-2\alpha/\pi}-\epsilon$; we fix them.

It remains to replace $F_h$ by a rational function. The map $\chi(z)=(z-1)/(z+1)$ sends $\Sbar\alpha\cup\{\infty\}$ onto the intersection $L$ of two closed disks, and $S_\alpha$ onto its interior, which contains $\chi(1)=0$; by convexity, $L/(1+\delta)\subset\chi(S_\alpha)$ for $\delta>0$. So $G_\delta(w)=F_h\big(\chi^{-1}(w/(1+\delta))\big)$ is holomorphic near $L$, with $|G_\delta|\le1$ on $L$, and, by Runge's theorem, polynomials $p_j$ converge to $G_\delta$ uniformly on $L$, so that $\max_L|p_j|\to\max_L|G_\delta|\le1$. As $f(A)=P^{-1/2}f(D)P^{1/2}$ for $f$ holomorphic near $\sigma(A)$, and $G_\delta\circ\chi\to F_h$ on $\sigma(A)$ as $\delta\downarrow0$, the function $f=p_j\circ\chi/\max_L|p_j|$, for $\delta$ small and $j$ large, is rational with its only pole at $-1$, satisfies $|f|\le1$ on $\Sbar\alpha$, and has $\norm{f(A)}>2^{1-2\alpha/\pi}-\epsilon$. Since $A\in\Kc(\Bf_S)$ and $\epsilon$ is arbitrary, $K(\Kc(\Bf_S))\ge2^{1-2\alpha/\pi}$.
\end{proof}

For $0<\alpha\le\pi/4$, the ratio of $c_\alpha=\sqrt{4-8\alpha/\pi}$ to the lower bound of Proposition~\ref{prop:sectorlower} is $2^x\sqrt{1-x}$, with $x=2\alpha/\pi$, whose maximum, at $x=1-1/(2\log2)$, is $\sqrt{2/(e\log2)}<1.031$.

\begin{remark}\label{rem:sectorcompare}
For $0<\alpha<\pi/4$, the upper bound \eqref{eq:crouzeixbound} of Crouzeix \cite{Cr25} for the constants of $\Kc(\Bf_S)$, $1-\alpha/\pi+\sqrt{2-4\alpha/\pi+\alpha^2/\pi^2}$, is at least $(1+\sqrt2)\sqrt{1-2\alpha/\pi}$, which exceeds $c_\alpha=\sqrt{4-8\alpha/\pi}$. The lower bound $2^{1-2\alpha/\pi}$ of Proposition~\ref{prop:sectorlower} exceeds the lower bound $\pi\sin\alpha/(2\alpha)$ of \cite{CD03}, \cite[Theorem~1]{BC06} and \cite[Lemma~4.1]{Cr16} for $0<\alpha<\pi/4$, and is exceeded by it for $\pi/4<\alpha<\pi/2$, since $\log\big(\pi\sin\alpha/(2\alpha)\big)-(1-2\alpha/\pi)\log2$ is strictly concave on $(0,\pi/2)$ and vanishes at $\alpha=\pi/4$ and $\alpha=\pi/2$.

By \cite[Theorem~4.2]{Cr04}, the constant for $2\times2$ matrices $A$ with $W(A)\subseteq\Sbar\alpha$ and rational $f$ is
\[
\begin{cases}
\displaystyle\sup_{y>0}\,\tanh\frac{\pi y}4\,\frac{\sin\alpha+\sqrt{\sin^2\alpha+\sinh^2(\alpha y)}}{\sinh(\alpha y)}, & 0<\alpha<\pi/4,\\[3ex]
\sqrt2, & \alpha=\pi/4,\\
\displaystyle\sup_{0<t<1}\,\tan\frac{\pi t}4\,\frac{\sin\alpha+\sqrt{\sin^2\alpha-\sin^2(\alpha t)}}{\sin(\alpha t)}, & \pi/4<\alpha<\pi/2.
\end{cases}
\]
Both functions tend to $\pi\sin\alpha/(2\alpha)$ as $y\to0$ and $t\to0$, so, for $\pi/4<\alpha<\pi/2$, this constant is at least $\pi\sin\alpha/(2\alpha)$, and hence exceeds $2^{1-2\alpha/\pi}$. For $0<\alpha<\pi/4$, numerical computations indicate that it is below $2^{1-2\alpha/\pi}$.
\end{remark}

\begin{corollary}\label{cor:sectorclosed}
Let $0<\alpha<\pi/2$ and let $A\in\B(H)$ satisfy $W(A)\subseteq\Sbar\alpha$. For every $d\ge1$ and every $M_d$-valued rational function $F$ with poles off $\Sbar\alpha$,
\[
\norm{F(A)}\le c_\alpha\sup_{\Sbar\alpha}\norm F .
\]
For $\alpha=\pi/4$ the constant $\sqrt2$ is optimal, and for $\alpha<\pi/4$ the optimal constant is at least $2^{1-2\alpha/\pi}$.
\end{corollary}

\begin{proof}
We may assume $\sup_{\Sbar\alpha}\norm F<\infty$. For $\epsilon>0$ put $A_\epsilon=A+\epsilon I$. Adding $\epsilon$ maps $\Sbar\alpha$ into $S_\alpha$, so $\overline{W(A_\epsilon)}\subseteq\Sbar\alpha+\epsilon\subset S_\alpha$; in particular $\sigma(A_\epsilon)\subset S_\alpha$ and $W(A_\epsilon)\subseteq\Sbar\alpha$. Proposition~\ref{prop:sectorbound}, applied to the restriction of $F$ to $S_\alpha$, gives $\norm{F(A_\epsilon)}\le c_\alpha\sup_{\Sbar\alpha}\norm F$. Since $\sigma(A)\subseteq\overline{W(A)}\subseteq\Sbar\alpha$ and $F$ is holomorphic near $\Sbar\alpha$, $F(A_\epsilon)\to F(A)$ in norm as $\epsilon\downarrow0$. Optimality at $\alpha=\pi/4$ follows from \cite[Lemma~4.1]{Cr16}, and the lower bound from Proposition~\ref{prop:sectorlower}.
\end{proof}

\begin{corollary}\label{cor:lens}
Let $0<\alpha<\pi/2$. Let $D_1,D_2$ be open disks of the Riemann sphere whose boundary circles cross, such that $\Om=D_1\cap D_2$ is a nonempty subset of $\C$ with angle $2\alpha$ at the crossing points, and let $b_i$ be a M\"obius transformation mapping $D_i$ onto $\D$. For every $T\in\B(H)$ with $\sigma(T)\subset\Om$ and $\norm{b_1(T)}\le1$, $\norm{b_2(T)}\le1$, every $d\ge1$ and every bounded holomorphic $F:\Om\to M_d$,
\[
\norm{F(T)}\le c_\alpha\sup_\Om\norm F .
\]
The same bound, with $\max_{\overline\Om}\norm F$ in place of $\sup_\Om\norm F$, where $\overline\Om$ is the closure of $\Om$ in the Riemann sphere, holds for every $T\in\B(H)$ with $\sigma(T)\subseteq\overline\Om$ and $\norm{b_1(T)}\le1$, $\norm{b_2(T)}\le1$, and every $M_d$-valued rational function $F$ with poles off $\overline\Om$. For $\alpha=\pi/4$ the constant $\sqrt2$ is optimal in both statements, and for $\alpha<\pi/4$ the optimal constant in both is at least $2^{1-2\alpha/\pi}$.
\end{corollary}

\begin{proof}
By Lemma~\ref{lem:mobius}, $A=m(T)$ satisfies $\sigma(A)\subset S_\alpha$ and $W(A)\subseteq\Sbar\alpha$, and Proposition~\ref{prop:sectorbound}, applied to $F\circ m^{-1}$, gives the first statement. Conversely, for a matrix $A$ with $\sigma(A)\subset S_\alpha$ and $W(A)\subseteq\Sbar\alpha$ and a rational $f$ with poles off $\Sbar\alpha$ and $|f|\le1$ there, $T=m^{-1}(A)$ satisfies the hypotheses of both statements, $f\circ m$ is rational with poles off $\overline\Om$ and $|f\circ m|\le1$ there, and $(f\circ m)(T)=f(A)$. So the lower bounds follow from Proposition~\ref{prop:sectorlower} and, for $\alpha=\pi/4$, from $A=I+\sqrt2N$, with $N=\big(\begin{smallmatrix}0&1\\0&0\end{smallmatrix}\big)$, and $f(z)=(1-z^2)/(1+z^2)$, for which $\norm{f(A)}=\sqrt2$ \cite[Lemma~4.1]{Cr16}.

For the second statement, let $\delta>0$ be small, and let $D_i^\delta\supseteq\overline{D_i}\cap\C$ be obtained by increasing the radius of $D_i$ by $\delta$ if $D_i$ is a disk, decreasing by $\delta$ the radius of the closed disk whose complement is $D_i$, or replacing $D_i=\{\Rp\ell_i>0\}$, with $\ell_i$ affine, by $\{\Rp\ell_i>-\delta\}$. Then $\Om^\delta=D_1^\delta\cap D_2^\delta$ is as in the first statement, with an angle $2\alpha_\delta\to2\alpha$ as $\delta\downarrow0$, and $\sigma(T)\subset\Om^\delta$. A conformal map $b_i^\delta$ of $D_i^\delta$ onto $\D$ satisfies $|b_i^\delta|\le1$ on $\overline{D_i}$, a spectral set for $T$ by von Neumann's inequality, so $\norm{b_i^\delta(T)}\le1$. For $\delta$ so small that $F$ has no poles on $\overline{\Om^\delta}$, the first statement on $\Om^\delta$ gives
\[
\norm{F(T)}\le c_{\alpha_\delta}\sup_{\Om^\delta}\norm F\longrightarrow c_\alpha\max_{\overline\Om}\norm F\qquad(\delta\downarrow0),
\]
since $c_\alpha$ is continuous in $\alpha$ and $F$ is continuous near $\overline\Om$.
\end{proof}

The bound of Corollary~\ref{cor:sectorclosed} extends to unbounded operators. Let $A$ be a closed operator in $H$ with $\ip{x}{Ax}\in\Sbar\alpha$ for all unit vectors $x$ in its domain, that is, $W(A)\subseteq\Sbar\alpha$, and with $\C\setminus\Sbar\alpha$ contained in the resolvent set $\operatorname{res}(A)$; in the terminology of \cite[Section~V.3.10]{KaP}, $A$ is m-sectorial with vertex $0$ and semi-angle $\alpha$. Then $\norm{(\lambda-A)^{-1}}\le\mathrm{dist}(\lambda,\Sbar\alpha)^{-1}$ for $\lambda\notin\Sbar\alpha$ \cite[Theorem~V.3.2]{KaP}. A rational function $F$ bounded on $\Sbar\alpha$ has the form $F=F(\infty)+\sum_{j,m}C_{jm}(\lambda_j-z)^{-m}$ with $\lambda_j\notin\Sbar\alpha$, and we put $F(A)=F(\infty)\otimes I+\sum_{j,m}C_{jm}\otimes(\lambda_j-A)^{-m}$.

\begin{corollary}\label{cor:unbounded}
Let $0<\alpha<\pi/2$, and let $A$ be a closed operator with $W(A)\subseteq\Sbar\alpha$ and $\C\setminus\Sbar\alpha\subseteq\operatorname{res}(A)$. For every $d\ge1$ and every $M_d$-valued rational function $F$ that is bounded on $\Sbar\alpha$,
\[
\norm{F(A)}\le c_\alpha\sup_{\Sbar\alpha}\norm F .
\]
\end{corollary}

\begin{proof}
For $n\ge1$ we have $\mathrm{dist}(-n,\Sbar\alpha)=n$, so $\norm{(n+A)^{-1}}\le1/n$, and $A_n=A(I+A/n)^{-1}=n-n^2(n+A)^{-1}$ is bounded; this Yosida approximation is used in \cite[Section~3]{Cr25}, and a related one in the proof of \cite[Theorem~3.13]{MR25}. For $x\in H$ put $y=(I+A/n)^{-1}x$, which lies in the domain of $A$. Then $x=y+Ay/n$ and $A_nx=Ay$, so
\[
\ip{x}{A_nx}=\ip{y}{Ay}+\tfrac1n\norm{Ay}^2\in\Sbar\alpha .
\]
Hence $W(A_n)\subseteq\Sbar\alpha$, and Corollary~\ref{cor:sectorclosed} gives $\norm{F(A_n)}\le c_\alpha\sup_{\Sbar\alpha}\norm F$. For $\lambda\notin\Sbar\alpha$ put $\lambda_n=n\lambda/(n-\lambda)$, so that $\lambda_n\to\lambda$ and $\lambda_n\in\operatorname{res}(A)$ for large $n$. For such $n$, since $\lambda-A_n=(n-\lambda)(\lambda_n-A)(n+A)^{-1}$ and $n+A=n^2/(n-\lambda)-(\lambda_n-A)$,
\[
(\lambda-A_n)^{-1}=\frac{n^2}{(n-\lambda)^2}\,(\lambda_n-A)^{-1}-\frac1{n-\lambda}.
\]
Since the resolvent is norm-continuous on $\operatorname{res}(A)$, $(\lambda-A_n)^{-1}\to(\lambda-A)^{-1}$ in norm. Hence $F(A_n)\to F(A)$ in norm.
\end{proof}

\begin{remark}\label{rem:hinfty}
If $A$ in Corollary~\ref{cor:unbounded} is injective, the bound extends to the $H^\infty$-calculus of $A$: for every $\theta\in(\alpha,\pi)$ and every bounded holomorphic $F:S_\theta\to M_d$, where the entries of $F(A)$ are defined by the functional calculus of \cite{Ha06}, with the matrix convention of Section~\ref{ssec:calc},
\[
\norm{F(A)}\le c_\alpha\sup_{S_\alpha}\norm F .
\]
For $d=1$, this follows from Corollary~\ref{cor:sectorclosed} and the proof of \cite[Theorem~3.13]{MR25}, which transfers a scalar bound for bounded operators with numerical range in a sector to the $H^\infty$-calculus with the same constant; we give the argument, which follows theirs, for general $d$. Since $\norm{(\lambda-A)^{-1}}\le\mathrm{dist}(\lambda,\Sbar\alpha)^{-1}$ for $\lambda\notin\Sbar\alpha$, $A$ is sectorial of angle $\alpha$ in the sense of \cite[Section~2.1]{Ha06}. Suppose first that the entries of $F$ lie in the class $H^\infty_0(S_\theta)$ of bounded holomorphic $f$ on $S_\theta$ with $|f(z)|\le C\min\{|z|^\kappa,|z|^{-\kappa}\}$, $C,\kappa>0$, and let $\alpha<\delta<\theta$. Then $F(A)$ and $F(z)$, $z\in\Sbar\alpha\setminus\{0\}$, are the Cauchy integrals of $F(w)\otimes(w-A)^{-1}$ and $F(w)(w-z)^{-1}$ over $\partial S_\delta$. Since $|w-z|$ and $\norm{(w-A)^{-1}}^{-1}$ are at least $|w|\sin(\delta-\alpha)$ there, fine Riemann sums over $\{w\in\partial S_\delta:1/n\le|w|\le n\}$ give rational $F_n$ with poles on $\partial S_\delta$, $F_n\to F$ uniformly on $\Sbar\alpha\setminus\{0\}$ and $F_n(A)\to F(A)$ in norm; so Corollary~\ref{cor:unbounded} gives the bound.

For general $F$ we may assume $\theta<\pi/2$, since $F(A)$ depends only on the restriction of $F$ to $S_{\theta'}$ for any $\theta'\in(\alpha,\theta)$. Let $q_\epsilon(z)=z/((1+\epsilon z)(z+\epsilon))$ be as in the proof of Proposition~\ref{prop:sectorbound}; by \eqref{eq:qeps}, $|q_\epsilon|\le1$ on $S_\theta$ and the entries of $q_\epsilon F$ lie in $H^\infty_0(S_\theta)$, so $\norm{(q_\epsilon F)(A)}\le c_\alpha\sup_{S_\alpha}\norm F$. As $\epsilon\downarrow0$, $q_\epsilon\to1$ pointwise; since $H$ is reflexive and $A$ is injective, $A$ has dense domain and dense range \cite[Proposition~2.1.1(h)]{Ha06}, so the convergence lemma \cite[Proposition~5.1.4(b)]{Ha06}, applied to the entries, gives $(q_\epsilon F)(A)\to F(A)$ strongly, hence the bound for $F$.

For the $H^\infty$-calculus of injective operators $A$ as in Corollary~\ref{cor:unbounded}, the scalar bound $2+2/\sqrt3$ \cite{CD03} was known, already with the supremum over $S_\alpha$ \cite[Corollary~7.1.17]{Ha06}, and the scalar bound $2$ follows from \cite[Theorem~3.13]{MR25} combined with \cite{Jin,LS}; the complete bound $2+2/\sqrt3$, also with the supremum over $S_\alpha$, follows from Le Merdy's proof of \cite[Theorem~1.1]{LeM03} (see \cite[(3.8)]{LeM03} and the paragraph after \cite[(3.4)]{LeM03}). For $\alpha>\alpha_*$, the bound \eqref{eq:crouzeixbound} of \cite{Cr25}, which is then smaller than $c_\alpha=\sqrt2$, holds for the operators of Corollary~\ref{cor:unbounded} \cite[Section~3]{Cr25}, also for the complete constant, and the argument of this remark, with \eqref{eq:crouzeixbound} in place of $c_\alpha$, extends it to their $H^\infty$-calculus.
\end{remark}

\FloatBarrier
\section{Domains with real-analytic boundary}\label{sec:convex}

This section depends only on Sections~\ref{sec:setting}--\ref{sec:positivity}. Let $\Om$ be a bounded Jordan domain whose boundary is a regular real-analytic Jordan curve, that is, a Jordan curve with a real-analytic parametrization of nonvanishing derivative. Let $\psi:\{|t|>1\}\to\C\setminus\overline\Om$ be the exterior conformal map with $\psi(t)=ct+O(1)$, $c>0$; it extends univalently to $\{|t|>\tau_0\}$ for some $\tau_0<1$ \cite[Proposition~3.1]{Pom92}; for convex $\Om$ see also \cite[Lemma~A.1]{OAI}. For $|t|>\tau_0$ and $z\in\C$ with $z\ne\psi(t)$, in particular for $t\in\T$ and all $z\in\Om$, put
\[
r_t(z)=\frac{t\psi'(t)}{\psi(t)-z}.
\]
For $t\in\T$, with $\zeta=\psi(t)$, the outward unit normal $n(\zeta)$ and arclength $ds$ on $\partial\Om$, one has $ds=2\pi|\psi'(t)|\,d\omega(t)$ and $t\psi'(t)=|\psi'(t)|\,n(\zeta)$, since $it\psi'(t)$ is a tangent vector to $\partial\Om$ at $\zeta$ and $\partial\Om$ is traversed counterclockwise. So, for $\sigma(T)\subset\Om$,
\begin{equation}\label{eq:dlp}
\Rp r_t(T)=\pi|\psi'(t)|\,\mu(\zeta,T),\qquad \mu(\zeta,T)=\tfrac1\pi\Rp\big[n(\zeta)(\zeta-T)^{-1}\big],
\end{equation}
where $\mu$ is the operator-valued double-layer potential of \cite{DD,CP17}.

\begin{lemma}\label{lem:convexclass}
For every compact $C\subset\Om$ there is $\tau_C\in(\tau_0,1)$ such that $(t,z)\mapsto r_t(z)$ is bounded and continuous on $\{|t|>\tau_C\}\times C$, holomorphic in $t$, and tends to $1$ as $t\to\infty$, uniformly in $z\in C$. If $\Om$ is convex, $r_t$ is a Herglotz function on $\Om$ for every $t\in\T$, and for $T\in\B(H)$ with $\sigma(T)\subset\Om$ one has $\Rp r_t(T)\ge0$ for all $t\in\T$ if and only if $W(T)\subseteq\overline\Om$.
\end{lemma}

\begin{proof}
Let $\delta=\mathrm{dist}(C,\partial\Om)>0$, and fix $\tau_1\in(\tau_0,1)$. Since $\psi$ is uniformly continuous on $\{\tau_1\le|t|\le1\}$ and maps $\T$ onto $\partial\Om$, there is $\tau_C\in[\tau_1,1)$ with $|\psi(t)-\psi(t/|t|)|<\delta/2$ for $\tau_C\le|t|\le1$. As $\psi(t)\notin\Om$ for $|t|\ge1$, it follows that $|\psi(t)-z|\ge\delta/2$ for $|t|\ge\tau_C$ and $z\in C$; with $\psi(t)=ct+O(1)$ and $t\psi'(t)=ct+O(t^{-1})$, this gives the first statement. If $\Om$ is convex, the last statement is \cite[Proposition~3.4]{SdVIEOT}, by \eqref{eq:dlp}; the version with $W(T)\subset\Om$ and strict positivity is in \cite{BCD06}. Applied to $T=z$ on $H=\C$, with $z\in\Om$, it gives $\Rp r_t\ge0$ on $\Om$, and $\Rp r_t>0$ by the minimum principle, since $\Rp r_t$ is harmonic and nonconstant.
\end{proof}

By Lemma~\ref{lem:convexclass}, for $z$ in a compact $C\subset\Om$ the function $t\mapsto r_t(z)$ has a Laurent expansion
\[
r_t(z)=1+\sum_{n\ge1}\Phi_n(z)\,t^{-n},
\]
which converges absolutely and uniformly on $C\times\{|t|\ge\tau\}$ for every $\tau>\tau_C$. The coefficients $\Phi_n$ are the Faber polynomials of $\Om$ (the polynomials $b_n$ of \cite[(2.3)]{OAI}); we put $\Phi_0=1$. For $\rho>0$ let
\[
\Kc_\rho(\Om)=\big\{T\in\B(H):\ \sigma(T)\subset\Om,\ \Rp r_t(T)\ge\big(1-\tfrac\rho2\big)I\text{ for all }t\in\T\big\}.
\]
The functions $r_{t,\rho}=(1-\frac2\rho)+\frac2\rho r_t$ satisfy $\Rp r_{t,\rho}(T)\ge0$ exactly when $\Rp r_t(T)\ge(1-\frac\rho2)I$, and the classes $\Kc_\rho(\Om)$ increase with $\rho$. If $\Om$ is convex and $\rho\ge2$, the $r_{t,\rho}$ $(t\in\T)$ are Herglotz functions, $\Kc_\rho(\Om)$ is the class of the family $(r_{t,\rho})_{t\in\T}$, and $\Kc_2(\Om)=\{T:\sigma(T)\subset\Om,\ W(T)\subseteq\overline\Om\}$ by Lemma~\ref{lem:convexclass}. For general $\Om$, the least $\rho\ge2$ for which every $r_{t,\rho}$ has nonnegative real part on $\Om$ is
\begin{equation}\label{eq:rhoOm}
\rho_\Om=2+2\sup_{t\in\T,\,z\in\Om}\big(\Rp r_t(z)\big)_-,
\end{equation}
where $x_-=\max\{-x,0\}$; by Lemma~\ref{lem:rhoOm}(a) below, $\rho_\Om<\infty$, and $\rho_\Om=2$ exactly when $\Om$ is convex. Since, for normal $T$ with $\sigma(T)\subset\Om$, $\Rp r_t(T)\ge(1-\frac\rho2)I$ exactly when $\Rp r_t\ge1-\frac\rho2$ on $\sigma(T)$, for $\rho\ge2$ the class $\Kc_\rho(\Om)$ contains every normal operator with spectrum in $\Om$ exactly when $\rho\ge\rho_\Om$.

The boundary values of $r_t$ are described by the function
\[
S(t,u)=\frac{t\psi'(t)}{\psi(t)-\psi(u)}-\frac t{t-u}=t\,\partial_t\log\frac{\psi(t)-\psi(u)}{t-u}.
\]
Since $\psi$ is univalent, the difference quotient $(\psi(t)-\psi(u))/(t-u)$ is holomorphic and zero-free in $(1/t,1/u)$ on $(\{|t|>\tau_0\}\cup\{\infty\})^2$, with the values $\psi'(t)$ at $u=t$ and $c$ at $t=\infty$ or $u=\infty$. So $S$ is holomorphic there and vanishes when $t=\infty$ or $u=\infty$; hence $S(t,u)=\sum_{n,j\ge1}g_{nj}t^{-n}u^{-j}$, absolutely and uniformly on $\T\times\T$, where $-g_{nj}/n$ is the coefficient of $t^{-n}u^{-j}$ in $\log\frac{\psi(t)-\psi(u)}{c(t-u)}$, so that the numbers $b_{nj}=g_{nj}/n$ are the Grunsky coefficients of $\psi$ in the convention of \cite[Chapter~3]{Pom75}. Comparing the coefficients of $t^{-n}$ in $r_t(\psi(u))=t/(t-u)+S(t,u)$ for $\tau_0<|u|<1$, where $\psi(u)\in\Om$, and continuing analytically in $u$, gives
\begin{equation}\label{eq:faberboundary}
\Phi_n(\psi(u))=u^n+\sum_{j\ge1}g_{nj}u^{-j}\qquad(n\ge1,\ |u|>\tau_0),
\end{equation}
the boundary formula for Faber polynomials used by K\"ovari and Pommerenke \cite{KP67,Pom65}, see \cite{BC14}; this is the argument of \cite[proof of Theorem~2.1, (2.4)--(2.6)]{OAI}. On $\T\times\T$, $2\Rp[t/(t-u)]=1$ for $t\ne u$, so the function
\[
k_\partial(t,u)=2\Rp r_t(\psi(u))=1+2\Rp S(t,u)\qquad(t\ne u)
\]
extends continuously to $\T\times\T$. On the diagonal, $k_\partial(t,t)=1+\Rp\big(t\psi''(t)/\psi'(t)\big)$ is $|\psi'(t)|$ times the curvature of $\partial\Om$ at $\psi(t)$; and with $\zeta=\psi(t)$, $k_\partial(t,u)\,d\omega(t)=\pi^{-1}d_\zeta\arg(\zeta-\psi(u))$ is the double-layer potential kernel of $\partial\Om$, normalized as in \eqref{eq:dlp}. Let $\Theta$ be the integral operator $(\Theta f)(u)=\int_\T k_\partial(t,u)f(t)\,d\omega(t)$ on $L^2(\T)$, and let $V(\Om)$ be the norm of its restriction to the closed span of the functions $t^n$, $n\ge1$. By Lemma~\ref{lem:rhoOm}(c) below, $V(\Om)$ is the norm on $\ell^2$ of the matrix $(nb_{nj})_{n,j\ge1}$. Grunsky's inequality says that the symmetric matrix $(\sqrt{nj}\,b_{nj})_{n,j\ge1}$ has norm at most $1$ for every $\Om$ \cite[Chapter~3]{Pom75}; since $nb_{nj}=\sqrt{n/j}\,\sqrt{nj}\,b_{nj}$, it does not bound $V(\Om)$.

\begin{lemma}\label{lem:rhoOm}
\begin{enumerate}[(a)]
\item The constant $\rho_\Om$ of \eqref{eq:rhoOm} is finite and
\[
\rho_\Om=2+2\sup_{t,u\in\T,\ t\ne u}\big(\Rp r_t(\psi(u))\big)_-=2+\max\Big\{0,-\min_{\T\times\T}k_\partial\Big\};
\]
moreover, $\rho_\Om=2$ if and only if $\Om$ is convex.
\item Let $\rho\ge\rho_\Om$. Then $\Rp r_{t,\rho}>0$ on $\Om$ for every $t\in\T$, so that $(r_{t,\rho})_{t\in\T}$ is a family of Herglotz functions, and $\Kc_\rho(\Om)$ is its class. If $\rho>2$, put
\[
a_t=\psi(t)+\frac{t\psi'(t)}{\rho-2},\qquad R_t=\frac{|\psi'(t)|}{\rho-2}\qquad(t\in\T);
\]
then the disk $\{|z-a_t|<R_t\}$ does not meet $\overline\Om$, and $\mathrm{dist}(a_t,\overline\Om)=R_t$.
\item For $n\ge1$,
\[
\Theta t^n=\Phi_n\circ\psi-u^n=\sum_{j\ge1}g_{nj}u^{-j},
\]
so $V(\Om)$ is the norm of the matrix $(g_{nj})_{n,j\ge1}$ on $\ell^2$. If $\Om$ is convex, then $V(\Om)\le1$.
\item Write $\psi(t)=c\big(t+\sum_{k\ge0}\psi_kt^{-k}\big)$ for $|t|>\tau_0$. If $\sum_{k\ge1}k|\psi_k|\le\frac12$, then $V(\Om)\le1$.
\end{enumerate}
\end{lemma}

\begin{proof}
(a) Fix $t\in\T$, and put $p=\psi(t)$, $v=t\psi'(t)$ and $m_t=\inf_{u\ne t}\Rp r_t(\psi(u))=\frac12\min_{u\in\T}k_\partial(t,u)$. For $\kappa>0$ let $D$ be the disk $\{|z-p-v/(2\kappa)|<|v|/(2\kappa)\}$. For $z\ne p$ write $w=(p-z)/v=1/r_t(z)$; then $\Rp r_t(z)\ge-\kappa$ if and only if $\Rp w\ge-\kappa|w|^2$, that is, $|w+1/(2\kappa)|\ge1/(2\kappa)$, that is, $z\notin D$, since $z-p-v/(2\kappa)=-v\big(w+1/(2\kappa)\big)$; likewise, $\Rp r_t(z)>-\kappa$ if and only if $z\notin\overline D$. Now let $\kappa\ge-m_t$. Then $D$ does not meet $\partial\Om$; it is connected and contains the points $p+\epsilon v$, $0<\epsilon<1/\kappa$, which lie outside $\overline\Om$ for small $\epsilon$, since $v=|\psi'(t)|\,n(p)$. Hence $D\cap\overline\Om=\emptyset$; and since $\Om$ is open, it does not meet the boundary circle of $D$ either, so $\Rp r_t>-\kappa$ on $\Om$. Thus $\inf_\Om\Rp r_t\ge-(m_t)_-$; and $\inf_\Om\Rp r_t\le m_t$, since $\Rp r_t$ is continuous on $\overline\Om\setminus\{p\}$. So $\sup_{z\in\Om}(\Rp r_t(z))_-=(m_t)_-$, which gives the formula in (a); it is finite because $k_\partial$ is continuous on $\T\times\T$.

By the formula in (a), $\rho_\Om=2$ if and only if $k_\partial\ge0$ on $\T\times\T$, that is, if and only if the double-layer potential kernel of $\partial\Om$ is nonnegative; by \cite[Proposition~2.4]{SdVIEOT}, this holds if and only if $\Om$ is convex. This proves (a).

(b) If $\rho=2$, then $\Om$ is convex by (a), and Lemma~\ref{lem:convexclass} applies. If $\rho>2$, then $\kappa=\frac\rho2-1>0$ satisfies $\kappa\ge(m_t)_-\ge-m_t$ for every $t$, by (a), and $D$ is the disk $\{|z-a_t|<R_t\}$. By the first paragraph, $\Rp r_t>-\kappa$ on $\Om$, so that $\Rp r_{t,\rho}>1-\frac2\rho-\frac2\rho\kappa=0$; since $\Rp r_{t,\rho}(T)\ge0$ exactly when $\Rp r_t(T)\ge(1-\frac\rho2)I$, $\Kc_\rho(\Om)$ is the class of $(r_{t,\rho})$. Also by the first paragraph, $D$ does not meet $\overline\Om$, so $\mathrm{dist}(a_t,\overline\Om)=R_t$, since $p\in\overline\Om$ lies on its boundary circle.

(c) On $\T\times\T$, $k_\partial=1+S+\overline S$, where $S$ contains only the powers $t^{-n}u^{-j}$ and $\overline S$ only the powers $t^nu^j$, $n,j\ge1$. Integrating against $1$ and $t^n$ in $t$ gives $\Theta1=1$ and $\Theta t^n=\sum_jg_{nj}u^{-j}$, which is $\Phi_n\circ\psi-u^n$ by \eqref{eq:faberboundary}; since the $t^n$ and the $u^{-j}$ are orthonormal, $V(\Om)$ is the norm of $(g_{nj})$ on $\ell^2$. If $\Om$ is convex, then $k_\partial\ge0$, by Lemma~\ref{lem:convexclass} and continuity, and $\int_\T k_\partial\,d\omega(t)=\int_\T k_\partial\,d\omega(u)=1$; so $V(\Om)\le\norm\Theta\le1$ by the Schur test, as in the proof of \cite[Proposition~2.2]{OAI}.

(d) Put $x=\sum_{k\ge1}k|\psi_k|\le\frac12$. Since $(t^{-k}-u^{-k})/(t-u)=-\sum_{i=1}^kt^{-i}u^{i-k-1}$, we have, for $|t|,|u|\ge1$, $(\psi(t)-\psi(u))/(c(t-u))=1-X$ with $|X|\le x$, and $S=t\,\partial_t\log(1-X)=Y/(1-X)$, where
\[
X=\sum_{k\ge1}\psi_k\sum_{i=1}^kt^{-i}u^{i-k-1},\qquad Y=-t\,\partial_tX=\sum_{k\ge1}\psi_k\sum_{i=1}^ki\,t^{-i}u^{i-k-1}.
\]
Expanding $S=\sum_{l\ge0}YX^l$ in powers of $t^{-1}$ and $u^{-1}$, which is legitimate for large $|t|$ and $|u|$, shows that $|g_{nj}|$ is at most the coefficient of $t^{-n}u^{-j}$ in $\sum_{l\ge0}\hat Y\hat X^l$, where $\hat X$ and $\hat Y$ are $X$ and $Y$ with $|\psi_k|$ in place of $\psi_k$. Summing over $j$ amounts to setting $u=1$. Then the coefficients of $\hat X$, as a series in $t^{-1}$, sum to $x$, and the coefficient of $t^{-i}$ in $\hat Y$ is $i\sum_{k\ge i}|\psi_k|\le x$; so the coefficients of $\hat Y\hat X^l$ are at most $x^{l+1}$, and $\sum_j|g_{nj}|\le\sum_{l\ge0}x^{l+1}=x/(1-x)\le1$ for every $n$. Likewise, at $t=1$ the coefficients of $\hat X$, as a series in $u^{-1}$, sum to $x$, and the coefficient of $u^{-p}$ in $\hat Y$ is $\sum_{k\ge p}(k+1-p)|\psi_k|\le x$; so $\sum_n|g_{nj}|\le1$ for every $j$. By (c) and the Schur test, $V(\Om)\le1$.
\end{proof}

The convex case of (c) has an analogue in \cite{MMOR25}: for convex $\Om$, the norm $a(\Om)$ of the double-layer potential operator on the functions continuous on $\overline\Om$ and holomorphic in $\Om$, modulo constants, is less than $1$ \cite[Theorem~2]{MMOR25}.

Part (d) applies to nonconvex domains. For $m\ge3$ let
\[
\psi(t)=t+\frac{t^{-m}}{2m}.
\]
For $t\ne u$ with $|t|,|u|\ge r>2^{-1/(m+1)}$,
\[
\Big|\frac{\psi(t)-\psi(u)}{t-u}\Big|=\Big|1-\frac1{2m}\sum_{i=1}^mt^{-i}u^{i-m-1}\Big|\ge1-\frac{r^{-m-1}}2>0,
\]
so $\psi$ is univalent on $\{|t|>2^{-1/(m+1)}\}$, a neighbourhood of $\{|t|\ge1\}$. Hence $\psi(\T)$ bounds a Jordan domain $\Om_m$ with regular real-analytic boundary, whose exterior conformal map is $\psi$, and $V(\Om_m)\le1$ by (d), since $\sum_kk|\psi_k|=1/2$. As $|\psi|\le1+1/(2m)\le7/6$ on $\T$, $\Om_m$ lies in the disk $\{|z|\le7/6\}$. But $\Om_m$ is not convex: where $t^{m+1}=-1$,
\[
k_\partial(t,t)=1+\Rp\frac{t\psi''(t)}{\psi'(t)}=-\frac{m-2}3<0,
\]
so the curvature of $\partial\Om_m$ at $\psi(t)$ is negative. So the domains $\Om_m$, $m\ge3$, all lie in one disk and have $V(\Om_m)\le1$, while the minimum of their double-layer potential kernels tends to $-\infty$ as $m\to\infty$.

An operator $T\in\B(H)$ has a \emph{Faber $\rho$-dilation} if there is a unitary $U$ on a Hilbert space containing $H$ with
\[
\Phi_n(T)=\rho P_HU^n|_H\qquad(n\ge1),
\]
where $P_H$ is the orthogonal projection onto $H$; for $\Om=\D$ these are the $\rho$-contractions \cite[Section~I.11]{SzNF}.

\begin{lemma}\label{lem:rhoconvex}
Let $\rho>0$, and let $T\in\B(H)$ with $\sigma(T)\subset\Om$.
\begin{enumerate}[(a)]
\item $T\in\Kc_\rho(\Om)$ if and only if $T$ has a Faber $\rho$-dilation.
\item Let $\rho>2$ and $\rho\ge\rho_\Om$, and let $a_t$ be as in Lemma~\ref{lem:rhoOm}(b). Then $T\in\Kc_\rho(\Om)$ if and only if $\norm{(a_t-T)^{-1}}\le\mathrm{dist}(a_t,\overline\Om)^{-1}$ for all $t\in\T$.
\end{enumerate}
\end{lemma}

\begin{proof}
(a) Choose a compact $C\subset\Om$ that contains a cycle $\Gamma$ as in Section~\ref{ssec:calc}, and $\tau\in(\tau_C,1)$. Inserting the expansion of $r_t$ into the Riesz--Dunford integral over $\Gamma$ gives $r_t(T)=I+\sum_{n\ge1}\Phi_n(T)t^{-n}$ for $|t|\ge1$, with $\sum_n\norm{\Phi_n(T)}\tau^{-n}<\infty$. Hence, for $t\in\T$,
\[
B(t):=\Rp r_{t,\rho}(T)=I+\frac1\rho\sum_{n\ge1}\big(\Phi_n(T)t^{-n}+\Phi_n(T)^*t^n\big),
\]
an absolutely convergent series, and $T\in\Kc_\rho(\Om)$ if and only if $B(t)\ge0$ for all $t\in\T$. Put $\gamma_0=I$, $\gamma_n=\Phi_n(T)/\rho$ and $\gamma_{-n}=\gamma_n^*$ for $n\ge1$. Then $\gamma_n=\int_\T t^nB(t)\,d\omega(t)$ for all $n\in\Z$, so $\sum_{m,n}\ip{h_n}{\gamma_{m-n}h_m}=\int_\T\ip{h(t)}{B(t)h(t)}\,d\omega(t)$ for every finitely supported sequence $(h_m)$ in $H$ and $h(t)=\sum_mt^mh_m$. If $B(t)\ge0$ for all $t$, the sequence $(\gamma_n)$ is therefore positive definite on $\Z$, and Naimark's dilation theorem \cite[Chapter~I]{SzNF} gives a unitary $U$ on a Hilbert space containing $H$ with $\gamma_n=P_HU^n|_H$ for $n\ge1$. Conversely, if $\Phi_n(T)=\rho P_HU^n|_H$, then for $s>1$ and $t\in\T$,
\[
r_{st,\rho}(T)=I+2\sum_{n\ge1}(st)^{-n}P_HU^n|_H=P_H\big[(I+Y)(I-Y)^{-1}\big]\big|_H,\qquad Y=(st)^{-1}U,
\]
which has positive real part by Lemma~\ref{lem:cayley}. As $s\downarrow1$, $r_{st,\rho}(T)\to r_{t,\rho}(T)$ in norm, so $B(t)\ge0$.

(b) Put $\kappa=\frac\rho2-1>0$, $p=\psi(t)$, $v=t\psi'(t)$ and $Z=p-T$, so that $a_t=p+v/(2\kappa)$ and, by Lemma~\ref{lem:rhoOm}(b), $\mathrm{dist}(a_t,\overline\Om)=R_t=|v|/(2\kappa)$. Since $p$ and $a_t$ lie outside $\Om\supseteq\sigma(T)$, $Z$ and $a_t-T=Z+v/(2\kappa)$ are invertible. From $r_t(T)=vZ^{-1}$ we get $Z^*(\Rp r_t(T))Z=\Rp(\bar vZ)$, and expanding $\kappa(a_t-T)^*(a_t-T)$ gives
\[
Z^*\big(\Rp r_t(T)+\kappa I\big)Z=\kappa\big((a_t-T)^*(a_t-T)-R_t^2I\big).
\]
So $\Rp r_t(T)\ge-\kappa I=(1-\frac\rho2)I$ if and only if $\norm{(a_t-T)x}\ge R_t\norm x$ for all $x\in H$, that is, $\norm{(a_t-T)^{-1}}\le1/R_t$.
\end{proof}

\begin{remark}\label{rem:special}
Special cases of both parts of Lemma~\ref{lem:rhoconvex} are known. For $\Om=\D$, where $\psi(t)=t$ and $\Phi_n=z^n$, (a) is the characterization of $\rho$-contractions by $\Rp(I-\bar tT)^{-1}\ge(1-\frac\rho2)I$ \cite[Section~I.11]{SzNF}, \cite[Proposition~3.1]{BCK18}, and (b) says that $\norm{(\lambda-T)^{-1}}\le\rho-2$ on the circle $|\lambda|=(\rho-1)/(\rho-2)$; compare \cite[Section~3.4]{DEY17}. For $\rho=2$, (a) is the Naimark dilation of the positive operator-valued measure $\frac12\mu(\zeta,T)\,ds$, of total mass $I$ \cite{DD,PS05,BCD06}, combined with Beckermann's representation of $\Phi_n(T)/2$ as a moment of this measure \cite[(10)]{Be05}, \cite{BC14}; it gives Beckermann's bound $\norm{\Phi_n(T)}\le2$ \cite[Theorem~1.1]{Be05}. One direction of (b) is \cite[Lemma~2(3)]{CG26}, which goes back to \cite[Lemma~7]{CG19}, where it is proved by an identity of the same kind. For convex $\Om$, the resolvent bound $\norm{(\lambda-T)^{-1}}\le\mathrm{dist}(\lambda,\overline\Om)^{-1}$ for all $\lambda\notin\overline\Om$ is equivalent to $W(T)\subseteq\overline\Om$ (to deduce $W(T)\subseteq\overline\Om$, let $\lambda\to\infty$ along the outward normals); $\Kc_\rho(\Om)$ requires it only on the curve $\{a_t:t\in\T\}$. For $\rho>2$ with $\rho\ge\rho_\Om$, the disks $\{|z-a_t|<\mathrm{dist}(a_t,\overline\Om)\}$ lie outside $\overline\Om$ and are tangent to $\partial\Om$ at $\psi(t)$; shrinking them, tangent at the same points, to a common radius preserves the resolvent bounds \cite[Lemma~7(ii)]{DEY17}, so every $T\in\Kc_\rho(\Om)$ satisfies, on sufficiently short boundary arcs, the hypotheses of \cite[Theorem~2]{DEY17}, which make $\overline\Om$ a complete $K$-spectral set for $T$ with $K$ not explicit; since the classes $\Kc_\rho(\Om)$ increase with $\rho$, the same holds for every $\rho>0$.
\end{remark}

The kernel of Proposition~\ref{prop:pencil} for the family $(r_{t,\rho})$ is, by Parseval's identity,
\[
k_{\Om,\rho}(z,w)=1+\int_\T r_{t,\rho}(z)\overline{r_{t,\rho}(w)}\,d\omega(t)=2+\frac4{\rho^2}\sum_{n\ge1}\Phi_n(z)\overline{\Phi_n(w)} .
\]
Thus $\E_{\Om,\rho}=\E_{k_{\Om,\rho}}$ is the space of Faber series with norm
\[
\Big\|c_0+\sum_{n\ge1}c_n\Phi_n\Big\|^2=\tfrac12|c_0|^2+\tfrac{\rho^2}4\sum_{n\ge1}|c_n|^2,
\]
and $1\in\E_{\Om,\rho}$. This is \eqref{eq:intspace} with $\Lambda=\{0,1,2,\dots\}$, $\phi_n=\Phi_n$, $\beta(0)^2=\frac12$ and $\beta(n)^2=\rho^2/4$ for $n\ge1$. Its injectivity hypothesis holds because the coefficients are recovered by the functionals
\[
L_m(f)=\frac1{2\pi i}\int_{|t|=\tau}f(\psi(t))\,t^{-m-1}\,dt\qquad(m\ge0),
\]
with $\tau_0<\tau<1$ fixed. Indeed, $L_m(\Phi_n)=\delta_{mn}$, since $\Phi_n(\psi(t))-t^n$ is holomorphic on $\{|t|>\tau_0\}$ and vanishes at $\infty$, by \eqref{eq:faberboundary}, and $L_m$ is continuous for locally uniform convergence, under which $\sum_nc_n\Phi_n$ converges for square-summable $(c_n)$, since $\sum_{n\ge1}|\Phi_n(z)|^2=\int_\T|r_t(z)-1|^2\,d\omega(t)$ is locally bounded. Hence $\sqrt2,\frac2\rho\Phi_1,\frac2\rho\Phi_2,\dots$ is an orthonormal basis of $\E_{\Om,\rho}$, and by Lemma~\ref{lem:conv}(c) the moment map is the absolutely convergent series
\[
\Psi_T(X)=2X+\frac4{\rho^2}\sum_{n\ge1}\Phi_n(T)\,X\,\Phi_n(T)^* .
\]

The multiplier constant of $\E_{\Om,\rho}$ is bounded by comparing its norm with the boundary norm $\norm Q_\partial^2=\int_\T\hs{Q(\psi(t))}^2\,d\omega(t)$ of matrix polynomials $Q$. For convex $\Om$, where $V(\Om)\le1$, the comparison \eqref{eq:faber} below is \cite[Theorem~2.1, Proposition~2.2]{OAI} with $1+V(\Om)^2$ in place of $2$, and has the same proof: $\Theta$ is the integral operator of \cite[(2.8)]{OAI}, whose kernel $k_\partial$ is that of \cite[Theorem~2.1]{OAI}, and convexity enters there only through the sign of $k_\partial$.

\begin{lemma}\label{lem:faber}
Let $\rho\ge\sqrt2$. For $d\ge1$ and every matrix polynomial $Q=\sum_{n\ge0}Q_n\Phi_n$ with $Q_n\in M_d$,
\begin{equation}\label{eq:faber}
\sum_{n\ge0}\hs{Q_n}^2\le\norm{Q}_\partial^2\le\hs{Q_0}^2+\big(1+V(\Om)^2\big)\sum_{n\ge1}\hs{Q_n}^2 .
\end{equation}
Moreover,
\begin{equation}\label{eq:convexmult}
M(\E_{\Om,\rho})\le\max\big\{\tfrac\rho{\sqrt2},\sqrt{1+V(\Om)^2}\big\}.
\end{equation}
\end{lemma}

\begin{proof}
By \eqref{eq:faberboundary} and Lemma~\ref{lem:rhoOm}(c),
\[
Q\circ\psi=Q_0+Q_++\Theta Q_+\ \text{ on }\T,\qquad Q_+(u)=\sum_{n\ge1}Q_nu^n,
\]
where $\Theta$ acts entrywise on $L^2(\T;M_d)$, with the Hilbert--Schmidt norm on $M_d$. Since the entries of $Q_+$ lie in the closed span of the functions $t^n$, $n\ge1$, we have $\norm{\Theta Q_+}\le V(\Om)\norm{Q_+}$, and since the three terms have disjoint Fourier supports, Parseval's identity gives
\[
\norm Q_\partial^2=\hs{Q_0}^2+\norm{Q_+}^2+\norm{\Theta Q_+}^2,
\]
with $L^2(\T;M_d)$ norms; this proves \eqref{eq:faber}. Since the norm of $\E_d=\E_{\Om,\rho}\otimes M_d$ is $\norm Q_{\E_d}^2=\frac12\hs{Q_0}^2+\frac{\rho^2}4\norm{Q_+}^2$ and $\rho\ge\sqrt2$, it follows that
\[
\tfrac4{\rho^2}\norm Q_{\E_d}^2\le\norm Q_\partial^2\le\max\big\{2,\tfrac4{\rho^2}\big(1+V(\Om)^2\big)\big\}\norm Q_{\E_d}^2 .
\]
So $c_-=4/\rho^2$ and $c_+=\max\{2,4(1+V(\Om)^2)/\rho^2\}$ in the notation of Lemma~\ref{lem:compare}, whose argument we now repeat for matrix polynomials.

If $F$ is a matrix polynomial, then $FQ$ is again a matrix polynomial and $\hs{F(\psi(t))Q(\psi(t))}\le\max_{\overline\Om}\norm F\,\hs{Q(\psi(t))}$, so
\[
\norm{FQ}_{\E_d}^2\le\tfrac{\rho^2}4\norm{FQ}_\partial^2\le\tfrac{\rho^2}4\max_{\overline\Om}\norm F^2\,\norm Q_\partial^2\le\max\big\{\tfrac{\rho^2}2,1+V(\Om)^2\big\}\max_{\overline\Om}\norm F^2\,\norm Q_{\E_d}^2,
\]
and likewise for $QF$. Matrix polynomials are dense in $\E_d$, and convergence there implies pointwise convergence; so $m_\ell(F),m_r(F)\le\max\{\frac\rho{\sqrt2},\sqrt{1+V(\Om)^2}\}\max_{\overline\Om}\norm F$ for matrix polynomials $F$. Now let $F:\Om\to M_d$ be bounded and holomorphic. Let $\varphi$ be a conformal map of $\D$ onto $\Om$; since $\partial\Om$ is a regular real-analytic Jordan curve, $\varphi$ extends to a conformal map of a disk $\{|w|<1+\delta\}$ \cite[Proposition~3.1]{Pom92}. For $0<\tau<1$ put $F_\tau(z)=F(\varphi(\tau\varphi^{-1}(z)))$. Then $F_\tau$ is holomorphic near $\overline\Om$ with $\sup_{\overline\Om}\norm{F_\tau}\le\sup_\Om\norm F$, hence a uniform limit on $\overline\Om$ of matrix polynomials, by Runge's theorem, since $\C\setminus\overline\Om$ is connected; and $F_\tau\to F$ pointwise on $\Om$ as $\tau\uparrow1$. So there are matrix polynomials $G_j\to F$ pointwise on $\Om$ with $\limsup_j\max_{\overline\Om}\norm{G_j}\le\sup_\Om\norm F$. For $Q\in\E_d$ the sequence $(G_jQ)$ is then bounded in $\E_d$, so a subsequence converges weakly; since point evaluations are continuous, its limit is $FQ$, whose norm is at most $\limsup_j\norm{G_jQ}_{\E_d}$. Arguing likewise for $QF$, we obtain $m_\ell(F),m_r(F)\le\max\{\frac\rho{\sqrt2},\sqrt{1+V(\Om)^2}\}\sup_\Om\norm F$, which is \eqref{eq:convexmult}.
\end{proof}

\begin{theorem}\label{thm:convex}
Let $\Om$ be a bounded Jordan domain with regular real-analytic boundary, and let $\rho\ge\sqrt2$. For every $T\in\Kc_\rho(\Om)$, every $d\ge1$ and every bounded holomorphic $F:\Om\to M_d$,
\[
\norm{F(T)}\le\max\big\{\rho,\sqrt{2+2V(\Om)^2}\big\}\sup_\Om\norm F .
\]
For $\rho\ge2$, the constant is $\rho$ if $V(\Om)\le1$, in particular if $\Om$ is convex or $\sum_{k\ge1}k|\psi_k|\le\frac12$, with $\psi_k$ as in Lemma~\ref{lem:rhoOm}(d); for $\Om=\D$, it is $\rho$ for every $\rho\ge\sqrt2$, and it is optimal.
\end{theorem}

\begin{proof}
The family $(r_{t,\rho})$ satisfies the hypotheses of Proposition~\ref{prop:pencil}, by Lemma~\ref{lem:convexclass}, since $(1-\frac2\rho)+\frac2\rho\cdot1=1$; hence $\Psi_T$ is accretive for every $T\in\Kc_\rho(\Om)$, and $\E_{\Om,\rho}$ is nondegenerate, since it contains $1$. So Corollary~\ref{cor:main} and Lemma~\ref{lem:faber} give the bound; the maximum is $\rho$ if $\rho\ge2$ and $V(\Om)\le1$, which holds in the two cases named by Lemma~\ref{lem:rhoOm}(c), (d), and if $\Om=\D$, where $V(\Om)=0$. For $\Om=\D$, $T=\rho N$, with $N=\big(\begin{smallmatrix}0&1\\0&0\end{smallmatrix}\big)$, lies in $\Kc_\rho(\D)$, since $r_t(\rho N)=I+\bar t\rho N$ and $\norm{\Rp(\bar tN)}=\frac12$, and $F(z)=z$ gives equality.
\end{proof}

\begin{corollary}\label{cor:convexclosed}
Let $\Om$ and $\rho$ be as in Theorem~\ref{thm:convex}, and let $T\in\B(H)$.
\begin{enumerate}[(a)]
\item If $T$ has a Faber $\rho$-dilation, then $\sigma(T)\subseteq\overline\Om$, and $\overline\Om$ is a complete $C_{\Om,\rho}$-spectral set for $T$, where $C_{\Om,\rho}=\max\{\rho,\sqrt{2+2V(\Om)^2}\}$, which, for $\rho\ge2$, is $\rho$ if $V(\Om)\le1$, in particular if $\Om$ is convex.
\item If $\Om$ is convex and $W(T)\subseteq\overline\Om$, then $\overline\Om$ is a complete $2$-spectral set for $T$.
\end{enumerate}
\end{corollary}

\begin{proof}
For $s>1$ let $\Om_s$ be the bounded component of $\C\setminus\psi(s\T)$, which contains $\overline\Om$. Its boundary $\psi(s\T)$ is a regular real-analytic Jordan curve, since $\psi$ is univalent on $\{|t|>\tau_0\}$, its exterior conformal map is $t\mapsto\psi(st)$, its functions $r_t$ are the $r_{st}$, its Faber polynomials are $s^{-n}\Phi_n$, and its function $S$ is $(t,u)\mapsto S(st,su)$. As $s\downarrow1$, $S(st,su)\to S(t,u)$ uniformly on $\T\times\T$; so the kernels $k_\partial$ of the domains $\Om_s$ tend to that of $\Om$ uniformly, and their operators $\Theta$ tend to that of $\Om$ in Hilbert--Schmidt norm, whence $V(\Om_s)\to V(\Om)$.

(a) If $T$ has a Faber $\rho$-dilation $\Phi_n(T)=\rho P_HU^n|_H$, then $\norm{\Phi_n(T)}\le\rho$ for $n\ge1$. Since the $c^n\Phi_n$ are the monic Faber polynomials of $\overline\Om$ and $c$ is its capacity, the spectral inclusion theorem of Atzmon, Eremenko and Sodin \cite[Theorem~3 and pp.~727--728]{AES99} gives $\sigma(T)\subseteq\overline\Om\subset\Om_s$. By Sz.-Nagy's dilation theorem \cite{SzNF} the contraction $s^{-1}U$ has a unitary dilation $\widetilde U$, so $s^{-n}\Phi_n(T)=\rho P_H\widetilde U^n|_H$ is a Faber $\rho$-dilation for $\Om_s$, and $T\in\Kc_\rho(\Om_s)$ by Lemma~\ref{lem:rhoconvex}(a). Now let $F$ be $M_d$-valued and holomorphic on a neighbourhood $\mathcal O$ of $\overline\Om$. Since $\bigcap_{s>1}\overline{\Om_s}=\overline\Om$, we have $\overline{\Om_s}\subset\mathcal O$ for $s$ close to $1$, and Theorem~\ref{thm:convex} for $\Om_s$ gives $\norm{F(T)}\le\max\{\rho,\sqrt{2+2V(\Om_s)^2}\}\sup_{\Om_s}\norm F$, which tends to $C_{\Om,\rho}\max_{\overline\Om}\norm F$ as $s\downarrow1$.

(b) Let $z_0\in\Om$ and $0<\tau<1$, and put $T_\tau=z_0+\tau(T-z_0)$. Since $\Om$ is convex and open, $W(T_\tau)=z_0+\tau(W(T)-z_0)$ lies in the compact set $z_0+\tau(\overline\Om-z_0)\subset\Om$, so $\sigma(T_\tau)\subset\Om$ and $T_\tau\in\Kc_2(\Om)$ by Lemma~\ref{lem:convexclass}. For $F$ holomorphic on a neighbourhood of $\overline\Om$, Theorem~\ref{thm:convex} with $\rho=2$ gives $\norm{F(T_\tau)}\le2\max_{\overline\Om}\norm F$, and $F(T_\tau)\to F(T)$ in norm as $\tau\uparrow1$, since $\sigma(T)\subseteq\overline{W(T)}\subseteq\overline\Om$.
\end{proof}

For convex $\Om$ and $\rho=2$, Theorem~\ref{thm:convex} and Corollary~\ref{cor:convexclosed} give the complete Crouzeix inequality of \cite{OAI,OAI2} for these domains, and with the outer approximation of \cite[Section~4]{OAI} it extends to arbitrary compact convex sets. Both steps of the proof of \cite[Theorem~3.2]{OAI} (Remark~\ref{rem:engine}) rest on the positivity of $\Rp r_t$ for convex $\Om$: the norm estimates for Faber series on $k_\partial\ge0$, and the comparison step on $N_T(t)=2\Rp r_t(T)\ge0$, which follows from $W(T)\subseteq\overline\Om$ (Lemma~\ref{lem:convexclass}). For nonconvex $\Om$, $\Rp r_t$ takes negative values for some $t\in\T$ (Lemma~\ref{lem:rhoOm}(a)). Theorem~\ref{thm:convex} assumes the lower bound $\Rp r_t(T)\ge(1-\frac\rho2)I$ directly, which gives accretivity by Proposition~\ref{prop:pencil}, and \eqref{eq:faber} replaces the sign of $k_\partial$ by the norm $V(\Om)$, which is at most $1$ also for some nonconvex $\Om$ (Lemma~\ref{lem:rhoOm}(d) and the example after it). Besides the disk, the constant of Theorem~\ref{thm:convex} is optimal for ellipses and $\rho=2$, since the scalar constant already equals $2$, by \cite[Theorem~6]{Jin}; for $\rho\le2$ and $V(\Om)>1$ it is not (Remark~\ref{rem:rho2}). We have not attempted to determine whether it is optimal in the other cases.

\begin{remark}\label{rem:rho2}
For $\rho=2$, the complete Crouzeix inequality \cite[Corollary~6.3]{OAI2} gives the constant $2$ in Theorem~\ref{thm:convex} for every $\Om$, without the condition $V(\Om)\le1$. Let $T\in\B(H)$ with $\sigma(T)\subset\Om$, let $t\in\T$ and put $\zeta=\psi(t)$. Since $1/r_t$ is affine, $r_t(T)^{-1}=(1/r_t)(T)$, and since inverses of invertible accretive operators are accretive, $\Rp r_t(T)\ge0$ if and only if $W(T)$ lies in the closed half-plane $\{\Rp(1/r_t)\ge0\}=\{z:\Rp(\overline{n(\zeta)}(z-\zeta))\le0\}$. Each of these half-planes contains $\overline\Om$ if $\Om$ is convex. In general, their intersection over $t\in\T$ is contained in $\overline\Om$: if $z\notin\overline\Om$ and $\zeta$ is a point of $\overline\Om$ nearest to $z$, then $\zeta=\psi(t)$ for some $t\in\T$ and $z-\zeta=|z-\zeta|\,n(\zeta)$, since $\partial\Om$ is smooth, so $z$ does not lie in the half-plane for this $t$. Hence $\overline{W(T)}\subseteq\overline\Om$ for every $T\in\Kc_2(\Om)$, and \cite[Corollary~6.3]{OAI2}, applied to the functions $F_\tau$ of the proof of Lemma~\ref{lem:faber}, which tend to $F$ locally uniformly on $\Om$ as $\tau\uparrow1$, gives $\norm{F(T)}\le2\sup_\Om\norm F$ for every bounded holomorphic $F:\Om\to M_d$. Likewise, if $T$ has a Faber $2$-dilation, then $T\in\Kc_2(\Om_s)$ for every $s>1$, by the proof of Corollary~\ref{cor:convexclosed}(a), so $\overline{W(T)}\subseteq\bigcap_{s>1}\overline{\Om_s}=\overline\Om$, and $\overline\Om$ is a complete $2$-spectral set for $T$. Since $\Kc_\rho(\Om)\subseteq\Kc_2(\Om)$ for $\rho<2$, the bound $2$ also holds for $\sqrt2\le\rho<2$, where it is smaller than the constant of Theorem~\ref{thm:convex} exactly when $V(\Om)>1$. This uses $\rho\le2$: for $\rho>2$, $\Kc_\rho(\D)$ contains $\rho N$, with $N$ as in the proof of Theorem~\ref{thm:convex}, and $W(\rho N)\not\subseteq\overline\D$; if $\Om$ is not convex and $\rho\ge\rho_\Om$, $\Kc_\rho(\Om)$ contains normal operators $T$ with $W(T)\not\subseteq\overline\Om$.
\end{remark}

\begin{remark}\label{rem:rho}
For $\Om=\D$, Faber $\rho$-dilations are unitary $\rho$-dilations, and Theorem~\ref{thm:convex} and Corollary~\ref{cor:convexclosed}(a) show that the closed unit disk is a complete $\rho$-spectral set for every $\rho$-contraction with $\rho\ge\sqrt2$. By Paulsen's theorem \cite[p.~226]{Pa84b}, this is equivalent to the case $\rho\ge\sqrt2$ of the theorem of Okubo and Ando \cite{OA} that every $\rho$-contraction is similar to a contraction by an invertible $G$ with $\norm G\norm{G^{-1}}\le\rho$; they proved it for every $\rho\ge1$, and for $1\le\rho<\sqrt2$ our method does not reach this constant (Remark~\ref{rem:floor}). For refinements and other proofs, see \cite{Dr08,SdV,HdV,HM}, and \cite[Section~7.2]{Jin} for the scalar constant for matrices and $\rho\ge2$; for relations with operator radii, see \cite{BCK18}.
\end{remark}

\begin{remark}\label{rem:convex}
For $t\in\T$, $r_t(T)=t\psi'(t)\,(\psi(t)-T)^{-1}$ is the resolvent family of \cite[Theorem~3.1]{OAI}, with $h=\psi$ and $\lambda=t$. For convex $\Om$, the positivity of $N_T(t)=r_t(T)+r_t(T)^*$, that is, $T\in\Kc_2(\Om)$, is proved there for matrices $T$ with $W(T)\subset\Om$, and in general in \cite[Proposition~3.4]{SdVIEOT} (Lemma~\ref{lem:convexclass}). For unit vectors $x,y\in\C^d\otimes H$ with $F(T)x=\norm{F(T)}y$, the proof of \cite[Theorem~3.2]{OAI} ends with a trace inequality which, in the notation of Section~\ref{sec:engine}, reads $\ip{P_x-P_y}{\Co{N_T(t)}(P_x-P_y)}_{\mathrm{HS}}\ge0$; integrated over $\T$, by \eqref{eq:pencil} and Lemma~\ref{lem:momentineq}, it becomes
\[
2\ip{P_x-P_y}{(\Rp\Psi_T)(P_x-P_y)}_{\mathrm{HS}}=2\big(J(x,x)^2+J(y,y)^2-J(x,y)^2-J(y,x)^2\big)\ge0,
\]
a special case of the accretivity of the moment map $\Psi_T$ of $\E_{\Om,2}$. Replacing $x$ by $cx$ and optimizing over $c>0$ gives the inequality of Lemma~\ref{lem:momentineq}. The first step of the proof of \cite[Theorem~3.2]{OAI} uses that multiplication by an $F$ with $\sup_\Om\norm F\le1$ is a contraction from the second of the two norms of Remark~\ref{rem:engine} to the first; the norm of $\E_{\Om,2}$ is $1/\sqrt2$ times the second, and, for convex $\Om$, Lemma~\ref{lem:faber} gives instead $M(\E_{\Om,2})\le\sqrt2$ (see also Remark~\ref{rem:twospaces}).
\end{remark}

\begin{remark}\label{rem:CG}
By \eqref{eq:dlp}, the class $\Kc_\rho(\Om)$ is defined by the hypothesis $\mu(\zeta,T)\ge\lambda(\zeta)I$ $(\zeta\in\partial\Om)$ of \cite[Theorem~1]{CG26}, with $\lambda(\psi(t))=-(\rho-2)/(2\pi|\psi'(t)|)$, for which $\int_{\partial\Om}\lambda\,ds=2-\rho$. We compare Theorem~\ref{thm:convex} with that theorem and with earlier bounds on $\Kc_\rho(\Om)$.

(a) Scalar bounds. The theorem of Crouzeix and Greenbaum \cite[Theorem~1]{CG26}, which rests on a lemma of Lorist and Schwenninger in the form given in \cite[Remark~2(ii)]{LS}, needs no convexity and gives the scalar bound $\max\{\rho,1+\sqrt{\rho/2}\}$. For $\rho\ge2$ this is $\rho$, so for scalar functions Theorem~\ref{thm:convex} follows from it. For $\sqrt2\le\rho<2$ it is $1+\sqrt{\rho/2}>\rho$, and the constant of Theorem~\ref{thm:convex} is smaller when $2+2V(\Om)^2<(1+\sqrt{\rho/2})^2$, in particular when $V(\Om)\le1/\sqrt2$, since $1+\sqrt{\rho/2}\ge1+2^{-1/4}>\sqrt3$. For an ellipse with semi-axes $a\ge b$, $\psi(t)=c(t+qt^{-1})$ with $q=(a-b)/(a+b)$, so $S(t,u)=\sum_{n\ge1}q^nt^{-n}u^{-n}$ and $V(\Om)=q$ by Lemma~\ref{lem:rhoOm}(c); thus $V(\Om)\le1/\sqrt2$ when $a/b\le3+2\sqrt2$. Earlier scalar bounds in terms of the negative part of $\mu$ are in \cite{CGL18,CG19}, and for Faber polynomials in \cite{BC14}.

(b) Complete bounds. Extended to matrix-valued functions as in \cite[Section~2]{CK25}, the argument of \cite[Theorem~2]{CG19} gives the complete bound
\[
\frac\rho2+\sqrt{\frac{\rho^2}4+|\rho-2|+\kappa_\Om},\qquad\kappa_\Om=\max_{u\in\T}\int_\T|k_\partial(t,u)|\,d\omega(t);
\]
in the notation of \cite[Theorem~2]{CG19}, one can take $c_2=\rho/2$ and $\hat\gamma\le|\rho-2|$, and $c_1\le\kappa_\Om$. Since $\int_\T k_\partial(t,u)\,d\omega(t)=1$, $\kappa_\Om\ge1$, with equality exactly when $k_\partial\ge0$, that is, for convex $\Om$ (Lemma~\ref{lem:rhoOm}(a)). This bound is larger than $\rho$, and larger than $2$ for $\sqrt2\le\rho<2$. So the constant of Theorem~\ref{thm:convex} is smaller when $V(\Om)\le1$ or $\rho\ge\sqrt{2+2V(\Om)^2}$, and, for $\rho\ge2$, also when $V(\Om)^2\le\frac12+\sqrt2$: if $V(\Om)>1$, then $\Om$ is not convex (Lemma~\ref{lem:rhoOm}(c)), so $\kappa_\Om>1$ and the bound exceeds $1+\sqrt2$. Otherwise the two constants are not directly comparable: $V(\Om)$ is the norm of $\Theta$ on the closed span of the functions $t^n$, $n\ge1$, in $L^2(\T)$, while $\kappa_\Om$ is the norm of the same integral operator on $L^\infty(\T)$.
\end{remark}

\section{Examples}\label{sec:examples}

Sections~\ref{ssec:schrodinger} and~\ref{ssec:divform} apply Corollary~\ref{cor:unbounded} and Remark~\ref{rem:hinfty} to differential operators with numerical range in a sector, and Section~\ref{ssec:hyponormal} applies the results of Sections~\ref{sec:twodisks} and~\ref{sec:convex} to hyponormal operators.

\subsection{Schr\"odinger operators with bounded complex potentials}\label{ssec:schrodinger}
Let $V\in L^\infty(\RR^n)$ be complex-valued, let $A=-\Delta+V$ in $L^2(\RR^n)$, with domain the Sobolev space $W^{2,2}(\RR^n)$, and put $c_V=\norm{\Ip V}_\infty-\operatorname*{ess\,inf}\Rp V$. For $u$ in the domain, $\ip{u}{Au}=\norm{\nabla u}^2+\int V|u|^2$, so
\[
\Rp\ip{u}{(A+c_V)u}\ge\norm{\Ip V}_\infty\norm u^2\ge|\Ip\ip{u}{(A+c_V)u}|,
\]
that is, $W(A+c_V)\subseteq\Sbar{\pi/4}$. For $\mu>\max\{0,\norm V_\infty-c_V\}$, $A+c_V+\mu=(-\Delta+c_V+\mu)\big(I+(-\Delta+c_V+\mu)^{-1}V\big)$ with $\norm{(-\Delta+c_V+\mu)^{-1}V}\le\norm V_\infty/(c_V+\mu)<1$, so $-\mu\in\operatorname{res}(A+c_V)$; since $\C\setminus\Sbar{\pi/4}$ is connected, \cite[Theorem~V.3.2]{KaP} gives $\C\setminus\Sbar{\pi/4}\subseteq\operatorname{res}(A+c_V)$. Since $\Sbar{\pi/4}+t\subseteq\Sbar{\pi/4}$ for $t\ge0$, both inclusions also hold for $A+c$ with $c\ge c_V$. Moreover, $A+c$ is injective: if $(A+c)u=0$, then $\norm{\nabla u}^2\le\Rp\ip{u}{(A+c)u}=0$, so $u$ is constant, hence $u=0$. Hence Corollary~\ref{cor:unbounded} and Remark~\ref{rem:hinfty} apply to $A+c$, $c\ge c_V$, with $\alpha=\pi/4$: for every $d\ge1$,
\[
\norm{F(A+c)}\le\sqrt2\,\sup_{\Sbar{\pi/4}}\norm F
\]
for every $M_d$-valued rational $F$ bounded on $\Sbar{\pi/4}$, and
\[
\norm{F(A+c)}\le\sqrt2\,\sup_{S_{\pi/4}}\norm F
\]
for every $\theta\in(\pi/4,\pi)$ and every bounded holomorphic $F:S_\theta\to M_d$. The constant does not depend on $V$, which enters only through $c_V$; for instance, the bounds hold for $A+c$ whenever $\Rp V\ge0$ and $\norm{\Ip V}_\infty\le c$. For real $V$, $A+c_V$ is self-adjoint and nonnegative, and the spectral theorem gives the constant $1$. Finally, with $b(z)=(z-1)/(z+1)$, the closed lens $L=b(\Sbar{\pi/4}\cup\{\infty\})\subset\overline\D$, with vertices $\pm1$ and vertex angle $\pi/2$, is a complete $\sqrt2$-spectral set for the bounded operator $B=b(A+c_V)=I-2(A+c_V+1)^{-1}$. Indeed, for $\lambda\notin L$ the point $z=(1+\lambda)/(1-\lambda)$ lies outside $\Sbar{\pi/4}$, and
\[
(\lambda-B)^{-1}=\frac1{\lambda-1}\big(I+(1+z)(A+c_V-z)^{-1}\big),
\]
so $\sigma(B)\subseteq L$; and for every $M_d$-valued rational $F$ with poles off $L$, $F\circ b$ is rational and bounded on $\Sbar{\pi/4}$, with $\sup_{\Sbar{\pi/4}}\norm{F\circ b}=\max_L\norm F$, and $F(B)=(F\circ b)(A+c_V)$.

\subsection{Divergence form operators with complex coefficients}\label{ssec:divform}
Let $\pi/4\le\alpha<\pi/2$ and $\eta>0$, and let $a:\RR^n\to M_n$ be bounded and measurable, with $\ip{\xi}{a(x)\xi}\in\Sbar\alpha$ and $\Rp\ip{\xi}{a(x)\xi}\ge\eta|\xi|^2$ for all $\xi\in\C^n$ and almost every $x\in\RR^n$. The form
\[
\mathfrak a[u,v]=\int_{\RR^n}\ip{\nabla u}{a\nabla v}\,dx\qquad(u,v\in W^{1,2}(\RR^n))
\]
is densely defined, sectorial and, since $\eta\norm{\nabla u}^2\le\Rp\mathfrak a[u,u]\le\norm a_\infty\norm{\nabla u}^2$, closed, and the operator $A=-\operatorname{div}(a\nabla)$ associated with it by Kato's first representation theorem \cite[Theorem~VI.2.1]{KaP} is m-sectorial, with $\ip{u}{Au}=\mathfrak a[u,u]$ for $u$ in its domain; in particular, $-\mu\in\operatorname{res}(A)$ for large $\mu>0$. Since $\ip{\nabla u}{a\nabla u}\in\Sbar\alpha$ almost everywhere and $\Sbar\alpha$ is a closed convex cone, $W(A)\subseteq\Sbar\alpha$, and so, since $\C\setminus\Sbar\alpha$ is connected and contains $-\mu$, $\C\setminus\Sbar\alpha\subseteq\operatorname{res}(A)$ by \cite[Theorem~V.3.2]{KaP}. If $Au=0$, then $\eta\norm{\nabla u}^2\le\Rp\mathfrak a[u,u]=0$, so $u$ is constant, hence $u=0$; thus $A$ is injective. Corollary~\ref{cor:unbounded} and Remark~\ref{rem:hinfty} give, for every $d\ge1$,
\[
\norm{F(A)}\le\sqrt2\,\sup_{\Sbar\alpha}\norm F
\]
for every $M_d$-valued rational $F$ bounded on $\Sbar\alpha$, and
\[
\norm{F(A)}\le\sqrt2\,\sup_{S_\alpha}\norm F
\]
for every $\theta\in(\alpha,\pi)$ and every bounded holomorphic $F:S_\theta\to M_d$. The constant depends neither on $n$, nor on $\eta$, nor on $a$ beyond $\alpha$; for estimates of the least $\alpha$ for a given $a$, see \cite{MR25}. For these operators, the scalar bound with the supremum over $S_\theta$ and the constant $2$ follows from \cite[Theorem~3.13]{MR25} and \cite{Jin,LS}; for the corresponding operators on open sets with mixed boundary conditions, it holds after a shift by any $\delta>0$, and without the shift if the form is coercive \cite[Theorem~3.15]{MR25}. For the other known bounds, among them the bound \eqref{eq:crouzeixbound} of \cite{Cr25}, which is smaller for $\alpha>\alpha_*$, see Remark~\ref{rem:hinfty}. Corollary~\ref{cor:unbounded} applies in the same way to the operators associated with other densely defined closed forms whose numerical range lies in $\Sbar\alpha$, for instance with Dirichlet, Neumann or mixed boundary conditions on open sets, and Remark~\ref{rem:hinfty} applies to those that are injective.

\subsection{Hyponormal operators}\label{ssec:hyponormal}
An operator $T\in\B(H)$ is \emph{hyponormal} if $T^*T\ge TT^*$. Hyponormal operators are normaloid, that is, $\norm T=\max\{|z|:z\in\sigma(T)\}$ \cite[Proposition~III.1.3]{MP}, and the inverse of an invertible hyponormal operator is hyponormal \cite[Proposition~III.2.1]{MP}. Since $aT+b$ is hyponormal for $a,b\in\C$, and every M\"obius transformation $m$ is affine or of the form $m(z)=a+b(z-c)^{-1}$, it follows that $m(T)$ is hyponormal whenever the pole of $m$ is not in $\sigma(T)$.

\begin{proposition}\label{prop:hyponormal}
Let $T\in\B(H)$ be hyponormal.
\begin{enumerate}[(a)]
\item $\overline{W(T)}=\operatorname{conv}\sigma(T)$; hence, if $r$ is affine and $\Rp r\ge0$ on $\sigma(T)$, then $\Rp r(T)\ge0$.
\item Let $b$ be a M\"obius transformation that maps an open disk $D$ of the Riemann sphere onto $\D$. If $\sigma(T)\subseteq\overline D$, then $\norm{b(T)}\le1$.
\item Let $\Om$ be as in Section~\ref{sec:convex}, with $\sigma(T)\subset\Om$. Then $T\in\Kc_\rho(\Om)$ for every $\rho\ge\rho_\Om$.
\end{enumerate}
\end{proposition}

\begin{proof}
(a) The first statement is classical \cite[Theorem~2]{St65}, \cite[Exercise~III.12]{MP}, and the second follows, since $W(r(T))=r(W(T))$ for affine $r$.

(b) Since $b(\overline D)=\overline\D$, the pole of $b$ is not in $\overline D\supseteq\sigma(T)$. So $b(T)$ is hyponormal, hence normaloid, and $\norm{b(T)}=\max_{z\in\sigma(T)}|b(z)|\le1$ by the spectral mapping theorem; compare the paragraph before \cite[Corollary~8]{DEY17}.

(c) Let $\rho\ge\rho_\Om$ and $t\in\T$. The M\"obius transformation $r_{t,\rho}$ has positive real part on $\Om$ (Lemma~\ref{lem:rhoOm}(b)), so $\sigma(T)\subset\Om\subseteq D_t=\{\Rp r_{t,\rho}>0\}$, an open disk of the Riemann sphere, which the M\"obius transformation $(r_{t,\rho}-1)/(r_{t,\rho}+1)$ maps onto $\D$. By (b) and \eqref{eq:cayleyfact}, $\Rp r_{t,\rho}(T)\ge0$. So $T\in\Kc_\rho(\Om)$, the class of $(r_{t,\rho})_{t\in\T}$.
\end{proof}

Thus a hyponormal operator satisfies the disk norm condition for every open disk of the Riemann sphere whose closure contains its spectrum, and the bounds of Sections~\ref{sec:twodisks} and~\ref{sec:convex} apply to it. For $\Om$ as in Section~\ref{sec:convex} put
\begin{equation}\label{eq:COm}
C_\Om=\max\big\{\rho_\Om,\sqrt{2+2V(\Om)^2}\big\},
\end{equation}
which is the constant $C_{\Om,\rho}$ of Corollary~\ref{cor:convexclosed} for $\rho=\rho_\Om$; so $C_\Om=\rho_\Om$ if $V(\Om)\le1$, and $C_\Om=2$ if $\Om$ is convex.

\begin{corollary}\label{cor:hyponormal}
Let $T\in\B(H)$ be hyponormal, and let $\Om$ and $K$ be one of the following:
\begin{enumerate}[(i)]
\item $\Om=D_1\cap D_2$ as in Corollary~\ref{cor:lens}, with angle $2\alpha$ at the crossing points, and $K=c_\alpha$ (Theorem~C);
\item $\Om$ as in Section~\ref{sec:convex}, and $K=C_\Om$.
\end{enumerate}
If $\sigma(T)\subset\Om$, then, for every $d\ge1$ and every bounded holomorphic $F:\Om\to M_d$,
\[
\norm{F(T)}\le K\sup_\Om\norm F .
\]
If $\sigma(T)\subseteq\overline\Om$, then $\overline\Om$ is a complete $K$-spectral set for $T$.
\end{corollary}

\begin{proof}
(i) Each $b_i$ is a M\"obius transformation, so $\norm{b_i(T)}\le1$ by Proposition~\ref{prop:hyponormal}(b), and Corollary~\ref{cor:lens} gives both statements.

(ii) The first statement is Proposition~\ref{prop:hyponormal}(c) with $\rho=\rho_\Om$ and Theorem~\ref{thm:convex}. For the second, let $F$ be holomorphic on a neighbourhood of $\overline\Om$, and let $\Om_s$, $s>1$, be as in the proof of Corollary~\ref{cor:convexclosed}. Then $\sigma(T)\subseteq\overline\Om\subset\Om_s$, and $C_{\Om_s}\to C_\Om$ as $s\downarrow1$, since $V(\Om_s)\to V(\Om)$ and, by Lemma~\ref{lem:rhoOm}(a) and the uniform convergence of the kernels $k_\partial$, $\rho_{\Om_s}\to\rho_\Om$. For $s$ so close to $1$ that $F$ is holomorphic near $\overline{\Om_s}$, the first statement for $\Om_s$ gives $\norm{F(T)}\le C_{\Om_s}\sup_{\Om_s}\norm F\to C_\Om\max_{\overline\Om}\norm F$.
\end{proof}

\begin{remark}\label{rem:hyponormal}
(a) The spectrum of a hyponormal operator need not be a spectral set \cite{Cl70,Wa73,Ha82}, whereas for subnormal operators it is a complete spectral set \cite[Corollary~I.3.5]{MP}. Let $\Om=D_1\cap D_2$ be as in Theorem~B, bounded and with disjoint boundary circles. Then $\overline\Om$ is a complete $\sqrt2$-spectral set for every hyponormal $T$ with $\sigma(T)\subseteq\overline\Om$, and $\sqrt2$ cannot be lowered; the bound follows from a theorem of Bello and Yakubovich \cite[Theorem~4.2]{BY23} (see also \cite{Ts22}), as follows. The M\"obius transformation of Lemma~\ref{lem:mobius} maps $\overline\Om$ onto some $\overline{\Ann_R}$, with its pole off $\overline\Om$, since it maps both boundary circles onto circles centred at $0$, and, as its pole is not in $\sigma(T)$, it preserves hyponormality; so it suffices to take $\Om=\Ann_R$. For hyponormal $T$ with $\sigma(T)\subseteq\overline{\Ann_R}$, put $B=T/R$, $r=R^{-2}$ and $P=BB^*$. Since $T$ and $T^{-1}$ are normaloid, $r^2I\le P\le I$, so $r^2P^{-1}+P\le(1+r^2)I$, as $s+r^2/s\le1+r^2$ for $r^2\le s\le1$; and $B^{-1}(B^{-1})^*=(B^*B)^{-1}\le P^{-1}$, since $B^*B\ge P$. Hence
\[
r^2B^{-1}(B^{-1})^*+BB^*\le(1+r^2)I,
\]
and multiplying by $B$ on the left and by $B^*$ on the right shows that $X=B^*$ satisfies
\[
-X^{*2}X^2+(1+r^2)X^*X-r^2I\ge0.
\]
For invertible $X$ with this property, $\{r\le|z|\le1\}$ is a complete $\sqrt2$-spectral set \cite[Theorem~4.2]{BY23}; as this set is symmetric under complex conjugation, it is one for $B$ as well. For optimality, we rescale the example in the proof of \cite[Proposition~3.1]{Ts22}: the bilateral weighted shift $S\delta_n=w_n\delta_{n+1}$ $(n\in\Z)$ on $\ell^2(\Z)$, with $(\delta_n)$ the standard basis, $w_n=R^{-1}$ for $n<0$ and $w_n=R$ for $n\ge0$, is hyponormal, since its weights increase, with $\sigma(S)=\overline{\Ann_R}$ \cite{Sh}; and $f_n(z)=R^{-n}(z^n+z^{-n})$ satisfies $f_n(S)\delta_0=\delta_n+\delta_{-n}$, so $\norm{f_n(S)}\ge\sqrt2$, while $\max_{\overline{\Ann_R}}|f_n|=1+R^{-2n}\to1$. If instead $\Om$ is bounded with tangent boundary circles, then $\Om=\{|z-c_1|<R_1\}\cap\{|z-c_2|>R_2\}$ with $|c_1-c_2|=R_1-R_2$; for $0<\epsilon<R_2$ the circles $|z-c_1|=R_1+\epsilon$ and $|z-c_2|=R_2-\epsilon$ are disjoint, as $|c_1-c_2|+R_2-\epsilon<R_1+\epsilon$, and the bound $\sqrt2$ on the closures of the resulting domains, which contain $\overline\Om$ and decrease to it, passes to $\overline\Om$.

(b) We compare Corollary~\ref{cor:hyponormal} with known bounds. For $\Om$ as in Corollary~\ref{cor:hyponormal}(ii), Dritschel, Est\'evez and Yakubovich showed, also under more general hypotheses, that $\overline\Om$ is a complete $K$-spectral set for every hyponormal $T$ with $\sigma(T)\subseteq\overline\Om$, for some $K$ that is not explicit \cite[Theorem~2, Corollary~8]{DEY17}. Explicit complete constants follow from \cite[Theorem~2 and Lemma~7]{CG19} with the argument of \cite[Section~2]{CK25} (compare Remark~\ref{rem:CG}), and scalar constants at most $\rho_\Om$ from \cite[Theorem~1 and Lemma~2(3)]{CG26}; $C_\Om$ equals $\rho_\Om$ when $V(\Om)\le1$. For convex $\Om$, Corollary~\ref{cor:hyponormal}(ii) also follows from \cite[Corollary~6.3]{OAI2}, since $\overline{W(T)}=\operatorname{conv}\sigma(T)$ (Proposition~\ref{prop:hyponormal}(a)), which is then even a complete $2$-spectral set. For lenses and crescents of vertex angle $2\alpha$, the complete version of the bound \eqref{eq:crouzeixbound} of \cite{Cr25}, with Proposition~\ref{prop:hyponormal}(b) and Lemma~\ref{lem:mobius}, improves on the constant $c_\alpha$ of Corollary~\ref{cor:hyponormal}(i) exactly when $\alpha>\alpha_*$.
\end{remark}

\appendix

\section{Optimality of the space for acute sectors}\label{app:spaceopt}

We prove Proposition~\ref{prop:spaceopt} in the stronger form $m(k)\ge\sqrt{2-4\alpha/\pi}$, where $m(k)$ is the supremum of the norms of the multiplication operators $f\mapsto Ff$ on $\E_k$ over all holomorphic $F:S_\alpha\to\C$ with $\sup_{S_\alpha}|F|\le1$, so that $m(k)\le M(\E_k)$ (Definition~\ref{def:mult}); $\tilde k_{S_\alpha}$ is as in Section~\ref{ssec:sector}. We call $A\in M_2$ \emph{admissible} if $\sigma(A)\subset S_\alpha$ and $W(A)\subseteq\Sbar\alpha$. By \cite[Theorem~5.21]{PR}, $m(k)\le m$ if and only if $(m^2-F(z)\overline{F(w)})\,k(z,w)$ is a positive kernel for every holomorphic $F:S_\alpha\to\C$ with $\sup_{S_\alpha}|F|\le1$.

\begin{proof}[Proof of Proposition~\ref{prop:spaceopt}]
We may assume that $m:=m(k)<\infty$. Let $\mathcal C$ be the set of holomorphic kernels $k'$ on $S_\alpha$ with $k'(1,1)=1$, $m(k')\le m$ and accretive moment maps at all admissible $A\in M_2$. By \cite[Theorem~5.21]{PR} and Lemma~\ref{lem:cauchy}, $\mathcal C$ is convex.

Every nonzero holomorphic kernel $k'$ whose moment maps at admissible $A\in M_2$ are accretive satisfies $\Rp k'>0$ on $S_\alpha\times S_\alpha$. Indeed, for $z,w\in S_\alpha$, $A=\mathrm{diag}(z,w)$ is admissible, as $W(A)$ is the segment $[z,w]$, and $\Psi_A(N)=k'(z,w)N$ for $N=\big(\begin{smallmatrix}0&1\\0&0\end{smallmatrix}\big)$; so $\Rp k'(z,w)\ge0$. The function $\Phi(z,v)=k'(z,\bar v)$ is holomorphic on $S_\alpha\times S_\alpha$, by Hartogs' theorem \cite[Theorem~2.2.8]{Ho}, and $\Rp\Phi\ge0$. If $\Rp\Phi$ vanished at a point, it would vanish identically, by the minimum principle for harmonic functions, and $\Phi$ would be an imaginary constant, hence $0$, since $k'(z,z)\ge0$. In particular $k(1,1)>0$, and $k/k(1,1)\in\mathcal C$.

The set $\mathcal C$ is compact in the Fr\'echet space of functions on $S_\alpha\times S_\alpha$ that are holomorphic in the first variable and antiholomorphic in the second, with locally uniform convergence. Indeed, for $k'\in\mathcal C$ and $\Phi(z,v)=k'(z,\bar v)$, the function $(\Phi-1)/(\Phi+1)$ maps $S_\alpha\times S_\alpha$ into $\D$ and vanishes at $(1,1)$. By Montel's theorem, every sequence in $\mathcal C$ has a subsequence along which these functions converge locally uniformly, and the limit maps into $\D$, by the maximum principle, since it vanishes at $(1,1)$; so the kernels converge locally uniformly. Their limit is positive, holomorphic in the first variable and antiholomorphic in the second, hence a holomorphic kernel (Section~\ref{ssec:kernels}); it takes the value $1$ at $(1,1)$, and, by \cite[Theorem~5.21]{PR} and Lemma~\ref{lem:cauchy}, it lies in $\mathcal C$.

We next find a kernel in $\mathcal C$ that is homogeneous under the dilations $z\mapsto e^tz$. For $t\in\RR$ and $k'\in\mathcal C$ put $(T_tk')(z,w)=k'(e^tz,e^tw)/k'(e^t,e^t)$, where $k'(e^t,e^t)>0$. Since $z\mapsto e^tz$ maps $S_\alpha$ onto itself, and $e^tA$ is admissible with $A$, \cite[Theorem~5.21]{PR} and Lemma~\ref{lem:covariance} show that $T_t$ maps $\mathcal C$ into itself, and $T_t$ is continuous. Moreover $T_s\circ T_t=T_{s+t}$ and $T_0$ is the identity. By the Schauder--Tychonoff fixed point theorem, each $T_t$ has a fixed point in $\mathcal C$, and its fixed point set is closed. Since $T_{2^{-n}}=T_{2^{-n-1}}\circ T_{2^{-n-1}}$, the fixed point set of $T_{2^{-n-1}}$ is contained in that of $T_{2^{-n}}$, so, by compactness, some $k_0\in\mathcal C$ is fixed by every $T_{2^{-n}}$, hence by every $T_t$ with $t$ a dyadic rational, as $T_{-t}=T_t^{-1}$. Both sides of $k_0(e^tz,e^tw)=k_0(e^t,e^t)\,k_0(z,w)$ are continuous in $t$, so this identity holds for all $t\in\RR$, and $t\mapsto k_0(e^t,e^t)$ is continuous and multiplicative; thus $k_0(e^tz,e^tw)=e^{2\eta t}k_0(z,w)$ for some $\eta\in\RR$. The kernel $k_0^\#(z,w)=\overline{k_0(\bar z,\bar w)}$ also lies in $\mathcal C$: $F\mapsto\overline{F(\bar z)}$ preserves $\sup_{S_\alpha}|F|$, and, by Lemmas~\ref{lem:conv}(a) and~\ref{lem:intkernel}, since $e^\#(A)=e(A^*)^*$ for $e^\#(z)=\overline{e(\bar z)}$, the moment map of $k_0^\#$ at $A$ is the adjoint, on $\Sc(\C^2)$, of that of $k_0$ at the admissible matrix $A^*$. So $k_1=\frac12(k_0+k_0^\#)\in\mathcal C$, $k_1^\#=k_1$ and $k_1(e^tz,e^tw)=e^{2\eta t}k_1(z,w)$. For $\zeta,v$ in the strip $\Pi_\alpha$, the function $H(\zeta,v)=e^{-\eta(\zeta+v)}k_1(e^\zeta,e^{\bar v})$ is holomorphic and satisfies $H(\zeta+t,v+t)=H(\zeta,v)$ for $t\in\RR$. For fixed $(\zeta,v)$, the function $s\mapsto H(\zeta+s,v+s)$ is holomorphic on the strip $\{s:\zeta+s,\,v+s\in\Pi_\alpha\}$, which contains $\RR$, and is constant on $\RR$, hence constant; taking $s=-(\zeta+v)/2$ gives $H(\zeta,v)=\kappa(\zeta-v)$, with $\kappa(u)=H(u/2,-u/2)$ holomorphic on $\Pi_{2\alpha}$. Thus
\[
k_1(z,w)=e^{\eta(\log z+\overline{\log w})}\,\kappa(\log z-\log\bar w)\qquad(z,w\in S_\alpha).
\]
The identity $k_1(w,z)=\overline{k_1(z,w)}$ gives $\kappa(-\bar u)=\overline{\kappa(u)}$, and $k_1^\#=k_1$ gives $\kappa(\bar u)=\overline{\kappa(u)}$; so $\kappa$ is even and real on $\RR$, and $\kappa(s)=e^{-\eta s}k_1(e^s,1)>0$ for $s\in\RR$, since $\Rp k_1>0$. Moreover $(m^2-F(z)\overline{F(w)})\,\kappa(\log z-\log\bar w)$ is a positive kernel whenever $\sup_{S_\alpha}|F|\le1$, as it is obtained from the corresponding kernel for $k_1$ by dividing by $e^{\eta\log z}\,\overline{e^{\eta\log w}}$; so multiplication by such $F$ has norm at most $m$ on the space $\E_\kappa$ of $\kappa(\log z-\log\bar w)$.

We now describe $\E_\kappa$ by the Fourier transform. For real $s,s'$, $\kappa(s-s')=e^{-\eta(s+s')}k_1(e^s,e^{s'})$, so $\kappa$ is positive definite on $\RR$, and, by Bochner's theorem, $\kappa(t)=\int e^{i\xi t}\,d\mu(\xi)$ for a finite even positive measure $\mu$. For $\tau\in\RR$, $F(z)=e^{-\alpha|\tau|}z^{i\tau}$ satisfies $\sup_{S_\alpha}|F|=1$, as $|z^{i\tau}|=e^{-\tau\arg z}$, and positivity of $(m^2-F(e^s)\overline{F(e^{s'})})\kappa(s-s')$ means, by the uniqueness in Bochner's theorem, that $\mu(E-\tau)\le m^2e^{2\alpha|\tau|}\mu(E)$ for Borel sets $E$. Averaging $\mu(E)\le m^2e^{2\alpha|\tau|}\mu(E+\tau)$ over $|\tau|\le1$ gives $\mu(E)\le\frac12m^2e^{2\alpha}\mu(\RR)|E|$, so $d\mu=\nu\,d\xi$ with $\nu$ even and bounded, and
\begin{equation}\label{eq:translate}
\nu(\xi)\le m^2e^{2\alpha|\tau|}\,\nu(\xi+\tau)\qquad(\xi,\tau\in\RR),
\end{equation}
almost everywhere for each $\tau$. Since $\kappa$ is holomorphic on $\Pi_{2\alpha}$, Fatou's lemma, applied to the central difference quotients, with $\Delta_h\kappa(x)=\kappa(x+h/2)-\kappa(x-h/2)$,
\[
(-1)^nh^{-2n}\Delta_h^{2n}\kappa(0)=\int\Big(\frac{2\sin(\xi h/2)}h\Big)^{2n}\nu(\xi)\,d\xi\qquad(h>0),
\]
gives, as $h\to0$, $\int\xi^{2n}\nu\,d\xi\le(-1)^n\kappa^{(2n)}(0)$, and the Cauchy estimates on the disk $|u|\le r<2\alpha$ bound this by $C_r(2n)!\,r^{-2n}$; hence $\int\cosh(b\xi)\,\nu(\xi)\,d\xi<\infty$ for $0\le b<2\alpha$, and, by analytic continuation,
\[
\kappa(u)=\int e^{i\xi u}\,\nu(\xi)\,d\xi\qquad(u\in\Pi_{2\alpha}).
\]
As $\kappa\ge0$ on $\RR$, Fubini's theorem gives $\int\kappa(t)e^{-\epsilon t^2}dt=\sqrt{\pi/\epsilon}\int\nu(\xi)e^{-\xi^2/(4\epsilon)}d\xi\le2\pi\sup\nu$ for $\epsilon>0$, so $\kappa\in L^1(\RR)$ by monotone convergence. By Fourier inversion we may take $\nu(\xi)=\frac1{2\pi}\int\kappa(t)e^{-i\xi t}\,dt$, which is continuous, with $\nu(0)>0$; then \eqref{eq:translate} holds everywhere, and $\nu>0$. By \eqref{eq:intspace}, with $\phi_\xi(z)=z^{i\xi}$ and $\beta^2=\nu^{-1}$, $\E_\kappa$ consists of the functions $f=\int\hat f(\xi)z^{i\xi}\,d\xi$ with
\begin{equation}\label{eq:Hkappa}
\norm f_{\E_\kappa}^2=\int|\hat f(\xi)|^2\,\nu(\xi)^{-1}\,d\xi<\infty;
\end{equation}
the required injectivity holds because such $\hat f$ are integrable, by the Cauchy--Schwarz inequality, as $\nu$ is, and $s\mapsto f(e^s)$ is the inverse Fourier transform of $\hat f$.

We now test the accretivity of the moment maps of $k_1$ on a family of $2\times2$ matrices. Let $t\in\RR$ and $0<\theta<\alpha$, and put $\lambda_1=e^{t/2+i\theta}$, $\lambda_2=e^{-t/2-i\theta}$, $b=2(\sin^2\alpha-\sin^2\theta)^{1/2}$, $c=b/|\lambda_1-\lambda_2|$ and $A=\big(\begin{smallmatrix}\lambda_1&c(\lambda_1-\lambda_2)\\0&\lambda_2\end{smallmatrix}\big)$. With $\gamma=\pi/2-\alpha$, the matrices $\Rp(e^{\pm i\gamma}A)$ have the positive diagonal entries $e^{t/2}\sin(\alpha\mp\theta)$ and $e^{-t/2}\sin(\alpha\pm\theta)$, off-diagonal entries of modulus $b/2$ and determinant $\sin(\alpha-\theta)\sin(\alpha+\theta)-b^2/4=0$, so they are positive; thus $W(A)\subseteq\Sbar\alpha$, and $A$ is admissible. Moreover $A=\lambda_1P_1+\lambda_2P_2$ with $P_1=\big(\begin{smallmatrix}1&c\\0&0\end{smallmatrix}\big)$ and $P_2=I-P_1$, and, by Lemma~\ref{lem:conv}(a), $\Psi_A(Y)=\sum_{r,s=1}^2k_1(\lambda_r,\lambda_s)P_rYP_s^*$, where
\begin{gather*}
k_1(\lambda_1,\lambda_1)=e^{\eta t}\kappa(2i\theta),\qquad k_1(\lambda_2,\lambda_2)=e^{-\eta t}\kappa(2i\theta),\\
k_1(\lambda_1,\lambda_2)=\overline{k_1(\lambda_2,\lambda_1)}=e^{2i\eta\theta}\kappa(t).
\end{gather*}
For $X=\big(\begin{smallmatrix}-c&1\\1&c\end{smallmatrix}\big)$, $\ip{X}{P_rXP_r^*}_{\mathrm{HS}}=-c^2(1+c^2)$ for $r=1,2$, and $\ip{X}{P_1XP_2^*}_{\mathrm{HS}}=\ip{X}{P_2XP_1^*}_{\mathrm{HS}}=(1+c^2)^2$, so
\[
\Rp\ip X{\Psi_A(X)}_{\mathrm{HS}}=2(1+c^2)\big((1+c^2)\cos(2\eta\theta)\,\kappa(t)-c^2\cosh(\eta t)\,\kappa(2i\theta)\big).
\]
Since $|\lambda_1-\lambda_2|^2=2(\cosh t-\cos2\theta)$ and $b^2=2(\cos2\theta-\cos2\alpha)$, we have $1+c^2=(\cosh t-\cos2\alpha)/(\cosh t-\cos2\theta)$, and the accretivity of $\Psi_A$ gives
\begin{equation}\label{eq:twobytwo}
\cos(2\eta\theta)\,\kappa(t)\,(\cosh t-\cos2\alpha)\ge\cosh(\eta t)\,\kappa(2i\theta)\,(\cos2\theta-\cos2\alpha).
\end{equation}
For $\tilde k_{S_\alpha}$, where $\eta=0$ and $\kappa(u)=\sin(2\alpha)/(\cosh u-\cos2\alpha)$, both sides of \eqref{eq:twobytwo} equal $\sin2\alpha$.

Letting $\theta\uparrow\alpha$ in \eqref{eq:twobytwo} bounds $\kappa$ from below. Let $g(\xi)=e^{2\alpha\xi}\nu(\xi)$ and $g_\infty=\liminf_{\xi\to\infty}g(\xi)$. By \eqref{eq:translate} with $\xi=0$, $g(\tau)\ge\nu(0)/m^2$ for $\tau\ge0$, so $g_\infty>0$. For $0<L<g_\infty$ choose $R$ with $g\ge L$ on $[R,\infty)$. Since $\nu$ is even, for $0<\theta<\alpha$,
\[
\kappa(2i\theta)=\int\nu(\xi)e^{2\theta\xi}\,d\xi\ge\int_R^\infty g(\xi)e^{-2(\alpha-\theta)\xi}\,d\xi\ge\frac{L\,e^{-2(\alpha-\theta)R}}{2(\alpha-\theta)},
\]
and $(\cos2\theta-\cos2\alpha)/(2(\alpha-\theta))\to\sin2\alpha$ as $\theta\uparrow\alpha$. So \eqref{eq:twobytwo} gives, in the limit,
\[
\cos(2\eta\alpha)\,\kappa(t)\,(\cosh t-\cos2\alpha)\ge L\cosh(\eta t)\sin2\alpha\qquad(t\in\RR).
\]
At $t=0$ this shows that $g_\infty<\infty$ and $\cos(2\eta\alpha)>0$, and we may let $L\uparrow g_\infty$. As $\kappa>0$ and $\cos(2\eta\alpha)\le1\le\cosh(\eta t)$, it follows that $\kappa\ge g_\infty\tilde\kappa$ on $\RR$, where $\tilde\kappa(t)=\sin(2\alpha)/(\cosh t-\cos2\alpha)=\int e^{i\xi t}\nu_\alpha(\xi)\,d\xi$ (proof of Lemma~\ref{lem:prodweights}). Integrating, by Fourier inversion,
\begin{equation}\label{eq:nuzero}
2\pi\nu(0)=\int\kappa\,dt\ge g_\infty\int\tilde\kappa\,dt=2\pi g_\infty\nu_\alpha(0)=2\pi g_\infty\Big(1-\frac{2\alpha}\pi\Big).
\end{equation}

Finally, we test multipliers. For $\tau>0$, $F_\tau(z)=\cos(\tau\log z)/\cosh(\alpha\tau)=(z^{i\tau}+z^{-i\tau})/(2\cosh\alpha\tau)$ satisfies $\sup_{S_\alpha}|F_\tau|\le1$, since $|\cos(x+iy)|\le\cosh y$. For $0<\delta<\tau$ let $f_\delta=(2\delta)^{-1/2}\int_{-\delta}^\delta z^{i\xi}\,d\xi\in\E_\kappa$. Then $F_\tau f_\delta=\int\hat h(\xi)z^{i\xi}\,d\xi$ with $\hat h=(2\delta)^{-1/2}(2\cosh\alpha\tau)^{-1}\big(1_{[\tau-\delta,\tau+\delta]}+1_{[-\tau-\delta,-\tau+\delta]}\big)$, the two intervals being disjoint, so, by \eqref{eq:Hkappa} and the continuity and evenness of $\nu$, as $\delta\to0$,
\[
\norm{f_\delta}^2\to\frac1{\nu(0)},\qquad\norm{F_\tau f_\delta}^2\to\frac1{2\cosh^2(\alpha\tau)\,\nu(\tau)}=\frac2{g(\tau)\,(1+e^{-2\alpha\tau})^2}.
\]
Since multiplication by $F_\tau$ has norm at most $m$ on $\E_\kappa$, this gives $m^2\ge2\nu(0)/\big(g(\tau)(1+e^{-2\alpha\tau})^2\big)$. Letting $\tau\to\infty$ along a sequence with $g(\tau)\to g_\infty$, and using \eqref{eq:nuzero}, gives $m^2\ge2\nu(0)/g_\infty\ge2-4\alpha/\pi$.
\end{proof}

\par\textit{Acknowledgements}. This work was supported in part by a grant of access to OpenAI models through the ChatGPT for Academic Researchers program. In particular, ChatGPT was used as an auxiliary editorial tool for proofreading and preliminary literature search assistance. We have also benefited from proof checking, stylistic feedback and bibliographical guidance offered by Anthropic's Claude models. The authors verified all arguments and take full responsibility for the content of the paper.

\defbibheading{plain}[\refname]{\section*{#1}}
\printbibliography[heading=plain]

\end{document}